\documentclass{amsart}

\usepackage{amsfonts,amssymb,verbatim,amsmath,amsthm,latexsym,textcomp,amscd}
\usepackage{latexsym,amsfonts,amssymb,epsfig,verbatim}
\usepackage{amsmath,amsthm,amssymb,latexsym,graphics,textcomp}
\usepackage{graphicx}
\usepackage{color}
\usepackage{url}

\input xy
\xyoption{all}

\begin{document}

\newtheorem{theorem}{Theorem}[section]
\newtheorem{prop}[theorem]{Proposition}
\newtheorem{lemma}[theorem]{Lemma}
\newtheorem{cor}[theorem]{Corollary}
\newtheorem{definition}[theorem]{Definition}
\newtheorem{conj}[theorem]{Conjecture}
\newtheorem{rmk}[theorem]{Remark}
\newtheorem{claim}[theorem]{Claim}
\newtheorem{defth}[theorem]{Definition-Theorem}

\newcommand{\boundary}{\partial}
\newcommand{\C}{{\mathbb C}}
\newcommand{\integers}{{\mathbb Z}}
\newcommand{\natls}{{\mathbb N}}
\newcommand{\ratls}{{\mathbb Q}}
\newcommand{\bbR}{{\mathbb R}}
\newcommand{\proj}{{\mathbb P}}
\newcommand{\lhp}{{\mathbb L}}
\newcommand{\tube}{{\mathbb T}}
\newcommand{\cusp}{{\mathbb P}}
\newcommand\AAA{{\mathcal A}}
\newcommand\BB{{\mathcal B}}
\newcommand\CC{{\mathcal C}}
\newcommand\DD{{\mathcal D}}
\newcommand\EE{{\mathcal E}}
\newcommand\FF{{\mathcal F}}
\newcommand\GG{{\mathcal G}}
\newcommand\HH{{\mathcal H}}
\newcommand\II{{\mathcal I}}
\newcommand\JJ{{\mathcal J}}
\newcommand\KK{{\mathcal K}}
\newcommand\LL{{\mathcal L}}
\newcommand\MM{{\mathcal M}}
\newcommand\NN{{\mathcal N}}
\newcommand\OO{{\mathcal O}}
\newcommand\PP{{\mathcal P}}
\newcommand\QQ{{\mathcal Q}}
\newcommand\RR{{\mathcal R}}
\newcommand\SSS{{\mathcal S}}
\newcommand\TT{{\mathcal T}}
\newcommand\UU{{\mathcal U}}
\newcommand\VV{{\mathcal V}}
\newcommand\WW{{\mathcal W}}
\newcommand\XX{{\mathcal X}}
\newcommand\YY{{\mathcal Y}}
\newcommand\ZZ{{\mathcal Z}}
\newcommand\CH{{\CC\HH}}
\newcommand\PEY{{\PP\EE\YY}}
\newcommand\MF{{\MM\FF}}
\newcommand\RCT{{{\mathcal R}_{CT}}}
\newcommand\PMF{{\PP\kern-2pt\MM\FF}}
\newcommand\FL{{\FF\LL}}
\newcommand\PML{{\PP\kern-2pt\MM\LL}}
\newcommand\GL{{\GG\LL}}
\newcommand\Pol{{\mathcal P}}
\newcommand\half{{\textstyle{\frac12}}}
\newcommand\Half{{\frac12}}
\newcommand\Mod{\operatorname{Mod}}
\newcommand\Area{\operatorname{Area}}
\newcommand\ep{\epsilon}
\newcommand\hhat{\widehat}
\newcommand\Proj{{\mathbf P}}
\newcommand\U{{\mathbf U}}
 \newcommand\Hyp{{\mathbf H}}
\newcommand\D{{\mathbf D}}
\newcommand\Z{{\mathbb Z}}
\newcommand\R{{\mathbb R}}
\newcommand\Q{{\mathbb Q}}
\newcommand\E{{\mathbb E}}
\newcommand\til{\widetilde}
\newcommand\length{\operatorname{length}}
\newcommand\tr{\operatorname{tr}}
\newcommand\gesim{\succ}
\newcommand\lesim{\prec}
\newcommand\simle{\lesim}
\newcommand\simge{\gesim}
\newcommand{\simmult}{\asymp}
\newcommand{\simadd}{\mathrel{\overset{\text{\tiny $+$}}{\sim}}}
\newcommand{\ssm}{\setminus}
\newcommand{\diam}{\operatorname{diam}}
\newcommand{\pair}[1]{\langle #1\rangle}
\newcommand{\T}{{\mathbf T}}
\newcommand{\inj}{\operatorname{inj}}
\newcommand{\pleat}{\operatorname{\mathbf{pleat}}}
\newcommand{\short}{\operatorname{\mathbf{short}}}
\newcommand{\vertices}{\operatorname{vert}}
\newcommand{\collar}{\operatorname{\mathbf{collar}}}
\newcommand{\bcollar}{\operatorname{\overline{\mathbf{collar}}}}
\newcommand{\I}{{\mathbf I}}
\newcommand{\tprec}{\prec_t}
\newcommand{\fprec}{\prec_f}
\newcommand{\bprec}{\prec_b}
\newcommand{\pprec}{\prec_p}
\newcommand{\ppreceq}{\preceq_p}
\newcommand{\sprec}{\prec_s}
\newcommand{\cpreceq}{\preceq_c}
\newcommand{\cprec}{\prec_c}
\newcommand{\topprec}{\prec_{\rm top}}
\newcommand{\Topprec}{\prec_{\rm TOP}}
\newcommand{\fsub}{\mathrel{\scriptstyle\searrow}}
\newcommand{\bsub}{\mathrel{\scriptstyle\swarrow}}
\newcommand{\fsubd}{\mathrel{{\scriptstyle\searrow}\kern-1ex^d\kern0.5ex}}
\newcommand{\bsubd}{\mathrel{{\scriptstyle\swarrow}\kern-1.6ex^d\kern0.8ex}}
\newcommand{\fsubeq}{\mathrel{\raise-.7ex\hbox{$\overset{\searrow}{=}$}}}
\newcommand{\bsubeq}{\mathrel{\raise-.7ex\hbox{$\overset{\swarrow}{=}$}}}
\newcommand{\tw}{\operatorname{tw}}
\newcommand{\base}{\operatorname{base}}
\newcommand{\trans}{\operatorname{trans}}
\newcommand{\rest}{|_}
\newcommand{\bbar}{\overline}
\newcommand{\UML}{\operatorname{\UU\MM\LL}}
\newcommand{\EL}{\mathcal{EL}}
\newcommand{\tsum}{\sideset{}{'}\sum}
\newcommand{\tsh}[1]{\left\{\kern-.9ex\left\{#1\right\}\kern-.9ex\right\}}
\newcommand{\Tsh}[2]{\tsh{#2}_{#1}}
\newcommand{\qeq}{\mathrel{\approx}}
\newcommand{\Qeq}[1]{\mathrel{\approx_{#1}}}
\newcommand{\qle}{\lesssim}
\newcommand{\Qle}[1]{\mathrel{\lesssim_{#1}}}
\newcommand{\simp}{\operatorname{simp}}
\newcommand{\vsucc}{\operatorname{succ}}
\newcommand{\vpred}{\operatorname{pred}}
\newcommand\fhalf[1]{\overrightarrow {#1}}
\newcommand\bhalf[1]{\overleftarrow {#1}}
\newcommand\sleft{_{\text{left}}}
\newcommand\sright{_{\text{right}}}
\newcommand\sbtop{_{\text{top}}}
\newcommand\sbot{_{\text{bot}}}
\newcommand\sll{_{\mathbf l}}
\newcommand\srr{_{\mathbf r}}
\newcommand\geod{\operatorname{\mathbf g}}
\newcommand\mtorus[1]{\boundary U(#1)}
\newcommand\A{\mathbf A}
\newcommand\Aleft[1]{\A\sleft(#1)}
\newcommand\Aright[1]{\A\sright(#1)}
\newcommand\Atop[1]{\A\sbtop(#1)}
\newcommand\Abot[1]{\A\sbot(#1)}
\newcommand\boundvert{{\boundary_{||}}}
\newcommand\storus[1]{U(#1)}
\newcommand\Momega{\omega_M}
\newcommand\nomega{\omega_\nu}
\newcommand\twist{\operatorname{tw}}
\newcommand\modl{M_\nu}
\newcommand\MT{{\mathbb T}}
\newcommand\Teich{{\mathcal T}}
\renewcommand{\Re}{\operatorname{Re}}
\renewcommand{\Im}{\operatorname{Im}}

\title{Quasi-metric antipodal spaces and maximal Gromov hyperbolic spaces}

\author{Kingshook Biswas}
\address{Indian Statistical Institute, Kolkata, India. Email: kingshook@isical.ac.in}

\begin{abstract} Hyperbolic fillings of metric spaces are a well-known tool for proving results on extending  
quasi-Moebius maps between boundaries of Gromov hyperbolic spaces to quasi-isometries between the spaces. 
For a hyperbolic filling $Y$ of the boundary of a Gromov hyperbolic space $X$, one has a quasi-Moebius 
identification between the boundaries $\partial Y$ and $\partial X$. 

\medskip

For CAT(-1) spaces, and more generally boundary continuous Gromov hyperbolic spaces, 
one can refine the quasi-Moebius structure on the boundary to a Moebius structure. It is 
then natural to ask whether there exists a functorial hyperbolic filling of the boundary by a boundary continuous 
Gromov hyperbolic space with an identification between 
boundaries which is not just quasi-Moebius, but in fact Moebius. The filling should be functorial in the sense 
that a Moebius homeomorphism between boundaries should induce an isometry between there fillings.

\medskip

We give a positive answer to this question for a large class of boundaries satisfying one crucial hypothesis, 
the {\it antipodal property}. This gives a class of compact spaces called {\it quasi-metric antipodal spaces}. 
For any such space $Z$, we give a functorial construction of a 
boundary continuous Gromov hyperbolic space $\mathcal{M}(Z)$ together with a Moebius identification of its 
boundary with $Z$. The space $\mathcal{M}(Z)$ is maximal amongst all fillings of $Z$. These spaces $\mathcal{M}(Z)$ give 
in fact all examples of a natural class of spaces called {\it maximal Gromov hyperbolic spaces}. 

\medskip

We prove an equivalence of categories between quasi-metric antipodal spaces and maximal Gromov hyperbolic spaces. 
This is part of a more general equivalence we prove between the larger categories of certain spaces called {\it antipodal spaces} 
and {\it maximal Gromov product spaces}. We prove that the injective hull of a Gromov product space $X$ is isometric to 
the maximal Gromov product space $\mathcal{M}(Z)$, where $Z$ is the boundary of $X$. We also show that a Gromov product space 
is injective if and only if it is maximal. 

\end{abstract}

\bigskip

\maketitle

\tableofcontents

\section{Introduction.}

\medskip

The Gromov inner product $(.|.)_x$ on the visual boundary $\partial X$ of a Gromov hyperbolic space $X$ gives a 
family of quasi-metrics $\rho_x = e^{-(.|)_x}$ on the boundary, which we call {\it visual quasi-metrics}. The visual 
quasi-metrics define cross-ratios $[.,.,.,.]_{\rho_x}$ on the boundary which differ from each other by a uniformly bounded 
multiplicative error, and hence define a canonical {\it quasi-Moebius structure} on the boundary. The correspondence 
between quasi-isometries of Gromov hyperbolic spaces and quasi-Moebius maps of their boundaries is well studied. 

\medskip

It is well known that a quasi-isometry $F : X \to Y$ between Gromov hyperbolic spaces extends to a boundary map 
$\partial F : \partial X \to \partial Y$ which is quasisymmetric, more specifically it is {\it power quasisymmetric}. 
Conversely, Bonk and Schramm proved that any power quasisymmetric map between boundaries of Gromov hyperbolic spaces (satisfying 
some mild hypotheses) extends to a quasi-isometry between the spaces (\cite{bonk-schramm}). They also 
proved that a bi-Lipschitz map between boundaries extends to an {\it almost-isometric} map between the spaces (with the same 
hypotheses on the spaces). Under the stronger hypothesis of uniformly perfect and complete boundaries, Buyalo-Schroeder  
showed that any quasisymmetric map, or more generally quasi-Moebius map, between boundaries extends to a quasi-isometry between the Gromov hyperbolic spaces \cite{buyalo-schroeder}. These results may all be obtained from general results of Jordi \cite{jordi}, who proved that a power-quasi-Moebius map 
between boundaries of visual, roughly geodesic, Gromov hyperbolic spaces extends to a quasi-isometry between the spaces.

\medskip

All of the above extension results above depend on the idea of constructing in a (more or less) functorial way a {\it hyperbolic filling} 
of a given metric (or more generally quasi-metric) space $Z$. A hyperbolic filling of $Z$ is a Gromov hyperbolic space $X$ equipped with 
a homeomorphism $f : \partial X \to Z$ from its boundary to $Z$, which one normally requires to be bi-Lipschitz with respect to some 
visual metric (or visual quasi-metric) on $\partial X$. In particular, the identification of $\partial X$ with $Z$ is quasi-Moebius. 

\medskip

In the special case of CAT(-1) spaces however, there is a finer structure on the boundary. In this case, the Gromov inner product extends continuously 
to the boundary. This has the consequence that the cross-ratios of all the visual quasi-metrics $\rho_x$ (which are metrics in this case) agree with no 
error, $[.,.,.,.]_{\rho_x} = [.,.,.,.]_{\rho_y}$ for all $x,y \in X$. There is a thus a canonical cross-ratio and a hence a canonical {\it Moebius structure} 
on the boundary of a CAT(-1) space, as opposed to the rougher quasi-Moebius structure on the boundary of general Gromov hyperbolic spaces. More 
generally, this Moebius structure on the boundary also makes sense for {\it boundary continuous} Gromov hyperbolic spaces, which are Gromov 
hyperbolic spaces where the Gromov inner product extends continuously to the boundary. As we show in section 5, proper, CAT(0) Gromov hyperbolic spaces are boundary continuous, giving a large class of examples of such spaces.   

\medskip

In this setting, a natural question is whether one can improve on the hyperbolic filling constructions of \cite{bonk-schramm} and \cite{buyalo-schroeder}, 
in such a way that the identification of the boundary $\partial X$ of the filling with the given quasi-metric space $Z$ is not just quasi-Moebius, 
but in fact Moebius. We then say that the filling is a {\it Moebius} hyperbolic filling. One could even require the stronger condition that the identification $f : (\partial X, \rho_x) \to (Z, \rho)$ is
an isometry of quasi-metric spaces. We remark that this problem of finding a Moebius hyperbolic filling is quite non-trivial; the hyperbolic fillings of \cite{bonk-schramm} and \cite{buyalo-schroeder} are too rough as constructions to recover anything more than the coarse quasi-Moebius structure, and cannot recover the finer Moebius structure.

\medskip

For CAT(-1) spaces, a hyperbolic filling of their boundaries was constructed in \cite{biswas3}, where it was used to show that a 
Moebius homeomorphism between the boundaries of proper, geodesically complete CAT(-1) spaces extends to a $(1, \log 2)$-quasi-isometry 
between the spaces.  More precisely, given a proper, geodesically complete CAT(-1) space $X$, in \cite{biswas3} one defined the space 
$\mathcal{M}(\partial X)$ as the space of diameter one antipodal Moebius metrics on the boundary $\partial X$, equipped with a 
certain metric similar to the Thurston metric on Teichmuller space. Here a metric $\rho$ on $\partial X$ is said to be Moebius if it has the 
same cross-ratio as the visual metrics, and if it has diameter one then it is said to be {\it antipodal} if every $\xi \in \partial X$ 
has an "antipode", i.e. a point $\eta \in \partial X$ such that $\rho(\xi, \eta) = 1$. The space $\mathcal{M}(\partial X)$ comes equipped 
with a natural map $i_X : X \to \mathcal{M}(\partial X), x \mapsto \rho_x$; in \cite{biswas3} this map is shown to be an 
isometric embedding with image $\frac{1}{2}\log 2$-dense in $\mathcal{M}(\partial X)$. The functoriality of the filling $\mathcal{M}(Z)$ 
with respect to Moebius maps is easily seen. If $f : \partial X \to \partial Y$ is a Moebius homeomorphism between boundaries
of proper, geodesically complete CAT(-1) spaces $X, Y$, then one can push-forward by $f$ antipodal functions on $\partial X$ to antipodal functions on 
$\partial Y$, and since $f$ is Moebius it takes Moebius antipodal functions to Moebius antipodal functions. Thus $f$ canonically induces a map 
$f_* : \mathcal{M}(\partial X) \to \mathcal{M}(\partial Y)$, which it is not hard to show is an isometry.   

\medskip

It turns out that the crucial property required to solve the problem of finding  Moebius hyperbolic fillings (as opposed to fillings which are just quasi-Moebius) is the antipodal property. 
In the present paper, we generalize the construction of \cite{biswas3} to a much larger class of boundaries, which we call
{\it antipodal spaces}. Given a compact metrizable space $Z$, an {\it antipodal function} on $Z$ is a continuous function 
$\rho : Z \times Z \to [0,1]$ which satisfies symmetry, positivity, and the antipodal property: every $\xi \in Z$ has an antipode 
$\eta \in Z$, i.e. $\eta$ satisfies $\rho(\xi, \eta) = 1$. An antipodal space is then defined to be a compact Hausdorff space $Z$ equipped with an 
antipodal function $\rho$. The motivating example of an antipodal space is the boundary of a proper, geodesically complete, 
boundary continuous Gromov hyperbolic space equipped with a visual quasi-metric $\rho_x$. In this case, it is the 
geodesic completeness of the Gromov hyperbolic space which implies that the visual quasi-metrics are antipodal.    

\medskip

Given an antipodal space $(Z, \rho)$, we then define the {\it Moebius space} $\mathcal{M}(Z)$ to be the set of Moebius antipodal functions on $Z$, 
i.e. antipodal functions with the same cross-ratio as the that of the given antipodal function on $Z$. Note the difference 
from the construction in \cite{biswas3} is that we do not require the antipodal functions to be metrics, i.e. they do not have to 
satisfy the triangle inequality (in fact the space constructed in \cite{biswas3} is a closed subspace of the space constructed here). 

\medskip

It is quite non-trivial to construct points and geodesics in the Moebius space $\mathcal{M}(Z)$. We use an ODE on the Banach space $C(Z)$ 
(continuous functions on $Z$ with the sup norm), which we call the {\it antipodal flow}, in order to construct points and geodesics in $\mathcal{M}(Z)$. 
We show that any solution of the antipodal flow converges uniformly as time tends to infinity, and Moebius antipodal functions on $Z$ are 
constructed from the limits of solutions. This also allows us to construct geodesics in $\mathcal{M}(Z)$, and to show that $\mathcal{M}(Z)$ is 
contractible. We have:

\medskip

\begin{theorem} \label{mainthm1} Let $Z$ be an antipodal space. Then the metric space $\mathcal{M}(Z)$ is
unbounded, contractible, proper, geodesic and geodesically complete.
\end{theorem}

\medskip

If the antipodal space $(Z, \alpha)$ is quasi-metric, i.e. $\alpha$ is a quasi-metric, then we obtain a Moebius hyperbolic filling of $Z$ by $\mathcal{M}(Z)$:

\medskip

\begin{theorem} \label{mainthm2} Let $Z$ be an antipodal space. 

\medskip

Then $Z$ is quasi-metric if and only if $\mathcal{M}(Z)$ is Gromov hyperbolic.

\medskip

 In this case the space $\mathcal{M}(Z)$ is a boundary continuous Gromov hyperbolic space which is a Moebius hyperbolic filling of $Z$.
 
 \medskip
  
In fact, for the basepoint $\alpha \in \mathcal{M}(Z)$, the identification $f : (\partial\mathcal{M}(Z), \rho_{\alpha}) \to (Z, \alpha)$ 
is an isometry of quasi-metric spaces. 


\end{theorem}

\medskip

We remark that for $Z$ a quasi-metric antipodal space, the Moebius hyperbolic filling $\MM(Z)$ has the desirable property that it is contractible, as opposed to the hyperbolic filling of \cite{buyalo-schroeder}, which is not necessarily contractible.

\medskip

Theorem \ref{mainthm2}, together with further results below, generalize one of the main results of Beyrer-Schroeder \cite{beyrer-schroeder} (see also Hughes \cite{hughes}), who consider ultrametric antipodal spaces, 
to the much larger class of quasi-metric antipodal spaces. Beyrer-Schroeder construct for any ultrametric antipodal space $Z$ a geodesically complete tree $X$ which is a 
Moebius hyperbolic filling of $Z$ (we remark that their construction is done not only for compact ultrametric spaces, but also for possibly 
non-compact but still complete ultrametric spaces). Using this, they prove an equivalence of categories between trees and ultrametric spaces, which 
implies that a proper, geodesically complete tree is determined by the Moebius structure on its boundary as an ultrametric antipodal space. We show that for an ultrametric antipodal space $Z$, the space $\mathcal{M}(Z)$ is a tree 
(which is necessarily isometric to the tree constructed by Beyrer-Schroeder).    

\medskip

In our case, we extend the equivalence between ultrametric antipodal spaces and trees to 
an equivalence between the larger categories of quasi-metric antipodal spaces and a certain class of Gromov hyperbolic spaces which we call {\it maximal Gromov hyperbolic spaces}. Roughly speaking, a Gromov hyperbolic space $X$ (which is assumed to be proper, geodesically complete, and boundary continuous) is {\it maximal} if it is the "largest" Moebius hyperbolic filling of its boundary $\partial X$.
More precisely, for any Moebius hyperbolic filling $(Y, f: \partial Y \to \partial X)$ of $\partial X$ (where $Y$ is a proper, geodesically complete, 
boundary continuous Gromov hyperbolic space), the Moebius map $f : \partial Y \to \partial X$ extends to an isometric embedding $F : Y \to X$. 

\medskip

It is natural to ask for examples of maximal Gromov hyperbolic spaces. We show: 

\medskip

\begin{theorem} \label{mainthm3} Let $Z$ be a quasi-metric antipodal space. 
Then the Moebius space $\MM(Z)$ is a maximal Gromov hyperbolic space.
\end{theorem}

\medskip

In particular, any proper, geodesically complete trees is a maximal Gromov hyperbolic space, as we show in section 8 that it is of the form $\mathcal{M}(Z)$ for $Z$ an ultrametric antipodal space.
The Moebius spaces $\mathcal{M}(Z)$ for $Z$ quasi-metric in fact give all examples of maximal spaces. We note that just as in the case of CAT(-1) spaces, for a boundary continuous Gromov hyperbolic space $X$ (which is also 
assumed to be proper and geodesically complete), there is a natural map $i_X : X \to \mathcal{M}(\partial X), x \mapsto \rho_x$. We show that as before, 
this map is an isometric embedding, which we call the {\it visual embedding} of $X$. We then have:

\medskip

\begin{theorem} \label{mainthm4} Let $X$ be a proper, geodesically complete, boundary continuous Gromov hyperbolic space. Let $\delta \geq 0$ be such that
all geodesic triangles in $X$ are $\delta$-thin. 

\medskip

Then the image of the visual
embedding $i_X : X \to \mathcal{M}(\partial X)$ is $7\delta$-dense in $\mathcal{M}(\partial X)$.

\medskip

The visual embedding $i_X$ is an isometry if and only if $X$ is a maximal Gromov hyperbolic space.
\end{theorem}
  
\medskip

Thus any maximal Gromov hyperbolic space $X$ is isometric to the 
Moebius space $\MM(\partial X)$. The functoriality of the assignment 
$Z \mapsto \MM(Z)$ then implies the following, which completes the equivalence 
between quasi-metric antipodal spaces and maximal Gromov hyperbolic spaces:

\medskip

\begin{theorem} \label{mainthm5} Let $X, Y$ be maximal Gromov hyperbolic spaces, 
amd let $f : \partial X \to \partial Y$ be a Moebius homeomorphism. Then $f$ 
extends to an isometry $F : X \to Y$.
\end{theorem}

\medskip

The formalism of a cross-ratio on the boundary 
derived from Gromov inner products makes sense also for a larger class 
of spaces which are not necessarily Gromov hyperbolic, and which includes 
the Moebius spaces $\MM(Z)$ for any antipodal space $Z$ (not necessarily quasi-metric). 

\medskip

Given a proper, geodesic and geodesically complete metric space $X$, we say that $X$ is a {\it Gromov product space} if there exists a second countable, Hausdorff
compactification $\widehat{X}$ of $X$ such that for any $x \in X$, the Gromov inner product $(.|.)_x : X \times X \to [0, +\infty)$ extends to a
continuous function $(.|.)_x : \widehat{X} \times \widehat{X} \to [0, +\infty]$ such that $(\xi|\eta)_x = +\infty$ if and only if $\xi = \eta \in \widehat{X} \setminus X$.
It is not hard to show that such a compactification, if it exists, is unique up to equivalence, and we call it the {\it Gromov product compactification}
of the Gromov product space $X$, and denote it by $\widehat{X}^{P}$. The compact metrizable space $\partial_{P}X := \widehat{X}^{P} \setminus X$ is called
the {\it Gromov product boundary} of $X$. The motivating example of Gromov product 
spaces is given by proper, geodesically complete, boundary continuous Gromov 
hyperbolic spaces, where the Gromov product boundary is given by the visual boundary. 

\medskip

The functions $\rho_x = e^{-(.|.)_x}, x \in X$ on the Gromov product 
boundary are then antipodal functions with the same cross-ratio, so the 
Gromov product boundary $\partial_{P} X$ is an antipodal space. We can then define a {\it filling} of an antipodal space $Z$ to be a Gromov product space $X$ 
equipped with a Moebius homeomorphism $f : \partial_{P} X \to Z$. As before, we define a Gromov product space $X$ to be {\it maximal} if it is the "largest" 
filling of its boundary $\partial_{P} X$: for any filling $(Y, f : \partial_{P} Y \to \partial_{P} X)$, the map $f$ extends to an isometric embedding $F : Y \to X$. 
We then have:

\medskip

\begin{theorem} \label{mainthm6} Let $Z$ be an antipodal space. Then the 
Moebius space $\MM(Z)$ is a maximal Gromov product space which is a filling of $Z$.
\end{theorem}

\medskip

In particular, the spaces $\MM(Z)$, for $Z$ antipodal but not quasi-metric, 
give examples of Gromov product spaces which are not Gromov hyperbolic.

\medskip

As before, the visual embedding $i_X : X \to \MM(\partial_{P} X), x \mapsto \rho_x$ makes 
sense for Gromov product spaces, and we show that it is an isometric embedding. Moreover:

\medskip

\begin{theorem} \label{mainthm7} Let $X$ be a Gromov product space. Then $X$ is maximal if and only if the visual embedding $i_X : X \to \MM(\partial_{P} X)$ is 
an isometry.
\end{theorem}

\medskip

Thus all maximal Gromov product spaces are of the form $\mathcal{M}(Z)$ for some 
antipodal space $Z$. We then obtain an equivalence between antipodal spaces and maximal Gromov product spaces, generalizing that between quasi-metric antipodal spaces and maximal Gromov hyperbolic spaces.

\medskip

\begin{theorem} \label{mainthm8} Let $X, Y$ be maximal Gromov product spaces and let $f : \partial_{P} X \to \partial_{P} Y$ be a Moebius homeomorphism. Then $f$ extends to an isometry $F : X \to Y$.
\end{theorem}

\medskip

There is also a connection with injective metric spaces and injective hulls. We recall that a metric space $X$
is said to be injective if for any metric space $B$, if $f : A \subset B \to X$ is a $1$-Lipschitz map from a subset $A$ of $B$ to $X$,
then $f$ extends to a $1$-Lipschitz map $\overline{f} : B \to X$. Injective metric spaces have many features in common with CAT(0) spaces,
for example they are geodesic, contractible, and admit a convex geodesic bi-combing. Examples include metric trees and $l_{\infty}(I)$
for any index set $I$. For an introduction to basic properties and the literature
on injective metric spaces, we refer to Lang \cite{langinjective}. 

\medskip

For any metric space $X$, there is a "smallest" injective metric space $E(X)$ into which $X$ embeds isometrically, called the {\it injective hull} of $X$. 
More precisely, there is an isometric embedding $e : X \to E(X)$ such that any isometric embedding of $X$ into an injective metric space factors through $e$. 
On the other hand, if $X$ is a Gromov product space, then the visual embedding gives an isometric embedding of $X$ into a Gromov product space $\MM(\partial_{P} X)$ which is the "largest" filling of the boundary $\partial_{P}$ of $X$. We show:

\medskip

\begin{theorem} \label{mainthm9} Let $X$ be a Gromov product space. 

\medskip

\noindent Then the injective hull $E(X)$ of $X$ is isometric to the Moebius space $\MM(\partial_{P} X)$.
\end{theorem}

\medskip

As an immediate corollary we obtain:

\begin{theorem}\label{mainthm10} Let $X$ be a Gromov product space. 

\medskip

\noindent Then $X$ is maximal if and only if $X$ is injective.

\medskip

\noindent In particular, for any antipodal space $Z$, the Moebius space $\mathcal{M}(Z)$ is injective.
\end{theorem}

\medskip

We also study the smooth structure of the Moebius spaces $\MM(Z)$. 
This arises from a family of canonical isometric embeddings of $\mathcal{M}(Z)$ into the Banach
space $C(Z)$ of continuous functions on $Z$ equipped with the sup norm. 
Any two such embeddings differ by a translation of $C(Z)$.
Thus we can talk of smoothness of maps $f : M \to \mathcal{M}(Z)$, where $M$ is a smooth manifold (such a map is defined to be
smooth if the map $f$ composed with one of these isometric embeddings is a smooth map into $C(Z)$), in particular one can talk
of smooth curves in $\mathcal{M}(Z)$. Motivated by this, for an antipodal space $Z$ and $\rho \in \mathcal{M}(Z)$, we say
that a curve $\gamma : t \in (-\epsilon, \epsilon) \mapsto \rho(t) \in \mathcal{M}(Z)$ is
{\it admissible at $\rho$} if $\rho(0) = \rho$ and the curve $\gamma$ composed 
with one of the embeddings into $C(Z)$ is differentiable at $t = 0$.  We define the tangent space to $\mathcal{M}(Z)$ at the point $\rho \in \mathcal{M}(Z)$ to be the subset $T_{\rho} \mathcal{M}(Z)$
of $C(Z)$ given by the derivatives at $t = 0$ of all admissible 
curves at $\rho$. 
It is not clear a priori that the tangent space $T_{\rho} \mathcal{M}(Z)$ is even a vector space. We show that in fact
$T_{\rho} \mathcal{M}(Z)$ is a closed linear subspace of $C(Z)$. For this, we define the notion of $\rho$-odd functions on $Z$.
Given $\rho \in \mathcal{M}(Z)$, we say that a continuous function $v \in C(Z)$ is {\it $\rho$-odd} if, whenever
$\xi, \eta \in Z$ are such that $\rho(\xi, \eta) = 1$, then $v(\xi) + v(\eta) = 0$. We denote by $\mathcal{O}_{\rho}(Z)$
the set of $\rho$-odd functions, then it is clear that $\mathcal{O}_{\rho}(Z)$ is a closed linear subspace of $C(Z)$.
We have the following:

\medskip

\begin{theorem} \label{mainthm11} Let $Z$ be an antipodal space and let $\rho \in \mathcal{M}(Z)$. Then the tangent space to
$\mathcal{M}(Z)$ at $\rho$ is equal to the set of $\rho$-odd functions,
$$
T_{\rho} \mathcal{M}(Z) = \mathcal{O}_{\rho}(Z).
$$
In particular, the tangent space $T_{\rho} \mathcal{M}(Z)$ is a closed linear subspace of $C(Z)$, and is hence a Banach space.
\end{theorem}

\medskip

Finally, we remark that the study of the Moebius spaces $\mathcal{M}(Z)$ is 
motivated by rigidity problems for negatively curved spaces. It is still an 
open problem whether a Moebius homeomorphism $f : \partial X \to \partial Y$ between the boundaries of proper, geodesically complete CAT(-1) spaces extends to an isometry $F : X \to Y$. As explained in section 5 of \cite{biswas3}, a positive answer to this problem would resolve the longstanding {\it marked length spectrum 
rigidity conjecture} of Burns-Katok \cite{burns-katok}  for closed, negatively curved manifolds. We refer to Hamenstadt \cite{hamenstadt1} for the connection of the marked length spectrum with the geodesic flow, and to Otal \cite{otal1} for the connection with the Moebius structure on the boundary. 

\medskip

Using the Moebius spaces $\MM(Z)$, the problem of extending Moebius maps to isometries can be reformulated as follows. After identifying $X$ and $Y$ with subspaces of $\MM(\partial X)$ and $\MM(\partial Y)$ respectively via the visual embeddings, the Moebius map $f : \partial X \to \partial Y$ extends to an isometry $F : X \to Y$ if and only if the induced isometry $f_* : \MM(\partial X) \to \MM(\partial Y)$ maps the subspace $X$ onto the subspace $Y$, $f_*(X) = Y$. 

\medskip

\medskip

\noindent{\bf Acknowledgements.} The author thanks Mahan Mj for the suggestion to consider whether CAT(0) Gromov hyperbolic
spaces are boundary continuous. The author thanks Arkajit Palchaudhuri for a careful reading of this article and pointing out typos.
The author thanks Urs Lang for suggesting a possible connection between maximal spaces and injective spaces, and for helpful
correspondence regarding injective spaces.

\medskip

\section{Antipodal spaces $Z$ and their associated Moebius spaces $\mathcal{M}(Z)$.}

\medskip

While the latter half of this paper will be mainly concerned with quasi-metric antipodal spaces, the first half will study
the more general class of {\it antipodal spaces} and their associated Moebius spaces. We start with
some generalities about cross-ratios.

\medskip

\begin{definition} {\bf (Separating functions)} Let $Z$ be a compact metrizable space (containing at least four points).
A {\it separating function} on $Z$ is a continuous function $\rho_0 : Z \times Z \to [0, \infty)$ such that:

\medskip

\noindent (1) $\rho_0$ is symmetric, i.e. $\rho_0(\xi, \eta) = \rho_0(\eta, \xi)$ for all $\xi, \eta \in Z$.

\medskip

\noindent (2) $\rho_0$ satisfies positivity, i.e. $\rho_0(\xi, \eta) = 0$ if and only if $\xi = \eta$.

\medskip

For any separating function $\rho_0$, the {\it cross-ratio} with respect to the function $\rho_0$ is
the function of quadruples of distinct points in $Z$ defined by
$$
[\xi, \xi', \eta, \eta']_{\rho_0} := \frac{\rho_0(\xi, \eta)\rho_0(\xi', \eta')}{\rho_0(\xi, \eta')\rho_0(\xi', \eta)}.
$$

Two separating functions $\rho_1, \rho_2$ on $Z$ are said to be {\it Moebius equivalent} if their
cross-ratios agree, $[.,.,.,.]_{\rho_1} \equiv [.,.,.,.]_{\rho_2}$.

\medskip

For a separating function $\rho_0$, the $\rho_0$-diameter of a subset $K \subset Z$ is defined to be
$$
\hbox{diam}_{\rho_0}(K) := \sup_{\xi, \eta \in K} \rho_0(\xi, \eta).
$$
\end{definition}

\medskip

\begin{lemma} \label{derivatives} {\bf (Derivatives and Geometric Mean-Value Theorem)}
Let $\rho_1, \rho_2$ be separating functions on a compact metrizable space $Z$ (containing at least four points). Then $\rho_1$ and $\rho_2$
are Moebius equivalent if and only if there exists a unique positive continuous function
on $Z$, called the derivative of $\rho_2$ with respect to $\rho_1$, denoted by $\frac{d\rho_2}{d\rho_1} : Z \to (0, \infty)$, such that
the following "Geometric Mean-Value Theorem" holds:
\begin{equation} \label{gmvt}
\rho_2(\xi, \eta)^2 = \frac{d\rho_2}{d\rho_1}(\xi)\frac{d\rho_2}{d\rho_1}(\eta) \rho_1(\xi,\eta)^2
\end{equation}
for all $\xi, \eta \in Z$.
Explicitly, the function $\frac{d\rho_2}{d\rho_1}$ is given, for $\xi \in Z$, by
\begin{equation} \label{derivdefn}
\frac{d\rho_2}{d\rho_1}(\xi) = \frac{\rho_2(\xi, \eta) \rho_2(\xi, \eta') \rho_1(\eta, \eta')}{\rho_1(\xi, \eta)\rho_1(\xi, \eta')\rho_2(\eta, \eta')}
\end{equation}
where $\eta, \eta'$ are any two distinct points in $Z$ different from $\xi$.
\end{lemma}

\medskip

\noindent{\bf Proof:} ("If" part:) It is clear from the definition of the cross-ratio that if there exists a function $\frac{d\rho_2}{d\rho_1}$
such that the Geometric Mean-Value Theorem (\ref{gmvt}) holds, then $[.,.,.,.]_{\rho_2} = [.,.,.,.]_{\rho_1}$ and $\rho_1, \rho_2$ are
Moebius equivalent.

\medskip

\noindent ("Only if" part:) Given that $\rho_1, \rho_2$ are Moebius equivalent, let us denote the right-hand side of equation (\ref{derivdefn})
by $A(\xi, \eta, \eta')$, and define a function $\frac{d\rho_2}{d\rho_1}$ on $Z$ by
the equation (\ref{derivdefn}), i.e. $\frac{d\rho_2}{d\rho_1}(\xi) := A(\xi, \eta, \eta')$, where $\eta, \eta' \in Z$ are any two
distinct points different from $\xi$.
Then $\frac{d\rho_2}{d\rho_1}(\xi)$ is well-defined, independent of the choice of points $\eta, \eta'$ different
from $\xi$, because if we change one of the points $\eta$ to some $\zeta$, then the equality
$[\xi, \eta', \eta, \zeta]_{\rho_2} = [\xi, \eta', \eta, \zeta]_{\rho_1}$ implies that
$$
A(\xi, \eta, \eta') = A(\xi, \zeta, \eta'),
$$
and if we then change $\eta'$ to some $\zeta'$ then a similar argument gives
$$
A(\xi, \zeta, \eta') = A(\xi, \zeta, \zeta'),
$$
hence
$$
A(\xi, \eta, \eta') = A(\xi, \zeta, \zeta'),
$$
and so $\frac{d\rho_2}{d\rho_1}$ is well-defined. It is clear from equation (\ref{derivdefn}) that $\frac{d\rho_2}{d\rho_1}$ is
positive and continuous.

\medskip

Given $\xi, \eta \in Z$, if $\xi = \eta$ then (\ref{gmvt}) holds trivially. Otherwise, for $\xi \neq \eta$, choose $\eta'$ a point distinct
from both $\xi, \eta$, then by definition of the derivative $\frac{d\rho_2}{d\rho_1}$, we have
$$
\frac{d\rho_2}{d\rho_1}(\xi) = \frac{\rho_2(\xi, \eta) \rho_2(\xi, \eta') \rho_1(\eta, \eta')}{\rho_1(\xi, \eta)\rho_1(\xi, \eta')\rho_2(\eta, \eta')}
$$
and
$$
\frac{d\rho_2}{d\rho_1}(\eta) = \frac{\rho_2(\eta, \xi) \rho_2(\eta, \eta') \rho_1(\xi, \eta')}{\rho_1(\eta, \xi)\rho_1(\eta, \eta')\rho_2(\xi, \eta')},
$$
multiplying the above two equations gives
$$
\frac{d\rho_2}{d\rho_1}(\xi) \frac{d\rho_2}{d\rho_1}(\eta) = \frac{\rho_2(\xi, \eta)^2}{\rho_1(\xi, \eta)^2},
$$
which is the Geometric Mean-Value Theorem, as required.

\medskip

Finally, for the uniqueness of the function $\frac{d\rho_2}{d\rho_1}$, note that if a function $f$ satisfies the Geometric Mean-Value Theorem
$\rho_2(\xi, \eta)^2 = f(\xi) f(\eta) \rho_1(\xi, \eta)^2$, then it must be given by
$$
f(\xi) = \frac{\rho_2(\xi, \eta) \rho_2(\xi, \eta') \rho_1(\eta, \eta')}{\rho_1(\xi, \eta)\rho_1(\xi, \eta')\rho_2(\eta, \eta')}
$$
where $\eta, \eta'$ are any two distinct points different from $\eta$, hence $f = \frac{d\rho_2}{d\rho_1}$.
$\diamond$.

\medskip

It follows from the Geometric Mean-Value Theorem and the continuity of $\frac{d\rho_2}{d\rho_1}$ that
$$
\frac{d\rho_2}{d\rho_1}(\xi) = \lim_{\eta \to \xi} \frac{\rho_2(\xi, \eta)}{\rho_1(\xi, \eta)}
$$
for all non-isolated points $\xi$ of $Z$, hence justifying the name "derivative" for the function $\frac{d\rho_2}{d\rho_1}$.

\medskip

\begin{definition} {\bf (Unrestricted Moebius space of a separating function)} Let $\rho_0$ be a separating function on
a compact metrizable space $Z$ containing at least four points. The Unrestricted Moebius space of the separating function $\rho_0$ is defined to be
$$
\mathcal{UM}(Z, \rho_0) := \{ \rho | \ \rho \ \hbox{is a separating function on } Z \ \hbox{Moebius equivalent to } \rho_0 \},
$$
with the metric $d_{\mathcal{M}}$ on $\mathcal{UM}(Z, \rho_0)$ defined by
$$
d_{\mathcal{M}}(\rho_1, \rho_2) := \left|\left| \log \frac{d\rho_2}{d\rho_1} \right|\right|_{\infty}.
$$
\end{definition}

\medskip

The Geometric Mean-Value Theorem implies that for $\rho_1, \rho_2 \in \mathcal{UM}(Z, \rho_0)$ we have
$$
\rho_1 = \rho_2 \ \hbox{ if and only if } \ \frac{d\rho_2}{d\rho_1} \equiv 1,
$$
hence $d_{\mathcal{M}}(\rho_1, \rho_2) = 0$ if and only if $\rho_1 = \rho_2$.
From equation (\ref{derivdefn}), we also obtain easily the Chain rule for derivatives: if $\rho_1, \rho_2, \rho_3 \in \mathcal{UM}(Z, \rho_0)$,
then
$$
\frac{d\rho_3}{d\rho_1}(\xi) = \frac{d\rho_3}{d\rho_2}(\xi) \cdot \frac{d\rho_2}{d\rho_1}(\xi)
$$
for all $\xi \in Z$, from which the triangle inequality for the function $d_{\mathcal{M}}$ follows. The Chain Rule also implies
$$
\frac{d\rho_2}{d\rho_1} \cdot \frac{d\rho_1}{d\rho_2} \equiv 1
$$
for all $\rho_1, \rho_2 \in \mathcal{UM}(Z, \rho_0)$,
from which symmetry of the function $d_{\mathcal{M}}$ follows. Thus the function $d_{\mathcal{M}}$ is indeed
a metric on $\mathcal{UM}(Z, \rho_0)$.

\medskip

Lemma \ref{derivatives} gives an easy way to construct separating functions Moebius equivalent to a given separating function $\rho_0$, i.e.
points of $\mathcal{UM}(Z, \rho_0)$. Indeed, given any $\rho_1 \in \mathcal{UM}(Z, \rho_0)$ and any continuous function $\tau$ on $Z$, we obtain a
separating function $\rho = E_{\rho_1}(\tau) \in \mathcal{UM}(Z, \rho_0)$ defined by
\begin{equation} \label{Esubrho}
E_{\rho_1}(\tau)(\xi, \eta) := e^{\frac{1}{2}\tau(\xi)} e^{\frac{1}{2}\tau(\eta)} \rho_1(\xi, \eta), \ \xi, \eta \in Z,
\end{equation}
which has the property that
$$
\log \frac{d\rho}{d\rho_1} = \tau.
$$
It is easy to see from this and the Chain Rule that, denoting by $C(Z)$ the Banach space of continuous functions on $Z$
equipped with the sup norm, the map
\begin{align*}
E_{\rho_1} : C(Z) & \to \mathcal{UM}(Z, \rho_0) \\
             \tau & \mapsto E_{\rho_1}(\tau)
\end{align*}
is a surjective isometry, with inverse given by the isometry
\begin{align*}
\phi_{\rho_1} : \mathcal{UM}(Z, \rho_0) & \to C(Z) \\
                              \rho               & \mapsto \log \frac{d\rho}{d\rho_1}.
\end{align*}

Thus our definition of $\mathcal{UM}(Z, \rho_0)$ is too flexible to obtain an interesting space reflecting the properties of the separating function $\rho_0$,
since the metric space $\mathcal{UM}(Z, \rho_0)$ is always naturally isometric to the Banach space $C(Z)$, independent of the choice of $\rho_0$. Note that
consequently $\mathcal{UM}(Z, \rho_0)$ is also never proper when $Z$ is infinite.

\medskip

In order to get interesting, proper spaces of Moebius functions, we now introduce the key property which will be crucial for the rest of this
article:

\medskip

\begin{definition} {\bf (Antipodal function)} Let $Z$ be a compact metrizable space. An {\it antipodal function} on $Z$
is a separating function $\rho_0$ which satisfies the following two properties:

\medskip

\noindent (1) $\rho_0$ has diameter one, i.e. $\sup_{\xi, \eta \in Z} \rho_0(\xi, \eta) = 1$.

\medskip

\noindent (2) $\rho_0$ is antipodal, i.e. for all $\xi \in Z$ there exists $\eta \in Z$ such that $\rho_0(\xi, \eta) = 1$.

\medskip

An {\it antipodal space} is defined to be a compact metrizable space equipped with an antipodal function, $(Z, \rho_0)$.
\end{definition}

\medskip

The spaces of interest for us will be the following:

\medskip

\begin{definition} {\bf (Moebius space of an antipodal function)} Let $(Z, \rho_0)$ be an antipodal space. The Moebius space
associated to $(Z, \rho_0)$ is defined to be the space
$$
\mathcal{M}(Z, \rho_0) := \{ \rho | \ \rho \ \hbox{is an antipodal function Moebius equivalent to} \ \rho_0 \}
$$
equipped with the metric $d_{\mathcal{M}}$ defined as before (note that $\mathcal{M}(Z, \rho_0)$ is a subset of $\mathcal{UM}(Z, \rho_0)$).
\end{definition}

\medskip

We will often write just $\mathcal{M}(Z), \mathcal{UM}(Z)$ instead of $\mathcal{M}(Z, \rho_0), \mathcal{UM}(Z, \rho_0)$.
We have:

\medskip

\begin{lemma} \label{imageclosed} Let $\rho_1 \in \mathcal{M}(Z)$. The map
\begin{align*}
\phi_{\rho_1} : \mathcal{M}(Z) & \to C(Z) \\
                     \rho      & \mapsto \log \frac{d\rho}{d\rho_1}
\end{align*}
is an isometric embedding onto a closed subset of $C(Z)$.
\end{lemma}

\medskip

\noindent{\bf Proof:} Since this map is just the restriction of the map $\phi_{\rho_1}$ to the subspace
$\mathcal{M}(Z, \rho_0)$ of $\mathcal{UM}(Z, \rho_0)$, it is clearly an isometric embedding.

\medskip

To show that the image in $C(Z)$ is closed, let $\rho_n \in \mathcal{M}(Z)$ be a
sequence such that $\tau_n = \log \frac{d\rho_n}{d\rho_1} \in C(Z)$ converges in $C(Z)$ to some $\tau \in C(Z)$.
Let $\rho = E_{\rho_1}(\tau) \in \mathcal{UM}(Z, \rho_0)$, then $\log \frac{d\rho}{d\rho_1} = \tau$,
so it's enough to show that $\rho \in \mathcal{M}(Z, \rho_0)$.

\medskip

Since $\tau_n \to \tau$ uniformly on $Z$, by the Geometric Mean-Value Theorem we have $\rho_n(\xi, \eta) \to \rho(\xi, \eta)$
as $n \to \infty$ for all $\xi, \eta \in Z$. Since $\rho_n(\xi, \eta) \leq 1$ for all $n$, it follows that $\rho(\xi, \eta) \leq 1$
for all $\xi, \eta$. Also, given any $\xi \in Z$, for any $n$ there exists $\eta_n \in Z$ such that $\rho_n(\xi, \eta_n) = 1$
(since $\rho_n$ is antipodal), passing to a subsequence we may assume $\eta_n \to \eta$ for some $\eta$, then by the Geometric Mean-Value Theorem
and the uniform convergence of $\tau_n$ to $\tau$ we have
$$
\rho(\xi, \eta) = \lim_{n \to \infty} \rho_n(\xi, \eta_n) = 1,
$$
hence $\rho$ is diameter one and antipodal, i.e. $\rho \in \mathcal{M}(Z)$. $\diamond$

\medskip

Unlike the space $\mathcal{UM}(Z, \rho_0)$, the space $\mathcal{M}(Z, \rho_0)$
is always proper:

\medskip

\begin{lemma} \label{properness} The metric space $(\mathcal{M}(Z), d_{\mathcal{M}})$ is proper,
i.e. closed and bounded balls in $\mathcal{M}(Z)$ are compact. In particular $\mathcal{M}(Z)$ is complete.
\end{lemma}

\medskip

\noindent{\bf Proof:} Let $\{\rho_n\}$ be a bounded sequence in $\mathcal{M}(Z)$, say $d_{\mathcal{M}}(\rho_n, \rho_0) \leq M$
for all $n$, for some $M > 0$. Let $\tau_n = \log \frac{d\rho_n}{d\rho_0}$.
It's enough to show that $\{ \rho_n \}$ has a subsequence which is convergent in $\mathcal{M}(Z)$, for which, by the
previous Lemma, it is enough to show that $\{\tau_n\}$ has a subsequence which converges in $C(Z)$.

\medskip

Let $\alpha$ be a metric on $Z$ inducing the topology of $Z$.
Since $||\tau_n||_{\infty} \leq M$ for all $n$, by the Arzela-Ascoli Theorem it is enough to show that the
sequence $\{\tau_n\}$ is equicontinuous. Suppose not. Then there exists $\epsilon > 0$, such that after passing to a
subsequence, we have sequences $\xi_{1,n}, \xi_{2,n} \in Z$ such that $\alpha(\xi_{1,n}, \xi_{2,n}) \to 0$, but
$\tau_n(\xi_{2,n}) - \tau_n(\xi_{1,n}) \geq \epsilon$ for all $n$.

\medskip

For each $n$ choose $\eta_n \in Z$ such that $\rho_n(\xi_{1,n}, \eta_n) = 1$. Since $||\tau_n||_{\infty} \leq M$, by the
Geometric Mean-Value Theorem we have $\rho_0(\xi_{1,n}, \eta_n) \geq e^{-M}$ for all $n$. Thus passing to a subsequence,
we may assume there are $\xi, \eta \in Z$ such that $\rho_0(\xi, \eta) \geq e^{-M}$ (in particular $\xi \neq \eta$) and
such that $\xi_{1,n} \to \xi, \eta_n \to \eta$ as $n \to \infty$. Since $\alpha(\xi_{1,n}, \xi_{2,n}) \to 0$, we also
have $\xi_{2,n} \to \xi$ as $n \to \infty$. Then, by the Geometric Mean-Value Theorem,
\begin{align*}
\tau_n(\xi_{2,n}) + \tau_n(\eta_n) + \log \rho_0(\xi_{2,n}, \eta_n)^2 & = \log \rho_n(\xi_{2,n}, \eta_n)^2 \\
                                                                      & \leq 0 = \log \rho_n(\xi_{1,n}, \eta_n)^2 \\
                                                                      & = \tau_n(\xi_{1,n}) + \tau_n(\eta_n) + \log \rho_0(\xi_{1,n}, \eta_n)^2,
\end{align*}
hence
\begin{align*}
\tau_n(\xi_{2,n}) - \tau_n(\xi_{1,n}) & \leq \log \left(\frac{\rho_0(\xi_{1,n}, \eta_n)^2}{\rho_0(\xi_{2,n}, \eta_n)^2}\right) \\
                                      & \to \log \left(\frac{\rho_0(\xi, \eta)^2}{\rho_0(\xi, \eta)^2}\right) = 0, \ \hbox{as} \ n \to \infty,
\end{align*}
since $\xi_{1,n}, \xi_{2,n} \to \xi$, and $\eta_n \to \eta$, and $\rho_0$ is continuous. But this contradicts
$\tau_n(\xi_{2,n}) - \tau_n(\xi_{1,n}) \geq \epsilon$ for all $n$. $\diamond$

\medskip

For a continuous function $f$ on a compact space, we will denote by argmax $f$ (resp. argmin $f$) the closed set where the function
$f$ attains its maximum value (resp. minimum value). We will need the following Lemma:

\medskip

\begin{lemma} \label{maxmin} For any $\rho_1, \rho_2 \in \mathcal{M}(Z)$ we have
$$
\max_{\xi \in Z} \frac{d\rho_2}{d\rho_1}(\xi) \cdot \min_{\xi \in Z}
\frac{d\rho_2}{d\rho_1}(\xi) = 1.
$$
In particular,
$$
d_{\mathcal{M}}(\rho_1, \rho_2) = \max_{\xi \in Z} \log \frac{d\rho_2}{d\rho_1}(\xi) = - \min_{\xi \in Z} \log \frac{d\rho_2}{d\rho_1}(\xi).
$$

Moreover:

\medskip

\noindent (1) If $\xi \in \hbox{argmax} \ \frac{d\rho_2}{d\rho_1}$ and $\eta \in Z$ is such that $\rho_1(\xi, \eta) = 1$ then $\eta \in \hbox{argmin} \ \frac{d\rho_2}{d\rho_1}$, and $\rho_2(\xi, \eta) = \rho_1(\xi, \eta) = 1$.

\medskip

\noindent (2) If $\xi' \in \hbox{argmin} \ \frac{d\rho_2}{d\rho_1}$ and $\eta' \in Z$ is such that $\rho_2(\xi', \eta') = 1$ then $\eta' \in \hbox{argmax} \ \frac{d\rho_2}{d\rho_1}$, and $\rho_2(\xi', \eta') = \rho_1(\xi', \eta') = 1$.
\end{lemma}

\medskip

\noindent{\bf Proof:} Let $\lambda, \mu$ denote the maximum and
minimum values of $\frac{d\rho_2}{d\rho_1}$ respectively, and let $\xi \in \hbox{argmax} \ \frac{d\rho_2}{d\rho_1},
\xi' \in \hbox{argmin} \ \frac{d\rho_2}{d\rho_1}$. Given $\eta \in Z$ such that $\rho_1(\xi, \eta) =
1$, we have, using the Geometric Mean-Value Theorem, and $\frac{d\rho_2}{d\rho_1}(\eta) \geq \mu$, that
$$
1 \geq \rho_2(\xi, \eta)^2 = \frac{d\rho_2}{d\rho_1}(\xi)
\frac{d\rho_2}{d\rho_1}(\eta) \geq \lambda \cdot
\mu,
$$
while letting $\eta' \in Z$ be such that $\rho_2(\xi', \eta') = 1$, we have, again by the Geometric Mean-Value Theorem, and $\frac{d\rho_2}{d\rho_1}(\eta') \leq \lambda$, that
$$
1 \geq \rho_1(\xi', \eta')^2 =
1/\left(\frac{d\rho_2}{d\rho_1}(\xi')
\frac{d\rho_2}{d\rho_1}(\eta')\right) \geq 1/(\lambda
\mu),
$$
hence combining the two inequalities gives $\lambda \cdot \mu = 1$. It follows that we have equality everywhere in the above inequalities, and hence
$\frac{d\rho_2}{d\rho_1}(\eta) = \mu$, so $\eta \in \hbox{argmin} \ \frac{d\rho_2}{d\rho_1}$, and also
$\rho_2(\xi, \eta) = \rho_1(\xi, \eta) = 1$. Similarly we have $\frac{d\rho_2}{d\rho_1}(\eta') = \lambda$, so $\eta' \in \hbox{argmax} \ \frac{d\rho_2}{d\rho_1}$, and also
$\rho_1(\xi', \eta') = \rho_2(\xi', \eta') = 1$.
$\diamond$

\medskip

\section{The antipodal flow and the antipodalization map $\mathcal{P}_{\infty}$.}

\medskip

As remarked in the previous section, given a separating function $\rho_0$ on a compact metrizable space $Z$, it is trivial
to construct separating functions $\rho$ Moebius equivalent to $\rho_0$, one simply takes any continuous function $\tau \in C(Z)$
and puts $\rho = E_{\rho_0}(\tau)$, where $E_{\rho_0}(\tau)$ is defined by equation (\ref{Esubrho}).

\medskip

On the other hand, given an antipodal function $\rho_0$, of course by definition $\rho_0$ belongs to the Moebius space $\mathcal{M}(Z, \rho_0)$,
but constructing even a single point of $\mathcal{M}(Z, \rho_0)$ different from $\rho_0$ is quite a non-trivial matter, as one has to construct
separating functions Moebius equivalent to $\rho_0$ which have to be antipodal functions, i.e. they have to have diameter one, and they have to be antipodal.
Making a given separating function diameter one is easily achieved just by a rescaling, what is non-trivial is to ensure the antipodal condition.
The purpose of this section, which is the technical heart of the paper, is to describe a method for constructing antipodal Moebius equivalent functions
starting from arbitrary Moebius equivalent functions. These are obtained as limits as $t \to \infty$ of solutions $(\tau_t)_{t \geq 0}$ of a certain
ODE in the Banach space $C(Z)$, an approach reminiscent of the classical heat flow method for constructing harmonic functions.

\medskip

The key object is the following:

\medskip

\begin{definition} {\bf (The discrepancy function)} Let $(Z, \rho_0)$ be an antipodal space.
Given $\rho_1 \in \mathcal{M}(Z)$ and $\tau \in C(Z)$, the {\it discrepancy function} of $\tau$
with respect to $\rho_1$, is the continuous function $D_{\rho_1}(\tau) \in C(Z)$ defined by
$$
D_{\rho_1}(\tau)(\xi) := \max_{\eta \in Z - \{\xi\}} \tau(\xi) + \tau(\eta) + \log \rho_1(\xi, \eta)^2.
$$
\end{definition}

It is easy to see that the function $D_{\rho_1}(\tau)$ is indeed a continuous function on $Z$. The function $D_{\rho_1}(\tau)$
measures the failure of the Moebius equivalent separating function $\rho = E_{\rho_1}(\tau) \in \mathcal{UM}(Z)$ to be an
antipodal function, as it is easy to see that
$$
D_{\rho_1}(\tau)(\xi) = \max_{\eta \in Z - \{\xi\}} \log \rho(\xi, \eta)^2,
$$
and hence
$$
\rho \ \hbox{is an antipodal function, i.e.} \ \rho \in \mathcal{M}(Z), \ \hbox{if and only if} \ D_{\rho_1}(\tau) \equiv 0.
$$
For $\rho_1$ fixed, we can consider the map $D_{\rho_1} : C(Z) \to C(Z), \tau \mapsto D_{\rho_1}(\tau)$,  as
a vector field on the Banach space $C(Z)$. The remarks above mean that under the isometry $\phi_{\rho_1} : \mathcal{UM}(Z) \to C(Z)$,
the subspace $\mathcal{M}(Z) \subset \mathcal{UM}(Z)$ corresponds to the subspace $\{ \tau \in C(Z) \ | \ D_{\rho_1}(\tau) = 0 \} \subset C(Z)$.

\medskip

The map $D_{\rho_1} : C(Z) \to C(Z)$ is continuous, in fact $2$-Lipschitz. We have
the following estimates which will be useful:

\medskip

\begin{lemma} \label{discrepancylipschitz} Let $\rho_1 \in \mathcal{M}(Z)$. Then we have:
\begin{equation} \label{discrepbound}
||D_{\rho_1}(\tau)||_{\infty} \leq 2 ||\tau||_{\infty}\ , \ \hbox{for all } \ \tau \in C(Z),
\end{equation}
and
\begin{equation} \label{discreplip}
||D_{\rho_1}(\tau_0) - D_{\rho_1}(\tau_1)||_{\infty} \leq 2 ||\tau_0 - \tau_1||_{\infty}\ , \ \hbox{for all } \ \tau_0, \tau_1 \in C(Z).
\end{equation}
\end{lemma}

\medskip

\noindent{\bf Proof:} Given $\tau \in C(Z)$, and $\xi \in Z$, since $\rho_1$ is of diameter one we have, for all $\eta \in Z - \{\xi\}$,
\begin{align*}
\tau(\xi) + \tau(\eta) + \log \rho_1(\xi, \eta)^2 & \leq \tau(\xi) + \tau(\eta) \\
                                                  & \leq 2 ||\tau||_{\infty},
\end{align*}
and so
$$
D_{\rho_1}(\tau)(\xi) \leq 2||\tau||_{\infty},
$$
while choosing a point $\eta' \in Z$ such that $\rho_1(\xi, \eta') = 1$ gives
$$
D_{\rho_1}(\tau)(\xi) \geq \tau(\xi) + \tau(\eta') \geq - 2||\tau||_{\infty},
$$
and thus $||D_{\rho_1}(\tau)||_{\infty} \leq 2||\tau||_{\infty}$.

\medskip

To prove (\ref{discreplip}) above, let $\tau_0, \tau_1 \in C(Z)$. Let $\mu$ be a probability measure on $Z$ whose support equals all of $Z$.
Define, for any $\tau \in C(Z)$, and $1 \leq p < \infty$, functions
$$
D_{\rho_1,p}(\tau)(\xi) := \frac{1}{p} \log\left(\int_Z e^{p\tau(\xi)}e^{p\tau(\eta)} \rho_1(\xi, \eta)^{2p} d\mu(\eta)\right) \ , \ \xi \in Z.
$$
Then it is easy to see that $D_{\rho_1,p}(\tau)$ converges pointwise on $Z$ to $D_{\rho_1}(\tau)$ as $p \to \infty$.

\medskip

For $t \in [0,1]$, let $\tau_t = (1 - t) \tau_0 + t \tau_1 \in C(Z)$. Then it is not hard to see that for fixed $\xi \in Z$ and fixed $p \in [1, \infty)$,
the function $t \mapsto D_{\rho_1,p}(\tau_t)(\xi)$ is $C^1$, and its derivative satisfies
\begin{align*}
\left| \frac{d}{dt} D_{\rho_1,p}(\tau_t)(\xi) \right| & = \frac{\left| \int_Z ((\tau_1 - \tau_0)(\xi) + (\tau_1 - \tau_0)(\eta))e^{p\tau_t(\xi)}e^{p\tau_t(\eta)} \rho_1(\xi, \eta)^{2p} d\mu(\eta)\right|}{\int_Z e^{p\tau_t(\xi)}e^{p\tau_t(\eta)} \rho_1(\xi, \eta)^{2p} d\mu(\eta)} \\
                                                     & \leq 2 ||\tau_1 - \tau_0||_{\infty},
\end{align*}
from which we obtain by integrating
$$
|D_{\rho_1,p}(\tau_0)(\xi) - D_{\rho_1,p}(\tau_1)(\xi)| \leq 2 ||\tau_1 - \tau_0||_{\infty},
$$
for all $\xi \in Z, p \in [1,\infty)$. Letting $p \to \infty$ in this inequality gives
$$
|D_{\rho_1}(\tau_0)(\xi) - D_{\rho_1}(\tau_1)(\xi)| \leq 2 ||\tau_1 - \tau_0||_{\infty},
$$
for all $\xi \in Z$, and so (\ref{discreplip}) holds. $\diamond$

\medskip

We now define the ODE previously alluded to:

\medskip

\begin{definition} {\bf (Antipodal flow)} Fix an element $\rho_1 \in \mathcal{M}(Z)$. The {\it $\rho_1$-antipodal flow} on $C(Z)$ is
the ODE on $C(Z)$ defined by the vector field $-D_{\rho_1} : C(Z) \to C(Z)$ on $C(Z)$, namely the ODE
\begin{equation} \label{basicode}
\frac{d}{dt} \tau_t = - D_{\rho_1}(\tau_t).
\end{equation}
By a solution of the ODE (\ref{basicode}) we mean a $C^1$ curve in the Banach space $C(Z)$, $t \in I \subset \mathbb{R} \mapsto \tau_t \in C(Z)$, which
satisfies equation (\ref{basicode}).
\end{definition}

\medskip

It is well known that for locally bounded, locally Lipschitz vector fields on a Banach space $E$, we have the same local Existence-Uniqueness Theorems
for the corresponding ODE, as well as the existence and uniqueness of maximal integral curves for any initial condition, just as in $\mathbb{R}^n$.
For the sake of completeness, we state the relevant theorems in a form suitable for our use:

\medskip

\begin{theorem} \label{banachode1} {\bf (Local Existence and Uniqueness)} Let $E$ be a Banach space and let $X : E \to E$ be a map such that $X$ is locally bounded and locally Lipschitz. Then for any $x_0 \in E$, there exists $\epsilon_0 > 0$ and $\delta > 0$, such that for all $x \in B(x_0, \delta)$, and for all $\epsilon \in (0, \epsilon_0]$,
there exists a unique $C^1$ curve $\gamma : (-\epsilon, \epsilon) \to E$ such that $\gamma(0) = x$ and
\begin{equation} \label{generalode}
\gamma'(t) = X(\gamma(t)) \ ,
\end{equation}
for all $t \in (-\epsilon, \epsilon)$.
\end{theorem}

\medskip

For a proof of the above theorem, which follows the same lines as the proof in $\mathbb{R}^n$, one may look at for example \cite{mukherjea}, section 5.2.
Recall that an integral curve $\gamma$ of (\ref{generalode}) defined on an interval $I$ is called {\it maximal} if there does not exist an interval $J$ strictly containing $I$
and an integral curve $\tilde{\gamma}$ on $J$ such that $\tilde{\gamma} = \gamma$ on $I$. From the above theorem it is straightforward to prove the existence and uniqueness of maximal integral curves of such an ODE, and we have:

\medskip

\begin{theorem} \label{maximalintegral} {\bf (Maximal integral curves)} Let $E$ be a Banach space and let $X : E \to E$ be a map such that $X$ is locally bounded and locally Lipschitz. Then for any $x_0 \in E$, there is a unique maximal integral curve $\gamma : I \to E$
of the vector field $X$ defined on an open interval $I$ containing $0$ such that $\gamma(0) = x_0$. This means that if $\tilde{\gamma} : J \to E$ is an
integral curve of $X$ defined on an open interval $J$ containing $0$ such that $\tilde{\gamma}(0) = x_0$, then $J \subset I$ and $\tilde{\gamma} = \gamma$ on $J$.
\end{theorem}

\medskip

As an immediate corollary of the above theorems, we have:

\medskip

\begin{lemma} \label{antiflowmax} Let $\rho_1 \in \mathcal{M}(Z)$. Then for any $\tau_0 \in C(Z)$, there exists a unique maximal solution
$\gamma_{\tau_0} : I_{\tau_0} \to C(Z)$ of the $\rho_1$-antipodal flow, i.e. the ODE (\ref{basicode}), such that $\gamma_{\tau_0}(0) = \tau_0$. Here
$I_{\tau_0}$ is an open interval containing $0$, and is the unique maximal interval on which a solution satisfying $\gamma(0) = \tau_0$ exists.
\end{lemma}

\medskip

\noindent{\bf Proof:} From (\ref{discrepbound}) and (\ref{discreplip}) of Lemma \ref{discrepancylipschitz}, it follows that the vector field
$X = -D_{\rho_1}$ on the Banach space $E = C(Z)$ is locally bounded and locally Lipschitz, hence the Lemma follows immediately from Theorem \ref{maximalintegral}
above. $\diamond$

\medskip

We now proceed to show that the maximal solutions $\gamma_{\tau_0}$ are defined for all positive time, i.e. $[0,\infty) \subset I_{\tau_0}$.

\medskip

\begin{lemma} \label{gronwallbound} Fix $\rho_1 \in \mathcal{M}(Z)$. Let $\gamma : (a, b) \to C(Z)$ be a solution of the ODE (\ref{basicode}),
where $-\infty < a < 0 < b < \infty$. Then for all $t \in [0, b)$, we have
$$
||\gamma(t)||_{\infty} \leq ||\gamma(0)||_{\infty} e^{2t}.
$$
\end{lemma}

\medskip

\noindent{\bf Proof:} For $t \in [0, b)$, let $u(t) = \int_{0}^{t} ||\gamma(s)||_{\infty} ds$. Since $\gamma$ is a solution of (\ref{basicode}), we have,
for $t \in [0, b)$, using (\ref{discrepbound}), that
\begin{align*}
||\gamma(t)||_{\infty} & \leq ||\gamma(0)||_{\infty} + \int_{0}^{t} ||D_{\rho_1}(\gamma(s))||_{\infty} ds \\
                       & \leq ||\gamma(0)||_{\infty} + 2 \int_{0}^{t} ||\gamma(s)||_{\infty} ds,
\end{align*}
and hence
\begin{equation} \label{uanduprime}
u'(t) \leq C + 2u(t)
\end{equation}
(where $C = ||\gamma(0)||_{\infty}$), or
$$
\frac{d}{dt}\left(e^{-2t}u(t)\right) \leq Ce^{-2t},
$$
which gives on integrating (and noting that $u(0) = 0$)
$$
e^{-2t}u(t) \leq \frac{C}{2}(1 - e^{-2t}),
$$
so using (\ref{uanduprime}) we get finally
\begin{align*}
||\gamma(t)||_{\infty} = u'(t) & \leq C + 2u(t) \\
                               & \leq C + 2 \cdot (C/2) (e^{2t} - 1) \\
                               & = C e^{2t}
\end{align*}
as required. $\diamond$

\medskip

The long-term existence of the antipodal flow is a straightforward consequence of the previous lemma:

\medskip

\begin{prop} \label{longterm} {\bf (Long-term existence of antipodal flow)} Let $\rho_1 \in \mathcal{M}(Z)$. Given
$\tau_0 \in C(Z)$, let $\gamma : I_{\tau_0} \to C(Z)$ be the unique maximal solution of the $\rho_1$-antipodal flow
such that $\gamma(0) = \tau_0$. Then $\gamma$ is defined for all time $t \geq 0$, i.e. $[0, \infty) \subset I_{\tau_0}$.
\end{prop}

\medskip

\noindent{\bf Proof:} Let $a = \inf I_{\tau_0}$ and $b = \sup I_{\tau_0} \in (0, \infty]$. It is enough to show that $b = +\infty$. Suppose not, then
$0 < b < +\infty$, and $\gamma$ is defined on the interval $[0,b)$. By the previous Lemma \ref{gronwallbound} and (\ref{discrepbound}), we have
\begin{align*}
\lim_{t \to b^{-}} \int_{0}^{t} ||\gamma'(s)||_{\infty} \ ds & =  \lim_{t \to b^{-}} \int_{0}^{t} ||D_{\rho_1}(\gamma(s))||_{\infty} \ ds \\
                                                           & \leq \lim_{t \to b^{-}} 2 \int_{0}^{t} ||\gamma(s)||_{\infty} \ ds \\
                                                           & \leq \lim_{t \to b^{-}} 2 \int_{0}^{t} ||\gamma(0)||_{\infty} \ e^{2s} \ ds \\
                                                           & = ||\gamma(0)||_{\infty} \ (e^{2b} - 1) \\
                                                           & < +\infty.
\end{align*}
Since $C(Z)$ is a Banach space, this implies that the limit of the $C(Z)$-valued Riemann integrals
$$
\lim_{t \to b^{-}} \int_{0}^{t} \gamma'(s) \ ds
$$
exists in $C(Z)$, hence by the Fundamental Theorem of Calculus the limit
$$
\lim_{t \to b^{-}} \gamma(t)
$$
exists in $C(Z)$, and equals some $\tau_1 \in C(Z)$ say.

\medskip

Now by the Local Existence and Uniqueness Theorem (Theorem \ref{banachode1}), there exists $\epsilon_0, \delta > 0$ such that for any
$\tau \in B(\tau_1, \delta)$ there exists a unique solution $\alpha_{\tau} : (-\epsilon_0, \epsilon_0) \to C(Z)$ of the ODE (\ref{basicode}) such that
$\alpha_{\tau}(0) = \tau$.

\medskip

Since $\gamma(t) \to \tau_1$ as $t \to b^{-}$, we can
choose a $t_1 \in (b - \epsilon_0, b)$ such that $\gamma(t_1) \in B(\tau_1, \delta)$. Let $\tau = \gamma(t_1)$ and let
$\alpha_{\tau} : (-\epsilon_0, \epsilon_0) \to C(Z)$ be the unique solution of (\ref{basicode}) such that $\alpha_{\tau}(0) = \tau$.
Then by uniqueness of solutions to (\ref{basicode}), since $\gamma(t_1) = \tau = \alpha_{\tau}(0)$, we must have
\begin{equation} \label{solnequal}
\gamma(t_1 + s) = \alpha_{\tau}(s)
\end{equation}
for all $s$ such that $|s| < \epsilon_0$ and $t_1 + s \in I_{\tau_0}$.
Hence we can define a solution $\tilde{\gamma}$ of (\ref{basicode}) on the interval $J = (a, t_1+\epsilon_0)$ by
putting
\begin{align*}
\tilde{\gamma}(t) & := \gamma(t) \ , \ t \in (a, t_1) \\
                  & := \alpha_{\tau}(t - t_1) \ , \ t \in (t_1 - \epsilon_0, t_1+\epsilon_0),
\end{align*}
which is well-defined by (\ref{solnequal}). Since $t_1 + \epsilon_0 > b$, the interval $J$ strictly contains the interval $I_{\tau_0}$,
a contradiction to the maximality of the solution $\gamma$. $\diamond$

\medskip

We now proceed towards an analysis of the long-term behaviour of the antipodal flow.

\medskip

For a real-valued function $u$ on an interval in $\mathbb{R}$, we will denote by $\frac{d^+ u}{dt}$ the right-hand derivative
$$
\frac{d^+ u}{dt} := \lim_{h \to 0^+} \frac{u(t+h) - u(t)}{h}
$$
whenever it exists. We state without proof the following elementary lemma:

\medskip

\begin{lemma} \label{righthandderiv} Let $u : [a, b] \to \mathbb{R}$ be such that for all $t \in [a, b)$, the right-hand derivative
$\frac{d^+ u}{dt}$ exists and satisfies
$$
\frac{d^+ u}{dt} \leq 0.
$$
Then
$$
u(b) \leq u(a).
$$
\end{lemma}

\medskip

Given a separating function $\rho$ on $Z$ and $\xi, \eta \in Z$, we will write $\xi \sim_{\rho} \eta$ if
$$
\eta \in \ \hbox{argmax} \ \rho(\xi, .),
$$
i.e.
$$
\rho(\xi, \eta) = \max_{\zeta \in Z} \rho(\xi, \zeta).
$$
With this notation, given $\rho_1 \in \mathcal{M}(Z)$ and $\tau \in C(Z)$, letting $\rho = E_{\rho_1}(\tau)$,
it is clear that for any $\xi \in Z$, we have
$$
D_{\rho_1}(\tau)(\xi) = \log \rho(\xi, \eta)^2 = \tau(\xi) + \tau(\eta) + \log \rho_1(\xi, \eta)^2,
$$
for any $\eta \in Z$ such that $\xi \sim_{\rho} \eta$.

\medskip

\begin{lemma} \label{discrepderiv} Fix $\rho_1 \in \mathcal{M}(Z)$. Let $t \in (a, b) \mapsto \tau_t \in C(Z)$ be a
differentiable curve in $C(Z)$, and denote by $\tau'_t \in C(Z), t \in (a,b)$ its derivative. Then for any $\xi \in Z$ and
$t \in (a,b)$, the right-hand derivative $\frac{d^+}{dt} D_{\rho_1}(\tau_t)(\xi)$ exists and is given by
$$
\frac{d^+}{dt} D_{\rho_1}(\tau_t)(\xi) = \tau'_t(\xi) + \sup_{\eta \in A_t} \tau'_t(\eta)
$$
where $A_t = \{ \eta \in Z | \xi \sim_{\rho_t} \eta \}$ and $\rho_t = E_{\rho_1}(\tau_t)$.
\end{lemma}

\medskip

\noindent{\bf Proof:} Fix $\xi \in Z$ and $t \in (a,b)$.
For any $\eta \in Z$ such that $\xi \sim_{\rho_t} \eta$, for $h > 0$ small we have, by definition of $D_{\rho_1}$,
\begin{align*}
D_{\rho_1}(\tau_{t+h})(\xi) & \geq \tau_{t+h}(\xi) + \tau_{t+h}(\eta) + \log \rho_1(\xi, \eta)^2 \\
                            & = (\tau_{t+h} - \tau_t)(\xi) + (\tau_{t+h} - \tau_t)(\eta) + \log \rho_t(\xi, \eta)^2 \\
                            & = (\tau_{t+h} - \tau_t)(\xi) + (\tau_{t+h} - \tau_t)(\eta) + D_{\rho_1}(\tau_t)(\xi),
\end{align*}
from which it follows easily that
$$
\liminf_{h \to 0^+} \frac{D_{\rho_1}(\tau_{t+h})(\xi) - D_{\rho_1}(\tau_t)(\xi)}{h} \geq \tau'_t(\xi) + \tau'_t(\eta)
$$
for all $\eta$ such that $\xi \sim_{\rho_t} \eta$, and hence
\begin{equation} \label{lowerbound}
\liminf_{h \to 0^+} \frac{D_{\rho_1}(\tau_{t+h})(\xi) - D_{\rho_1}(\tau_t)(\xi)}{h} \geq \tau'_t(\xi) + \sup_{\eta \in A_t} \tau'_t(\eta).
\end{equation}

\medskip

For $h > 0$ small, let $f_h = (\tau_{t+h} - \tau_t)/h \in C(Z)$, then since $(\tau_s)_{s \in (a,b)}$ is a
differentiable curve in $C(Z)$, we have $f_h \to \tau'_t$ uniformly on $Z$ as $h \to 0^+$.
For each $h > 0$ let $\eta_h \in Z$ be such that $\xi \sim_{\rho_{t+h}} \eta_h$. Then we have
\begin{align*}
D_{\rho_1}(\tau_{t+h})(\xi) & = \tau_{t+h}(\xi) + \tau_{t+h}(\eta_h) + \log \rho_1(\xi, \eta_h)^2 \\
                              & = (\tau_{t+h} - \tau_t)(\xi) + (\tau_{t+h} - \tau_t)(\eta_h) + \log \rho_t(\xi, \eta_h)^2 \\
                              & \leq (\tau_{t+h} - \tau_t)(\xi) + (\tau_{t+h} - \tau_t)(\eta_h) + D_{\rho_1}(\tau_t)(\xi),
\end{align*}
thus
\begin{equation} \label{righthandbound}
\frac{D_{\rho_1}(\tau_{t+h})(\xi) - D_{\rho_1}(\tau_t)(\xi)}{h} \leq f_h(\xi) + f_h(\eta_h).
\end{equation}
Now as $h \to 0^+$ we have $\rho_{t+h} \to \rho_t$ uniformly on $Z \times Z$, from which it follows easily that any limit point $\eta$ of the points
$\eta_h \in A_{t+h}$ must be contained in $A_t$. Since $f_h \to \tau'_t$ uniformly as $h \to 0^+$, this implies that
$$
\limsup_{h \to 0^+} f_h(\eta_h) \leq \sup_{\eta \in A_t} \tau'_t(\eta).
$$
From this and (\ref{righthandbound}) it follows that
$$
\limsup_{h \to 0^+} \frac{D_{\rho_1}(\tau_{t+h})(\xi) - D_{\rho_1}(\tau_t)(\xi)}{h} \leq \tau'_t(\xi) + \sup_{\eta \in A_t} \tau'_t(\eta),
$$
which, together with (\ref{lowerbound}), completes the proof of the Lemma. $\diamond$

\medskip

As consequences we obtain the following:

\medskip

\begin{lemma} \label{discrepdecay} Let $\rho_1 \in \mathcal{M}(Z)$, and let $t \in [0, \infty) \mapsto \tau_t \in C(Z)$ be a
solution of the $\rho_1$-antipodal flow. Then for any $\xi \in Z$, the right-hand derivative
$\frac{d^+}{dt} D_{\rho_1}(\tau_t)(\xi)$ exists and satisfies
\begin{equation} \label{diffineq}
\frac{d^+}{dt} D_{\rho_1}(\tau_t)(\xi) \leq - 2 D_{\rho_1}(\tau_t)(\xi)
\end{equation}
for all $t \geq 0$. We also have
\begin{equation} \label{expdecay1}
D_{\rho_1}(\tau_t)(\xi) \leq D_{\rho_1}(\tau_s)(\xi) \cdot e^{-2(t-s)}
\end{equation}
for all $t \geq s \geq 0$.
\end{lemma}

\medskip

\noindent{\bf Proof:} Let $\xi \in Z$ and $t \geq 0$. Let $\rho_t = E_{\rho_1}(\tau_t)$, and let $A_t = \{ \eta \in Z | \xi \sim_{\rho_t} \eta \}$.
Since $(\tau_s)_{s \geq 0}$ is a differentiable curve in $C(Z)$ with $\tau'_t = - D_{\rho_1}(\tau_t)$, by the previous Lemma \ref{discrepderiv},
the right-hand derivative $\frac{d^+}{dt} D_{\rho_1}(\tau_t)(\xi)$ exists and is given by
$$
\frac{d^+}{dt} D_{\rho_1}(\tau_t)(\xi) = - D_{\rho_1}(\tau_t)(\xi) + \sup_{\eta \in A_t} (- D_{\rho_1}(\tau_t)(\eta)).
$$
Now for any $\eta \in A_t$, we have
\begin{align*}
D_{\rho_1}(\tau_t)(\eta) & \geq \tau_t(\eta) + \tau_t(\xi) + \log \rho_1(\eta, \xi)^2 \\
                         & = D_{\rho_1}(\tau_t)(\xi),
\end{align*}
thus
$$
\sup_{\eta \in A_t} (- D_{\rho_1}(\tau_t)(\eta)) \leq -D_{\rho_1}(\tau_t)(\xi),
$$
and hence
$$
\frac{d^+}{dt} D_{\rho_1}(\tau_t)(\xi) \leq - 2 D_{\rho_1}(\tau_t)(\xi)
$$
as required.

\medskip

From (\ref{diffineq}) we obtain
$$
\frac{d^+}{dt} \left( e^{2t} D_{\rho_1}(\tau_t)(\xi) \right) \leq 0,
$$
from which (\ref{expdecay1}) follows by applying Lemma \ref{righthandderiv}. $\diamond$

\medskip

\begin{lemma} \label{oncenegative} Let $\rho_1 \in \mathcal{M}(Z)$, and let $t \in [0, \infty) \mapsto \tau_t \in C(Z)$ be a
solution of the $\rho_1$-antipodal flow. Let $\xi \in Z$. If
$$
D_{\rho_1}(\tau_{t_0})(\xi) \leq 0
$$
for some $t_0 \geq 0$, then
$$
D_{\rho_1}(\tau_t)(\xi) \leq 0
$$
for all $t \geq t_0$.
\end{lemma}

\medskip

\noindent{\bf Proof:} Given $D_{\rho_1}(\tau_{t_0})(\xi) \leq 0$ for some $t_0 \geq 0$, then for all $t \geq t_0$, by
(\ref{expdecay1}) we have
\begin{align*}
D_{\rho_1}(\tau_t)(\xi) & \leq D_{\rho_1}(\tau_{t_0})(\xi) \cdot e^{-2(t - t_0)} \\
                        & \leq 0.
\end{align*}
$\diamond$

\medskip

\begin{lemma} \label{discderivest}  Let $\rho_1 \in \mathcal{M}(Z)$, and let $t \in [0, \infty) \mapsto \tau_t \in C(Z)$ be a
solution of the $\rho_1$-antipodal flow. Let $\xi \in Z$. If
$$
D_{\rho_1}(\tau_t)(\xi) \leq 0
$$
for some $t \geq 0$, then
\begin{equation} \label{discrepsq}
\frac{d^+}{dt} (D_{\rho_1}(\tau_t)(\xi))^2 \leq - 2 (D_{\rho_1}(\tau_t)(\xi))^2 + 2 \cdot ||D_{\rho_1}(\tau_0)||_{\infty} \cdot |D_{\rho_1}(\tau_t)(\xi)| \cdot e^{-2t}
\end{equation}
(for that same $t$).
\end{lemma}

\medskip

\noindent{\bf Proof:} Given $D_{\rho_1}(\tau_t)(\xi) \leq 0$, let $\rho_t = E_{\rho_1}(\tau_t)$ and $A_t = \{ \eta \in Z | \xi \sim_{\rho_t} \eta \}$.
Note that $A_t$ is compact, so we can choose $\eta_1 \in A_t$ such that
$$
-D_{\rho_1}(\tau_t)(\eta_1) = \sup_{\eta \in A_t} (-D_{\rho_1}(\tau_t)(\eta)).
$$
Then by Lemma \ref{discrepderiv} and the Chain Rule we have
\begin{equation} \label{sqderiv}
\frac{d^+}{dt} (D_{\rho_1}(\tau_t)(\xi))^2 = -2 D_{\rho_1}(\tau_t)(\xi) \cdot (D_{\rho_1}(\tau_t)(\xi) + D_{\rho_1}(\tau_t)(\eta_1)).
\end{equation}
On the other hand, using (\ref{expdecay1}), we have
\begin{align*}
D_{\rho_1}(\tau_t)(\xi) + D_{\rho_1}(\tau_t)(\eta_1) - e^{-2t}D_{\rho_1}(\tau_0)(\eta_1) & \leq D_{\rho_1}(\tau_t)(\xi) \\
                                                                                         & \leq 0.
\end{align*}
Since $D_{\rho_1}(\tau_t)(\xi) \leq 0$, this implies
$$
D_{\rho_1}(\tau_t)(\xi) \cdot (D_{\rho_1}(\tau_t)(\xi) + D_{\rho_1}(\tau_t)(\eta_1) - e^{-2t}D_{\rho_1}(\tau_0)(\eta_1)) \geq (D_{\rho_1}(\tau_t)(\xi))^2,
$$
and therefore
\begin{align*}
D_{\rho_1}(\tau_t)(\xi) \cdot (D_{\rho_1}(\tau_t)(\xi) + D_{\rho_1}(\tau_t)(\eta_1)) & \geq (D_{\rho_1}(\tau_t)(\xi))^2 + e^{-2t} \cdot D_{\rho_1}(\tau_0)(\eta_1) \cdot D_{\rho_1}(\tau_t)(\xi) \\
                                                                                     & \geq (D_{\rho_1}(\tau_t)(\xi))^2 - e^{-2t} \cdot ||D_{\rho_1}(\tau_0)||_{\infty} \cdot |D_{\rho_1}(\tau_t)(\xi)|,
\end{align*}
and the Lemma follows immediately from this last inequality and the equality (\ref{sqderiv}) above. $\diamond$

\medskip

We can now show that for any solution of the antipodal flow, the sup norm of the discrepancy decays exponentially as $t \to \infty$:

\medskip

\begin{prop} \label{expdecay2} Let $\rho_1 \in \mathcal{M}(Z)$, and let $t \in [0, \infty) \mapsto \tau_t \in C(Z)$ be a
solution of the $\rho_1$-antipodal flow. Then
$$
||D_{\rho_1}(\tau_t)||_{\infty} \leq 2 ||D_{\rho_1}(\tau_0)||_{\infty} \cdot e^{-t/2}
$$
for all $t \geq 0$.
\end{prop}

\medskip

\noindent{\bf Proof:} Suppose not. Then there is a $t \geq 0$ and a $\xi \in Z$ such that
\begin{equation} \label{forcontra}
|D_{\rho_1}(\tau_t)(\xi)| > 2 ||D_{\rho_1}(\tau_0)||_{\infty} \cdot e^{-t/2}.
\end{equation}
Note that the reverse inequality
$$
|D_{\rho_1}(\tau_s)(\xi)| \leq 2 ||D_{\rho_1}(\tau_0)||_{\infty} \cdot e^{-s/2}
$$
holds trivially for $s = 0$, hence we must have $t > 0$. Now if $D_{\rho_1}(\tau_t)(\xi) \geq 0$, then by Lemma \ref{discrepdecay} we have
\begin{align*}
|D_{\rho_1}(\tau_t)(\xi)| = D_{\rho_1}(\tau_t)(\xi) & \leq D_{\rho_1}(\tau_0)(\xi) \cdot e^{-2t} \\
                                                    & \leq 2 ||D_{\rho_1}(\tau_0)||_{\infty} \cdot e^{-t/2},
\end{align*}
which contradicts (\ref{forcontra}). Thus we must have $D_{\rho_1}(\tau_t)(\xi) < 0$.

\medskip

Now let
$$
t_1 = \sup \{ s \in [0,t] \ | \ |D_{\rho_1}(\tau_s)(\xi)| \leq 2 ||D_{\rho_1}(\tau_0)||_{\infty} \cdot e^{-s/2} \},
$$
and note that we have the following: $0 \leq t_1 < t$,
\begin{equation} \label{toneeq}
|D_{\rho_1}(\tau_{t_1})(\xi)| = 2 ||D_{\rho_1}(\tau_0)||_{\infty} \cdot e^{-t_1/2},
\end{equation}
and
\begin{equation} \label{sbig}
|D_{\rho_1}(\tau_s)(\xi)| > 2 ||D_{\rho_1}(\tau_0)||_{\infty} \cdot e^{-s/2}
\end{equation}
for all $s \in (t_1, t]$. Now note that we must have $D_{\rho_1}(\tau_{t_1})(\xi) \leq 0$, because otherwise if $D_{\rho_1}(\tau_{t_1})(\xi) > 0$
then, since $D_{\rho_1}(\tau_t)(\xi) < 0$, by continuity there would be an $s \in (t_1, t)$ such that $D_{\rho_1}(\tau_s)(\xi) = 0$, but this would contradict the inequality (\ref{sbig}) above.
Hence $D_{\rho_1}(\tau_{t_1})(\xi) \leq 0$, and so by Lemma \ref{oncenegative} we have
$$
D_{\rho_1}(\tau_s)(\xi) \leq 0
$$
for all $s \in [t_1, t]$. Then by Lemma \ref{discderivest}, we obtain
$$
\frac{d^+}{ds} (D_{\rho_1}(\tau_s)(\xi))^2 \leq - 2 (D_{\rho_1}(\tau_s)(\xi))^2 + 2 \cdot ||D_{\rho_1}(\tau_0)||_{\infty} \cdot |D_{\rho_1}(\tau_s)(\xi)| \cdot e^{-2s}
$$
for all $s \in [t_1, t]$. On the other hand by (\ref{sbig}) we have
$$
|D_{\rho_1}(\tau_s)(\xi)|^2 \geq 2 ||D_{\rho_1}(\tau_0)||_{\infty} \cdot |D_{\rho_1}(\tau_s)(\xi)| \cdot e^{-s/2}
$$
for $s \in [t_1,t]$, so substituting this in the previous inequality gives
\begin{align*}
\frac{d^+}{ds} (D_{\rho_1}(\tau_s)(\xi))^2 & \leq - (D_{\rho_1}(\tau_s)(\xi))^2 - (D_{\rho_1}(\tau_s)(\xi))^2 + 2 \cdot ||D_{\rho_1}(\tau_0)||_{\infty} \cdot |D_{\rho_1}(\tau_s)(\xi)| \cdot e^{-2s} \\
                                           & \leq - (D_{\rho_1}(\tau_s)(\xi))^2 - 2 ||D_{\rho_1}(\tau_0)||_{\infty} \cdot |D_{\rho_1}(\tau_s)(\xi)| \cdot e^{-s/2}
                                           + 2 \cdot ||D_{\rho_1}(\tau_0)||_{\infty} \cdot |D_{\rho_1}(\tau_s)(\xi)| \cdot e^{-2s} \\
                                           & \leq - (D_{\rho_1}(\tau_s)(\xi))^2,
\end{align*}
thus
$$
\frac{d^+}{ds} \left( e^s (D_{\rho_1}(\tau_s)(\xi))^2 \right) \leq 0
$$
for $s \in [t_1, t]$. Applying Lemma \ref{righthandderiv}, (\ref{forcontra}) and (\ref{toneeq}), we have
\begin{align*}
e^t \cdot \left(2 e^{-t/2} ||D_{\rho_1}(\tau_0)||_{\infty}\right)^2 & < e^t (D_{\rho_1}(\tau_t)(\xi))^2 \\
                                                                    & \leq e^{t_1} (D_{\rho_1}(\tau_{t_1})(\xi))^2 \\
                                                                    & = e^{t_1} \cdot \left( 2 e^{-t_1/2} ||D_{\rho_1}(\tau_0)||_{\infty} \right)^2 \\
                                                                    & = 4 \left(||D_{\rho_1}(\tau_0)||_{\infty}\right)^2,
\end{align*}
and so
$$
4 \left(||D_{\rho_1}(\tau_0)||_{\infty}\right)^2 < 4 \left(||D_{\rho_1}(\tau_0)||_{\infty}\right)^2,
$$
a contradiction. $\diamond$

\medskip

The previous Proposition allows us to prove the main Theorem of this section:

\medskip

\begin{theorem} \label{antipodallimit} {\bf (Limit of the antipodal flow)} Let $\rho_1 \in \mathcal{M}(Z)$.
Given $\tau_0 \in C(Z)$, let $t \in [0, \infty) \mapsto \tau_t \in C(Z)$ denote the unique solution of the $\rho_1$-antipodal
flow starting from $\tau_0$. Then the limit
$$
\tau_{\infty} := \lim_{t \to \infty} \tau_t
$$
exists in $C(Z)$. Moreover, we have the exponential convergence estimate
\begin{equation} \label{expconv}
|| \tau_t - \tau_{\infty} ||_{\infty} \leq 4 ||D_{\rho_1}(\tau_0)||_{\infty} \cdot e^{-t/2}
\end{equation}
for all $t \geq 0$.

\medskip

Also,
$$
D_{\rho_1}(\tau_{\infty}) \equiv 0.
$$
Thus if we put $\rho := E_{\rho_1}(\tau_{\infty})$, then $\rho$ is an antipodal
function Moebius equivalent to $\rho_1$, thus $\rho \in \mathcal{M}(Z)$, and we have
$$
\log \frac{d\rho}{d\rho_1} = \tau_{\infty}.
$$
\end{theorem}

\medskip

\noindent{\bf Proof:} By the previous Proposition \ref{expdecay2}, we have
\begin{align*}
\int_{0}^{\infty} \left|\left| \frac{d\tau_t}{dt} \right|\right|_{\infty} \ dt & = \int_{0}^{\infty} ||D_{\rho_1}(\tau_t)||_{\infty} \ dt \\
                                                                          & \leq \int_{0}^{\infty} 2 ||D_{\rho_1}(\tau_0)||_{\infty} \cdot e^{-t/2} \ dt \\
                                                                          & = 4 ||D_{\rho_1}(\tau_0)||_{\infty} \\
                                                                          & < +\infty,
\end{align*}
so since $C(Z)$ is a Banach space, absolute convergence of the improper integral implies convergence of the improper integral in $C(Z)$, i.e. the limit
$$
\lim_{t \to \infty} \int_{0}^{t} \frac{d\tau_s}{ds} \ ds
$$
exists in $C(Z)$, hence by the Fundamental Theorem of Calculus the limit
$$
\tau_{\infty} := \lim_{t \to \infty} \tau_t
$$
exists in $C(Z)$. Moreover, we have
\begin{align*}
|| \tau_t - \tau_{\infty} ||_{\infty} & \leq \int_{t}^{\infty} \left|\left| \frac{d\tau_s}{ds} \right|\right|_{\infty} \ ds \\
                                      & \leq \int_{t}^{\infty} 2 ||D_{\rho_1}(\tau_0)||_{\infty} \cdot e^{-s/2} \ ds \\
                                      & = 4 ||D_{\rho_1}(\tau_0)||_{\infty} \cdot e^{-t/2}.
\end{align*}

Note that by Proposition \ref{expdecay2}, $||D_{\rho_1}(\tau_t)||_{\infty} \to 0$ as $t \to \infty$. Since $\tau_t \to \tau_{\infty}$ in
$C(Z)$ as $t \to \infty$ and the map $D_{\rho_1} : C(Z) \to C(Z)$ is continuous, it follows that
$$
D_{\rho_1}(\tau_{\infty}) = 0.
$$
From this it follows immediately that if we put $\rho = E_{\rho_1}(\tau_{\infty})$, then $\rho$ is an antipodal function Moebius equivalent
to $\rho_1$, thus $\rho \in \mathcal{M}(Z, \rho_0)$, and moreover by definition of the map $E_{\rho_1}$ we have
$$
\log \frac{d\rho}{d\rho_1} = \tau_{\infty}.
$$
$\diamond$

\medskip

We now have all the tools necessary to construct non-trivial points of the space $\mathcal{M}(Z, \rho_0)$.

\medskip

\begin{definition} \label{antipodalization} {\bf (The antipodalization map)} Let $(Z, \rho_0)$ be an antipodal space.
The {\it antipodalization map} is the map
\begin{align*}
\mathcal{P}_{\infty} : C(Z) \times  \mathcal{M}(Z) & \to \mathcal{M}(Z) \\
                       (\tau_0 ,    \rho_1)                   & \mapsto \rho := E_{\rho_1}(\tau_{\infty}),
\end{align*}
where, for $\rho_1 \in \mathcal{M}(Z)$ and $\tau_0 \in C(Z)$, $\tau_{\infty} = \tau_{\infty}(\tau_0, \rho_1)$ denotes
the limit in $C(Z)$ as $t \to \infty$ of the unique solution $(\tau_t)_{t \geq 0}$ of the $\rho_1$-antipodal flow
starting from $\tau_0$.
\end{definition}

\medskip

We remark that for any $\rho_1, \rho_2 \in \mathcal{M}(Z)$ and $\tau_0 \in C(Z)$,
it follows immediately from the Geometric Mean-Value Theorem that
\begin{equation} \label{shiftd}
D_{\rho_2}\left(\tau_0 + \log \frac{d\rho_1}{d\rho_2}\right) = D_{\rho_1}(\tau_0) \in C(Z),
\end{equation}
and
\begin{equation} \label{shifte}
E_{\rho_2}\left(\tau_0 + \log \frac{d\rho_1}{d\rho_2}\right) = E_{\rho_1}(\tau_0) \in \mathcal{UM}(Z).
\end{equation}

These identities lead to the following invariance property of the antipodalization map with respect to the second variable:

\medskip

\begin{prop} \label{pinfinv} Let $\rho_1, \rho_2 \in \mathcal{M}(Z)$. Then for all $\tau_0 \in C(Z)$, we have
$$
\mathcal{P}_{\infty}\left( \tau_0 + \log \frac{d\rho_1}{d\rho_2}, \rho_2\right) = \mathcal{P}_{\infty}(\tau_0, \rho_1).
$$
\end{prop}

\medskip

\noindent{\bf Proof:} Given $\tau_0 \in C(Z)$, let $(\tau_t)_{t \geq 0}$ be the unique solution of the $\rho_1$-antipodal flow
starting from $\tau_0$, and let $\tau_{\infty} = \lim_{t \to \infty} \tau_t$. Let
$$
\widetilde{\tau_t} := \tau_t + \log \frac{d\rho_1}{d\rho_2}
$$
for all $t \geq 0$. Then clearly
$$
\frac{d \widetilde{\tau_t}}{dt} = \frac{d \tau_t}{dt}
$$
for all $t \geq 0$. Moreover, by (\ref{shiftd}) we have
$$
D_{\rho_2}(\widetilde{\tau_t}) = D_{\rho_1}(\tau_t)
$$
for all $t \geq 0$. It follows from the above two equalities that
$$
\frac{d \widetilde{\tau_t}}{dt} = - D_{\rho_2}(\widetilde{\tau_t})
$$
for all $t \geq 0$, and hence $(\widetilde{\tau_t})_{t \geq 0}$ is the unique solution of the
$\rho_2$-antipodal flow starting from $\widetilde{\tau_0} = \tau_0 + \log \frac{d\rho_1}{d\rho_2}$.
Since
$$
\lim_{t \to \infty} \widetilde{\tau_t} = \tau_{\infty} + \log \frac{d\rho_1}{d\rho_2},
$$
it follows that
\begin{align*}
\mathcal{P}_{\infty}\left(\tau_0 + \log \frac{d\rho_1}{d\rho_2}, \rho_2\right) & = E_{\rho_2}\left(\tau_{\infty} + \log \frac{d\rho_1}{d\rho_2}\right) \\
                                                                               & = E_{\rho_1}(\tau_{\infty}) \\
                                                                               & = \mathcal{P}_{\infty}(\tau_0, \rho_1)
\end{align*}
(where we have used (\ref{shifte}) above). $\diamond$

\medskip

The above Proposition lets us naturally define convex combinations of points in $\mathcal{M}(Z)$ as follows.
Let $\mu$ be a probability measure on $\mathcal{M}(Z)$ with compact support. Then for any $\rho_1 \in \mathcal{M}(Z)$, we have a
continuous function $\int_{\mathcal{M}(Z)} \log \frac{d\rho}{d\rho_1} \ d\mu(\rho) \in C(Z)$, defined by
$$
\left( \int_{\mathcal{M}(Z)} \log \frac{d\rho}{d\rho_1} \ d\mu(\rho) \right) (\xi) := \int_{\mathcal{M}(Z)} \log \frac{d\rho}{d\rho_1}(\xi) \ d\mu(\rho) \ , \xi \in Z.
$$
The continuity of the above function follows from an application of Lebesgue's Dominated Convergence Theorem, after noting that if $\xi_n \to \xi$ in $Z$
then $\log \frac{d\rho}{d\rho_1}(\xi_n) \to \log \frac{d\rho}{d\rho_1}(\xi)$ for all $\rho \in \mathcal{M}(Z)$, and
$$
\int_{\mathcal{M}(Z)} \left|\left| \log \frac{d\rho}{d\rho_1} \right|\right|_{\infty} \ d\mu(\rho) < +\infty
$$
since $\mu$ has compact support and is a finite measure. The fact that $\mu$ is a probability measure implies that for any $\rho_2 \in \mathcal{M}(Z)$ we have
$$
\int_{\mathcal{M}(Z)} \log \frac{d\rho}{d\rho_2} \ d\mu(\rho) = \int_{\mathcal{M}(Z)} \log \frac{d\rho}{d\rho_1} \ d\mu(\rho) + \log \frac{d\rho_1}{d\rho_2}.
$$
It follows from this and Proposition \ref{pinfinv} above that
\begin{equation} \label{convexeq}
\mathcal{P}_{\infty}\left( \int_{\mathcal{M}(Z)} \log \frac{d\rho}{d\rho_2} \ d\mu(\rho), \rho_2\right) = \mathcal{P}_{\infty}\left( \int_{\mathcal{M}(Z)} \log \frac{d\rho}{d\rho_1} \ d\mu(\rho), \rho_1\right).
\end{equation}
This leads naturally to the following definition:

\medskip

\begin{definition} \label{convexcomb} {\bf (Convex combinations in $\mathcal{M}(Z)$)} Let $\mu$ be a probability measure
on $\mathcal{M}(Z)$ with compact support. Then we define the $\mu$-convex combination to be the point
$$
\int_{\mathcal{M}(Z)} \rho \ d\mu(\rho) \in \mathcal{M}(Z)
$$
defined by
$$
\int_{\mathcal{M}(Z)} \rho \ d\mu(\rho) := \mathcal{P}_{\infty}\left( \int_{\mathcal{M}(Z)} \log \frac{d\rho}{d\rho_1} \ d\mu(\rho), \rho_1\right)
$$
where $\rho_1$ is any element of $\mathcal{M}(Z)$. Note that by (\ref{convexeq}), this is well-defined independent of the
choice of $\rho_1$.
\end{definition}

\medskip

In the next section we will see how these convex combinations can be used to construct geodesics in $\mathcal{M}(Z)$.

\medskip

We have continuity of the map $\mathcal{P}_{\infty}$ in the first variable:

\medskip

\begin{prop} \label{pinfcont} Fix $\rho_1 \in \mathcal{M}(Z)$. Then the map
\begin{align*}
\mathcal{P}_{\infty}(., \rho_1) : C(Z) & \to \mathcal{M}(Z) \\
                                  \tau_0 & \mapsto \mathcal{P}_{\infty}(\tau_0, \rho_1)
\end{align*}
is continuous.
%
%
\end{prop}

\medskip

\noindent{\bf Proof:} Let $\tau_0 \in C(Z)$ and let $\epsilon > 0$.

\medskip

Choose $T > 0$ such that
$$
4(||D_{\rho_1}(\tau_0)||_{\infty} + 2) e^{-T/2} < \epsilon/3.
$$
For $\tau \in C(Z)$, let $\gamma(.,\tau) : [0, \infty) \to C(Z)$ denote the unique solution of the $\rho_1$-antipodal flow
starting from $\tau$. Since the vector field $X = - D_{\rho_1}$ on $C(Z)$ is globally Lipschitz, we have continuous dependence of solutions
on initial conditions, thus for each $t > 0$ the time-$t$ map
$$
\tau \in C(Z) \mapsto \gamma(t, \tau) \in C(Z)
$$
is continuous. Thus we can choose $\delta \in (0,1)$ such that for all $\tau \in B(\tau_0, \delta)$ we have
$$
||\gamma(T, \tau) - \gamma(T, \tau_0)||_{\infty} < \epsilon/3.
$$
Also, since $\delta < 1$ and $D_{\rho_1}$ is $2$-Lipschitz, we have
$$
||D_{\rho_1}(\tau)||_{\infty} \leq ||D_{\rho_1}(\tau_0)||_{\infty} + 2
$$
for all $\tau \in B(\tau_0, \delta)$.

\medskip

Denoting by $\gamma(\infty, \tau)$ the limit as $t \to \infty$ of $\gamma(t, \tau)$, it follows from the estimate (\ref{expconv}) of
Theorem \ref{antipodallimit} that
\begin{align*}
||\gamma(\infty, \tau_0) - \gamma(T, \tau_0)||_{\infty} & \leq 4 ||D_{\rho_1}(\tau_0)||_{\infty} \cdot e^{-T/2} \\
                                                        & < \epsilon/3,
\end{align*}
and
\begin{align*}
||\gamma(\infty, \tau) - \gamma(T, \tau)||_{\infty} & \leq 4 (||D_{\rho_1}(\tau_0)||_{\infty} + 2) \cdot e^{-T/2} \\
                                                        & < \epsilon/3
\end{align*}
for all $\tau \in B(\tau_0, \delta)$. Thus we obtain
\begin{align*}
||\gamma(\infty, \tau) - \gamma(\infty, \tau_0)||_{\infty} & \leq ||\gamma(\infty, \tau) - \gamma(T, \tau)||_{\infty} + ||\gamma(T, \tau) - \gamma(T, \tau_0)||_{\infty}
                                                                + ||\gamma(T, \tau_0) - \gamma(\infty, \tau_0)||_{\infty} \\
                                                           & < \epsilon/3 + \epsilon/3 + \epsilon/3 = \epsilon
\end{align*}
for all $\tau \in B(\tau_0, \delta)$, which proves continuity of the map $\mathcal{P}_{\infty}(.,\rho_1)$,
once we observe that
$$
d_{\mathcal{M}}( \mathcal{P}_{\infty}(\tau, \rho_1), \mathcal{P}_{\infty}(\tau_0, \rho_1) ) = ||\gamma(\infty, \tau) - \gamma(\infty, \tau_0)||_{\infty}.
$$
$\diamond$

\medskip

The antipodalization map in fact defines a retraction of the space $\mathcal{UM}(Z)$ onto the closed subspace
$\mathcal{M}(Z) \subset \mathcal{UM}(Z)$:

\medskip

\begin{prop} \label{retraction} The map
\begin{align*}
R_{\infty} : \mathcal{UM}(Z) & \to \mathcal{M}(Z) \\
         \rho       & \mapsto \mathcal{P}_{\infty}\left(\log \frac{d\rho}{d\rho_1}, \rho_1\right),
\end{align*}
where $\rho_1$ is any element of $\mathcal{M}(Z)$, is well-defined independent of the choice of $\rho_1$,
and is a retraction of $\mathcal{UM}(Z)$ onto the closed subspace
$\mathcal{M}(Z) \subset \mathcal{UM}(Z)$.
\end{prop}

\medskip

\noindent{\bf Proof:} The fact that the map $R_{\infty}$ is well-defined independent of the choice of $\rho_1$ is a special
case of the fact that $\mu$-convex combinations are well-defined (Definition \ref{convexcomb}), namely the case when $\mu$ is
a Dirac mass at a point $\rho \in \mathcal{M}(Z)$. By the previous Proposition \ref{pinfcont},
the map $R_{\infty}$ is continuous.

\medskip

Now fix a $\rho_1 \in \mathcal{M}(Z)$. Given any $\rho \in \mathcal{M}(Z)$,
let $\tau_0 = \log \frac{d\rho}{d\rho_1}$, and let $(\tau_t)_{t \geq 0}$
denote the unique solution of the $\rho_1$-antipodal flow starting from $\tau_0$, and let $\tau_{\infty} = \lim_{t \to \infty} \tau_t$.
Since $\rho$ is an antipodal function we have
$$
D_{\rho_1}(\tau_0) = 0,
$$
and hence, by estimate (\ref{expconv}) we have
$$
||\tau_{\infty} - \tau_0||_{\infty} \leq 4 ||D_{\rho_1}(\tau_0)||_{\infty} \cdot e^{-0/2} = 0,
$$
thus $\tau_{\infty} = \tau_0$, and so
\begin{align*}
R_{\infty}(\rho) & = \mathcal{P}_{\infty}(\tau_0, \rho_1) \\
                 & = E_{\rho_1}(\tau_{\infty}) \\
                 & = E_{\rho_1}(\tau_0) \\
                 & = \rho,
\end{align*}
which shows that $R_{\infty}$ is a retraction $\diamond$

\medskip

In fact, it is not hard to show that if we denote by $\mathcal{P}_t(.,\rho_1) : C(Z) \to \mathcal{UM}(Z)$ the map
given by composing the time-$t$ map of the $\rho_1$-antipodal flow with the map $E_{\rho_1} : C(Z) \to \mathcal{UM}(Z)$,
then the family of maps
\begin{align*}
R_t : \mathcal{UM}(Z) & \to \mathcal{UM}(Z) \\
          \rho        & \mapsto \mathcal{P}_t\left(\log \frac{d\rho}{d\rho_1}, \rho_1\right),
\end{align*}
$0 \leq t \leq \infty$, define a deformation retraction of $\mathcal{UM}(Z)$ onto the closed subspace $\mathcal{M}(Z) \subset \mathcal{UM}(Z)$.

\medskip

It also follows immediately from the above Proposition, that conjugating the map $R_{\infty} : \mathcal{UM}(Z) \to \mathcal{UM}(Z)$
by the isometry $\phi_{\rho_1} : \mathcal{UM}(Z) \to C(Z)$ gives a map $\widetilde{R_{\rho_1}} : C(Z) \to C(Z)$ which is a retraction of
$C(Z)$ onto the closed subspace $\{ \tau \ | \ D_{\rho_1}(\tau) = 0 \} \subset C(Z)$.
\medskip

The antipodalization map also gives contractibility of the space $\mathcal{M}(Z)$:

\medskip

\begin{prop} \label{contractible} The space $\mathcal{M}(Z)$ is contractible. Explicitly, for any $\rho_1 \in \mathcal{M}(Z)$, the map
\begin{align*}
H : \mathcal{M}(Z) \times [0,1] & \to \mathcal{M}(Z) \\
         (\rho, t)              & \mapsto \mathcal{P}_{\infty}\left( t \log \frac{d\rho}{d\rho_1}, \rho_1\right)
\end{align*}
gives a homotopy from the identity to the constant map $c : \rho \mapsto \rho_1$.
\end{prop}

\medskip

\noindent{\bf Proof:} By Proposition \ref{pinfcont} the map $H$ is continuous. For
$\rho \in \mathcal{M}(Z)$, using Proposition \ref{retraction} we get
\begin{align*}
H(\rho, 1) & = \mathcal{P}_{\infty}\left( \log \frac{d\rho}{d\rho_1}, \rho_1\right) \\
           & = \rho,
\end{align*}
and, since $\frac{d\rho_1}{d\rho_1} \equiv 1$,
\begin{align*}
H(\rho, 0) & = \mathcal{P}_{\infty}( 0, \rho_1) \\
           & = \mathcal{P}_{\infty}\left( \log \frac{d\rho_1}{d\rho_1}, \rho_1\right) \\
           & = \rho_1.
\end{align*}
$\diamond$

\medskip

The following estimates on the map $\mathcal{P}_{\infty}$ will be essential:

\medskip

\begin{prop} \label{pinfest} Let $\rho_1 \in \mathcal{M}(Z)$. Let $\tau_0 \in C(Z)$, let $\rho = \mathcal{P}_{\infty}(\tau_0, \rho_1)$
and let $\tau = \log \frac{d\rho}{d\rho_1}$. Then for all $\xi \in Z$,
we have
\begin{equation} \label{lwrbound}
\tau(\xi) \geq \tau_0(\xi) - \frac{1}{2} D_{\rho_1}(\tau_0)(\xi),
\end{equation}
and
\begin{equation} \label{uprbound}
\tau(\xi) \leq \tau_0(\xi) - D_{\rho_1}(\tau_0)(\xi) + \frac{1}{2} D_{\rho_1}(\tau_0)(\eta)
\end{equation}
where $\eta \in Z$ is any point such that $\xi \sim_{\rho_0} \eta$, where $\rho_0 = E_{\rho_1}(\tau_0)$.
\end{prop}

\medskip

\noindent{\bf Proof:} Note that by definition $\tau = \tau_{\infty} = \lim_{t \to \infty} \tau_t$, where
$(\tau_t)_{t \geq 0}$ is the unique solution of the $\rho_1$-antipodal flow starting from $\tau_0$.
Using the estimate (\ref{expdecay1}), we obtain
\begin{align*}
\tau(\xi) & = \tau_0(\xi) + \int_{0}^{\infty} \frac{d\tau_t}{dt}(\xi) \ dt \\
          & = \tau_0(\xi) - \int_{0}^{\infty} D_{\rho_1}(\tau_t)(\xi) \ dt \\
          & \geq \tau_0(\xi) - \int_{0}^{\infty} D_{\rho_1}(\tau_0)(\xi) \cdot e^{-2t} \ dt \\
          & = \tau_0(\xi) - \frac{1}{2} D_{\rho_1}(\tau_0)(\xi),
\end{align*}
which proves the lower bound (\ref{lwrbound}).

\medskip

Now let $\eta \in Z$ be such that $\xi \sim_{\rho_0} \eta$. Then by definition
$$
D_{\rho_1}(\tau_0)(\xi) = \tau_0(\xi) + \tau_0(\eta) + \log \rho_1(\xi, \eta)^2.
$$
Now, using the fact that $\rho$ is of diameter one, the Geometric Mean-Value Theorem, the above equality, and the estimate
(\ref{lwrbound}) applied at the point $\eta$, we get
\begin{align*}
0 & \geq \log \rho(\xi, \eta)^2 \\
  & = \tau(\xi) + \tau(\eta) + \log \rho_1(\xi, \eta)^2 \\
  & = \tau(\xi) + \tau(\eta) + D_{\rho_1}(\tau_0)(\xi) - \tau_0(\xi) - \tau_0(\eta) \\
  & \geq \tau(\xi) + \left(\tau_0(\eta) - \frac{1}{2} D_{\rho_1}(\tau_0)(\eta)\right) + D_{\rho_1}(\tau_0)(\xi) - \tau_0(\xi) - \tau_0(\eta) \\
  & = \tau(\xi) - \frac{1}{2} D_{\rho_1}(\tau_0)(\eta) + D_{\rho_1}(\tau_0)(\xi) - \tau_0(\xi),
\end{align*}
and rearranging terms above gives the inequality
$$
\tau(\xi) \leq \tau_0(\xi) - D_{\rho_1}(\tau_0)(\xi) + \frac{1}{2} D_{\rho_1}(\tau_0)(\eta)
$$
as required. $\diamond$

\medskip

\section{Geodesics in $\mathcal{M}(Z)$ and the maxima compactification $\widehat{\mathcal{M}(Z)}$.}

\medskip

Let $Z$ be an antipodal space and let $\mathcal{M}(Z)$ be the associated Moebius space. In this section we
will show that $\mathcal{M}(Z)$ is a geodesic metric space which is geodesically complete. We will
also describe a natural compactification $\widehat{\mathcal{M}(Z)} = \mathcal{M}(Z) \sqcup Z$ of $\mathcal{M}(Z)$ in which the
boundary of $\mathcal{M}(Z)$ will be the space $Z$. For any $\rho \in \mathcal{M}(Z)$ and $\xi \in Z$, we will
construct a geodesic ray $\gamma : [0, \infty) \to \mathcal{M}(Z)$ such that $\gamma(0) = \rho$ and $\gamma(t) \to \xi$ in
$\widehat{\mathcal{M}(Z)}$ as $t \to \infty$. We will also show how points of $\mathcal{M}(Z)$ can be
reconstructed from limits of Gromov inner products in $\mathcal{M}(Z)$ as one tends to points of $Z$ in
$\widehat{\mathcal{M}(Z)}$.

\medskip

\begin{lemma} \label{tcomb} Let $\rho_1, \rho_2 \in \mathcal{M}(Z)$ and $t \in [0,1]$. Define
$$
\rho := (1-t) \rho_1 + t \rho_2 \in \mathcal{M}(Z),
$$
which is the $\mu$-convex combination corresponding to the probability measure $\mu = (1-t) \delta_{\rho_1} + t \delta_{\rho_2}$
(see Definition \ref{convexcomb}). Then
$$
d_{\mathcal{M}}( \rho_1, \rho ) =  t d_{\mathcal{M}}( \rho_1, \rho_2 ) \quad \hbox{and} \quad d_{\mathcal{M}}( \rho_2, \rho ) =  (1-t) d_{\mathcal{M}}( \rho_1, \rho_2 ).
$$
\end{lemma}

\medskip

\noindent{\bf Proof:} According to Definition \ref{convexcomb}, we can write
\begin{align*}
\rho & = \mathcal{P}_{\infty}\left( (1-t) \log \frac{d\rho_1}{d\rho_1} + t \log \frac{d\rho_2}{d\rho_1}, \rho_1\right) \\
     & = \mathcal{P}_{\infty}\left( t \log \frac{d\rho_2}{d\rho_1}, \rho_1\right).
\end{align*}
Let
$$
\tau = \log \frac{d\rho}{d\rho_1}.
$$
Now for any $\xi \in Z$ and $\eta \in Z - \{\xi\}$, by the Geometric Mean-Value Theorem we have
\begin{align*}
& t \log \frac{d\rho_2}{d\rho_1}(\xi) + t \log \frac{d\rho_2}{d\rho_1}(\eta) + \log \rho_1(\xi, \eta)^2 \\
& =  t \log \frac{d\rho_2}{d\rho_1}(\xi) + t \log \frac{d\rho_2}{d\rho_1}(\eta) + t \log \rho_1(\xi, \eta)^2 + (1-t) \log \rho_1(\xi, \eta)^2 \\
& = t \log \rho_2(\xi, \eta)^2 + (1-t) \log \rho_1(\xi, \eta)^2 \\
& \leq 0,
\end{align*}
thus
$$
D_{\rho_1}\left( t \log \frac{d\rho_2}{d\rho_1} \right)(\xi) \leq 0.
$$
It then follows from (\ref{lwrbound}) of Proposition \ref{pinfest} and Lemma \ref{maxmin} that
\begin{align*}
\tau(\xi) & \geq t \log \frac{d\rho_2}{d\rho_1}(\xi) - \frac{1}{2} D_{\rho_1}\left( t \log \frac{d\rho_2}{d\rho_1} \right)(\xi) \\
          & \geq t \log \frac{d\rho_2}{d\rho_1}(\xi) \\
          & \geq - t \left|\left| \log \frac{d\rho_2}{d\rho_1} \right|\right|_{\infty} \\
          & = - t d_{\mathcal{M}}( \rho_1, \rho_2)
\end{align*}
for all $\xi \in Z$, thus
$$
\min_{\xi \in Z} \tau(\xi) \geq - t d_{\mathcal{M}}( \rho_1, \rho_2),
$$
and so by Lemma \ref{maxmin}
$$
d_{\mathcal{M}}(\rho_1, \rho) = - \min_{\xi \in Z} \tau(\xi) \leq t d_{\mathcal{M}}( \rho_1, \rho_2).
$$
On the other hand, if we choose $\xi_0 \in \hbox{argmax} \frac{d\rho_2}{d\rho_1}$ then as above we have, using (\ref{lwrbound}) and Lemma \ref{maxmin}, that
\begin{align*}
\tau(\xi_0) & \geq t \log \frac{d\rho_2}{d\rho_1}(\xi_0) - \frac{1}{2} D_{\rho_1}\left( t \log \frac{d\rho_2}{d\rho_1} \right)(\xi_0) \\
            & \geq t \log \frac{d\rho_2}{d\rho_1}(\xi_0) \\
          & = t \left|\left| \log \frac{d\rho_2}{d\rho_1} \right|\right|_{\infty} \\
          & = t d_{\mathcal{M}}( \rho_1, \rho_2),
\end{align*}
and so
$$
d_{\mathcal{M}}(\rho_1, \rho) \geq \tau(\xi_0) \geq t d_{\mathcal{M}}( \rho_1, \rho_2),
$$
thus
$$
d_{\mathcal{M}}(\rho_1, \rho) = t d_{\mathcal{M}}( \rho_1, \rho_2)
$$
as required. A similar proof gives $d_{\mathcal{M}}(\rho_2, \rho) = (1-t) d_{\mathcal{M}}( \rho_1, \rho_2)$. $\diamond$

\medskip

As an immediate corollary we obtain:

\medskip

\begin{prop} \label{spacegeodesic} The space $\mathcal{M}(Z)$ is a geodesic metric space.
\end{prop}

\medskip

\noindent{\bf Proof:} By the previous Lemma \ref{tcomb}, for any $\rho_1, \rho_2 \in \mathcal{M}(Z)$ there exists
$\rho \in \mathcal{M}(Z)$ such that
$$
d_{\mathcal{M}}(\rho_1, \rho) = d_{\mathcal{M}}(\rho_2, \rho) = \frac{1}{2} d_{\mathcal{M}}( \rho_1, \rho_2),
$$
i.e. any pair of points in $\mathcal{M}(Z)$ has a midpoint. As is well known, the existence of midpoints
together with the completeness of $\mathcal{M}(Z)$ (Lemma \ref{properness}) imply that $\mathcal{M}(Z)$
is a geodesic metric space. $\diamond$

\medskip

\medskip

\begin{lemma} \label{maxdiam} Let $\rho_1, \rho_2 \in \mathcal{M}(Z)$. Then for any $R \in \mathbb{R}$,
$$
\hbox{diam}_{\rho_1}\left( \left\{ \xi \in Z \ | \ \log \frac{d\rho_2}{d\rho_1}(\xi) \geq R \right\} \right) \leq e^{- R}.
$$
\end{lemma}

\medskip

\noindent{\bf Proof:} For any $\xi, \eta$ such that
$$
\log \frac{d\rho_2}{d\rho_1}(\xi), \log \frac{d\rho_2}{d\rho_1}(\eta) \geq R,
$$
by the
Geometric Mean-Value Theorem we have
\begin{align*}
1 & \geq \rho_2(\xi, \eta)^2 \\
  & = \frac{d\rho_2}{d\rho_1}(\xi) \cdot \frac{d\rho_2}{d\rho_1}(\eta) \cdot \rho_1(\xi, \eta)^2 \\
  & \geq e^{R} \cdot e^{R} \cdot \rho_1(\xi, \eta)^2,
\end{align*}
thus $\rho_1(\xi, \eta) \leq e^{-R}$ as required. $\diamond$

\medskip

We now want to define a compactification of $\mathcal{M}(Z)$ in which $Z$ appears as the boundary of $\mathcal{M}(Z)$.
We let $\widehat{\mathcal{M}(Z)}$ be the set
$$
\widehat{\mathcal{M}(Z)} := \mathcal{M}(Z) \sqcup Z.
$$
Motivated by the previous Lemma, we define the following basic neighbourhoods in $\widehat{\mathcal{M}(Z)}$ of points in $Z$:

\medskip

\begin{definition} \label{basicnbhd} Let $U \subset Z$ be an open set, let $\rho_1 \in \mathcal{M}(Z)$, and let $r > 0$. We define
the set $N(U, \rho_1, r) \subset \widehat{\mathcal{M}(Z)}$ by
$$
N(U, \rho_1, r) := U \sqcup \{ \rho \in \mathcal{M}(Z) \ | \ \hbox{argmax} \ \frac{d\rho}{d\rho_1} \subset U, \ d_{\mathcal{M}}(\rho, \rho_1) > r \}.
$$
\end{definition}

\medskip

\begin{prop} \label{basistop} Let $\rho_1 \in \mathcal{M}(Z)$. The collection of sets
$$
\mathcal{B}(\rho_1) := \{ B(\rho, r) \ | \rho \in \mathcal{M}(Z), r > 0 \} \cup \{ N(U, \rho_1, r) \ | \ U \subset Z, U \ \hbox{open}, r > 0 \}
$$
form a basis for a topology $\mathcal{T}(\rho_1)$ on $\widehat{\mathcal{M}(Z)}$. With this topology the space
$\widehat{\mathcal{M}(Z)}$ is a compact, Hausdorff space.
\end{prop}

\medskip

\noindent{\bf Proof:} Clearly the given collection of sets covers $\widehat{\mathcal{M}(Z)}$. We need to check the intersection
property for the collection to be a basis for a topology. For a point in the intersection of two balls this is clear.

\medskip

Note that any set $N(U, \rho_1,r) \cap \mathcal{M}(Z)$ is open in $\mathcal{M}(Z)$; this follows from the elementary fact that if $f$ is a
continuous function on a compact space and $U$ is an open set such that $\hbox{argmax} \ f \subset U$, then $\hbox{argmax} \ g \subset U$
for all $g$ uniformly close enough to $f$.

\medskip

Now for a point $\rho \in B(\rho_1, r_1) \cap N(U, \rho_1, r_2)$, since
$N(U, \rho_1, r_2) \cap \mathcal{M}(Z)$ is open in $\mathcal{M}(Z)$, we can find a ball
$B(\rho, \epsilon) \subset B(\rho_1, r_1) \cap N(U, \rho_1, r_2)$.

\medskip

For a point $z \in N(U_1, \rho_1, r_1) \cap N(U_2, \rho_1, r_2)$, if $z = \rho \in \mathcal{M}(Z)$ then as before
 we can find a ball $B(\rho, \epsilon) \subset N(U_1, \rho_1, r_1) \cap N(U_2, \rho_1, r_2)$.

\medskip

The only remaining case is when $z = \xi \in Z$. Then $\xi \in U_1 \cap U_2$. Letting $r = \max(r_1, r_2)$,
it is clear that $N(U_1 \cap U_2, \rho_1, r) \subset N(U_1, \rho_1, r_1) \cap N(U_2, \rho_1, r_2)$, and
$\xi \in N(U_1 \cap U_2, \rho_1, r)$.

\medskip

Thus the collection $\mathcal{B}(\rho_1)$ forms a basis for a topology $\mathcal{T}(\rho_1)$ on $\widehat{\mathcal{M}(Z)}$.
To prove compactness of $\widehat{\mathcal{M}(Z)}$ with this topology, let $\mathcal{V} = \{ V_j \}_{j \in J}$ be an open cover of $\widehat{\mathcal{M}(Z)}$.
For each $\xi \in Z$, there exists a $V_{\xi} \in \mathcal{V}$ containing $\xi$, hence there is a basic open set
$N(U_{\xi}, \rho_1,  r_{\xi}) \subset V_{\xi}$ containing $\xi$. Then $\{ U_{\xi} \ | \ \xi \in Z \}$ is an open
cover of $Z$, hence has a finite subcover $\{ U_{\xi_i} \ | 1 \leq i \leq n \}$, for some $\xi_1, \dots, \xi_n \in Z$. Let $V_i = V_{\xi_i}$ be
the corresponding elements of $\mathcal{V}$, similarly let $U_i = U_{\xi_i}, r_i = r_{\xi_i}, 1 \leq i \leq n$.

\medskip

Recall that $Z$ is metrizable, let $\alpha$ be a metric on $Z$ compatible with its topology.
 Let $\delta_0 > 0$ be the Lebesgue number with respect to the metric $\alpha$ of the finite cover
$\{ U_i \ | 1 \leq i \leq n \}$ of $Z$. Then it is easy to show, using positivity and continuity of $\rho_1$,
that there is an $\epsilon > 0$ such that $\rho_1(\xi, \eta) < \epsilon$ implies $\alpha(\xi, \eta) < \delta_0$. Choose $r > \max \{ r_i \ | \ 1 \leq i \leq n \}$
such that $e^{-r} < \epsilon$.

\medskip

Then for any $\rho \in \mathcal{M}(Z)$ such that $d_{\mathcal{M}}(\rho, \rho_1) > r$, by Lemma \ref{maxdiam}
we have
$$
\hbox{diam}_{\rho_1}\left( \hbox{argmax} \ \frac{d\rho}{d\rho_1} \right) \leq e^{-r} < \epsilon,
$$
hence by our choice of $\epsilon$ the $\alpha$-diameter of $\hbox{argmax} \ \frac{d\rho}{d\rho_1}$ is less than $\delta_0$,
so there exists an $i$ such that $\hbox{argmax} \ \frac{d\rho}{d\rho_1} \subset U_i$, and hence
$\rho \in N(U_i, \rho_1, r_i) \subset V_i$. This shows that
$$
Z \sqcup \{ \rho \in \mathcal{M}(Z) \ | \ d_{\mathcal{M}}(\rho, \rho_1) > r \} \subset \cup_{i = 1}^n V_i.
$$
On the other hand the closed ball $B := \{ \rho \in \mathcal{M}(Z) \ | \ d_{\mathcal{M}}(\rho, \rho_1) \leq r \}$ is compact by Lemma \ref{properness}, hence can be
covered by finitely many elements of the cover $\mathcal{V}$, while as we have just seen the complement $\widehat{\mathcal{M}(Z)} - B$ is covered by
$V_1, \dots, V_n$, thus the cover $\mathcal{V}$ has a finite subcover covering $\widehat{\mathcal{M}(Z)}$. This proves compactness of
$\widehat{\mathcal{M}(Z)}$.

\medskip

Finally, to show that $\widehat{\mathcal{M}(Z)}$ is Hausdorff, it is enough to separate distinct points of $Z$ in $\widehat{\mathcal{M}(Z)}$.
Given distinct points $\xi, \eta \in Z$, choose disjoint open neighbourhoods $U_1, U_2$ of $\xi, \eta$ respectively, choose
any $r > 0$, then clearly $N(U_1, \rho_1, r)$ and $N(U_2, \rho_1, r)$ are disjoint open neighbourhoods
of $\xi, \eta$ respectively in $\widehat{\mathcal{M}(Z)}$. $\diamond$

\medskip

\begin{lemma} \label{topind} The topology $\mathcal{T}(\rho_1)$ on $\widehat{\mathcal{M}(Z)}$ is independent of the
choice of $\rho_1$, i.e. $\mathcal{T}(\rho_1) = \mathcal{T}(\rho_2)$ for all $\rho_1, \rho_2 \in \mathcal{M}(Z)$.
\end{lemma}

\medskip

\noindent{\bf Proof:} Let $\rho_1, \rho_2 \in \mathcal{M}(Z)$.
It is enough to show that given $z \in N(U_1, \rho_1, r_1)$, there exists $V \in \mathcal{T}(\rho_2)$
such that $z \in V$ and $V \subset N(U_1, \rho_1, r_1)$.

\medskip

If $z = \rho \in \mathcal{M}(Z)$, then we can take $V$ to be a ball $B(\rho, \epsilon) \in \mathcal{T}(\rho_2)$
such that $B(\rho, \epsilon) \subset N(U_1, \rho_1, r_1)$.

\medskip

Otherwise $z = \xi \in U_1 \subset Z$. Then let $U_2$ be an open neighbourhood of $\xi$ in $Z$ such that
$\overline{U_2} \subset U_1$. Since $\overline{U_2}$ and $Z - U_1$ are disjoint compact sets, it is not hard to show
using positivity and continuity of $\rho_2$ that there is an $\epsilon > 0$ such that if $\eta \in \overline{U_2}$ and
$\zeta$ is such that $\rho_2(\eta, \zeta) < \epsilon$, then $\zeta \in U_1$. Now let $M = d_{\mathcal{M}}(\rho_1, \rho_2)$,
and choose $r_2 > r_1 + M$ such that $e^{-(r_2 - 2M)} < \epsilon$.

\medskip

Then given $\rho \in N(U_2, \rho_2, r_2) \cap \mathcal{M}(Z)$, choose $\eta \in \hbox{argmax} \ \frac{d\rho}{d\rho_2} \subset U_2$.
Then $\log \frac{d\rho}{d\rho_2}(\eta) = d_{\mathcal{M}}( \rho, \rho_2) > r_2$. Now for any $\zeta \in \hbox{argmax} \ \frac{d\rho}{d\rho_1}$
we have
\begin{align*}
 \log \frac{d\rho}{d\rho_2}(\zeta) & = \log \frac{d\rho}{d\rho_1}(\zeta) + \log \frac{d\rho_1}{d\rho_2}(\zeta) \\
                                   & \geq d_{\mathcal{M}}( \rho, \rho_1) - M \\
                                   & \geq (r_2 - M) - M = r_2 - 2M,
\end{align*}
so it follows from Lemma \ref{maxdiam} that
$$
\rho_2( \eta, \zeta ) \leq e^{-(r_2 - 2M)} < \epsilon,
$$
and by our choice of $\epsilon$ this implies that $\zeta \in U_1$.
It follows that $\hbox{argmax} \ \frac{d\rho}{d\rho_1} \subset U_1$, also $d_{\mathcal{M}}( \rho, \rho_1) \geq r_2 - M > r_1$,
hence $\rho \in N(U_1, \rho_1, r_1)$.

\medskip

So we have $N(U_2, \rho_2, r_2) \in \mathcal{T}(\rho_2)$ such that $N(U_2, \rho_2, r_2) \subset N(U_1, \rho_1, r_1)$ and $\xi \in N(U_2, \rho_2, r_2)$
as required. $\diamond$

\medskip

We remark that it is not hard to show that the space $\widehat{\mathcal{M}(Z)}$ with this topology is in fact second countable,
and hence metrizable.

\medskip

\begin{lemma} \label{raypt} Let $\rho_1 \in \mathcal{M}(Z)$ and let $\xi_0 \in Z$. Then for any $r > 0$, there exists $\rho_2 \in \mathcal{M}(Z)$
such that
$$
d_{\mathcal{M}}(\rho_1, \rho_2) = r, \quad \hbox{and} \quad \xi_0 \in \hbox{argmax} \ \frac{d\rho_2}{d\rho_1}.
$$
\end{lemma}

\medskip

\noindent{\bf Proof:} Let $U$ be a neighbourhood of $\xi_0$ in $Z$ such that $\rho_1(\xi, \eta) < e^{-2r}$ for all $\xi, \eta \in U$.
Let $\tau_0 \in C(Z)$ be a continuous function
with support contained in $U$, such that $0 \leq \tau_0 \leq 2r$ on $U$, and such that $\tau_0(\xi_0) = 2r$.
Let $\rho_2 = \mathcal{P}_{\infty}(\tau_0, \rho_1)$, and let $\tau = \log \frac{d\rho_2}{d\rho_1}$.

\medskip

Let $\xi \in U$. Then for any $\eta \in U - \{\xi\}$, we have
$$
\tau_0(\xi) + \tau_0(\eta) + \log \rho_1(\xi, \eta)^2 < 2r + 2r + (-4r) = 0,
$$
while for $\eta \in Z - U$ we have
$$
\tau_0(\xi) + \tau_0(\eta) + \log \rho_1(\xi, \eta)^2 \leq 2r + 0 + 0 = 2r,
$$
hence
$$
D_{\rho_1}(\tau_0)(\xi) \leq 2r
$$
for all $\xi \in U$. On the other hand, if $\xi \in Z - U$, then for all $\eta \in Z - \{\xi\}$ we have
$$
\tau_0(\xi) + \tau_0(\eta) + \log \rho_1(\xi, \eta)^2 \leq 0 + 2r + 0 = 2r,
$$
and hence
$$
D_{\rho_1}(\tau_0)(\xi) \leq 2r
$$
also for all $\xi \in Z - U$. It follows that for all $\xi \in Z$, by estimate (\ref{lwrbound}) of Proposition \ref{pinfest}, we have
$$
\tau(\xi) \geq \tau_0(\xi) - \frac{1}{2} D_{\rho_1}(\tau_0)(\xi) \geq 0 - \frac{1}{2} \cdot 2r = -r,
$$
thus by Lemma \ref{maxmin}
$$
\max_{\xi \in Z} \tau(\xi) = - \min_{\xi \in Z} \tau(\xi) \leq r.
$$
On the other hand, by estimate (\ref{lwrbound})
$$
\tau(\xi_0) \geq \tau_0(\xi_0) - \frac{1}{2} D_{\rho_1}(\tau_0)(\xi_0) \geq 2r - \frac{1}{2} \cdot 2r = r,
$$
thus
$$
\max_{\xi \in Z} \tau(\xi) = \tau(\xi_0) = r,
$$
hence
$$
d_{\mathcal{M}}(\rho_1, \rho_2) = r, \quad \hbox{and} \quad \xi_0 \in \hbox{argmax} \ \frac{d\rho_2}{d\rho_1}.
$$
as required. $\diamond$

\medskip

\begin{lemma} \label{maximaconv} Let $\xi \in Z$ and $\rho \in \mathcal{M}(Z)$. Suppose $\rho_n \in \mathcal{M}(Z)$
and $\xi_n \in Z$ are sequences such that $\xi_n \in \hbox{argmax} \ \frac{d\rho_n}{d\rho}$ for all $n$, and such that
$\xi_n \to \xi$ and $d_{\mathcal{M}}( \rho, \rho_n) \to \infty$ as $n \to \infty$. Then $\rho_n \to \xi$ in
$\widehat{\mathcal{M}(Z)}$.
\end{lemma}

\medskip

\noindent{\bf Proof:} It is enough, given a basic open set $N(U, \rho, r)$ containing $\xi$, to
show that $\rho_n \in N(U, \rho, r)$ for all $n$ large enough.
Since $U$ is an open neighbourhood of $\xi$, by positivity and continuity of $\rho$ it is not hard to show that
we can choose $\epsilon > 0$ such that for any $\xi', \eta \in Z$, if $\rho(\xi, \xi') < \epsilon$
and $\rho(\xi', \eta) < \epsilon$, then $\eta \in U$.

\medskip

Let $r_n = d_{\mathcal{M}}( \rho, \rho_n)$.
Then there exists $N \geq 1$ such that for all $n \geq N$, we have $\rho(\xi, \xi_n) < \epsilon$,
$r_n > r$, and $e^{-r_n} < \epsilon$. Then it follows from Lemma \ref{maxdiam} that for all $\eta \in \hbox{argmax} \ \frac{d\rho_n}{d\rho}$,
we have $\rho(\xi_n, \eta) \leq e^{-r_n} < \epsilon$, hence $\eta \in U$. Thus $\hbox{argmax} \ \frac{d\rho_n}{d\rho} \subset U$ and so
$\rho_n \in N(U, \rho, r)$ for all $n \geq N$. $\diamond$

\medskip

\begin{prop} \label{opendense} $\mathcal{M}(Z)$ is an open, dense subset of $\widehat{\mathcal{M}(Z)}$.
\end{prop}

\medskip

\noindent{\bf Proof:} It is clear from the definition of the topology on $\widehat{\mathcal{M}(Z)}$ that $\mathcal{M}(Z)$ is
open in $\widehat{\mathcal{M}(Z)}$. Density of $\mathcal{M}(Z)$ follows immediately from the previous two Lemmas,
since given $\xi \in Z$, then fixing a $\rho \in \mathcal{M}(Z)$
by Lemma \ref{raypt} we can construct a sequence $\rho_n \in \mathcal{M}(Z)$ such that
$\xi \in \hbox{argmax} \ \frac{d\rho_n}{d\rho}$ for all $n$ and such that $d_{\mathcal{M}}( \rho, \rho_n) \to \infty$ as $n \to \infty$.
Then by Lemma \ref{maximaconv} we have $\rho_n \to \xi$ in $\widehat{\mathcal{M}(Z)}$ as $n \to \infty$. $\diamond$

\medskip

\begin{definition} \label{extcomp} {\bf (Maxima compactification)} The {\it maxima compactification} of $\mathcal{M}(Z)$ is
defined to be the space $\widehat{\mathcal{M}(Z)} = \mathcal{M}(Z) \sqcup Z$ equipped with the topology $\mathcal{T}(\rho_1)$
for any $\rho_1 \in \mathcal{M}(Z)$ (which is independent of the choice of $\rho_1$ by Lemma \ref{topind}).
\end{definition}

\medskip

Note that by Propositions \ref{basistop} and \ref{opendense}, the maxima compactification $\widehat{\mathcal{M}(Z)}$ is
indeed a compactification of $\mathcal{M}(Z)$.

\medskip

Fix a $\rho_1 \in \mathcal{M}(Z)$. Then it is easy to see that a sequence $\rho_n \in \mathcal{M}(Z)$
converges in $\widehat{\mathcal{M}(Z)}$ to a point $\xi \in Z$ if and only if
for any neighbourhood $U$ of $\xi$ in $Z$, we have $\hbox{argmax} \ \frac{d\rho_n}{d\rho_1} \subset U$ for all $n$ large,
and $\max \frac{d\rho_n}{d\rho_1}\to +\infty$ as $n \to \infty$. In other words, the maximum value of the derivatives $\frac{d\rho_n}{d\rho_1}$
tends to $+\infty$, and the set of maxima of $\frac{d\rho_n}{d\rho_1}$ becomes localized around the point $\xi \in Z$.

\medskip

We now proceed to show that $\mathcal{M}(Z)$ is geodesically complete.

\medskip

\begin{lemma} \label{trieq} Let $\rho_1, \rho_2, \rho_3 \in \mathcal{M}(Z)$ be such that
$$
d_{\mathcal{M}}(\rho_1, \rho_3) = d_{\mathcal{M}}(\rho_1, \rho_2) + d_{\mathcal{M}}(\rho_2, \rho_3).
$$
Then
$$
\hbox{argmax} \ \frac{d\rho_3}{d\rho_1} \subset \hbox{argmax} \ \frac{d\rho_2}{d\rho_1} \cap \hbox{argmax} \ \frac{d\rho_3}{d\rho_2}.
$$
\end{lemma}

\medskip

\noindent{\bf Proof:} Let $\xi \in \hbox{argmax} \ \frac{d\rho_3}{d\rho_1}$, then
\begin{align*}
d_{\mathcal{M}}(\rho_1, \rho_3) & = \log \frac{d\rho_3}{d\rho_1}(\xi) \\
                                & = \log \frac{d\rho_2}{d\rho_1}(\xi) + \log \frac{d\rho_3}{d\rho_2}(\xi) \\
                                & \leq \max_{\zeta \in Z} \log \frac{d\rho_2}{d\rho_1}(\zeta) + \max_{\zeta \in Z} \log \frac{d\rho_3}{d\rho_2}(\zeta) \\
                                & = d_{\mathcal{M}}(\rho_1, \rho_2) + d_{\mathcal{M}}(\rho_2, \rho_3) \\
                                & = d_{\mathcal{M}}(\rho_1, \rho_3),
\end{align*}
so we must have equality throughout, which implies
$$
\xi \in \hbox{argmax} \ \frac{d\rho_2}{d\rho_1} \cap \hbox{argmax} \ \frac{d\rho_3}{d\rho_2}.
$$
$\diamond$

\medskip

\begin{lemma} \label{concatenate} Let $\gamma_1 : [a, b] \to \mathcal{M}(Z)$ and $\gamma_2 : [b,c] \to \mathcal{M}(Z)$ be
geodesics such that $\gamma_1(b) = \gamma_2(b)$. Let $\gamma = \gamma_1 \cdot \gamma_2 : [a,c] \to \mathcal{M}(Z)$ denote
the path given by concatenating $\gamma_1$ and $\gamma_2$ (so $\gamma_{|[a,b]} = \gamma_1$ and $\gamma_{|[b,c]} = \gamma_2$).
Let $\rho_1 = \gamma_1(a), \rho_2 = \gamma_1(b)$ and $\rho_3 = \gamma_2(c)$. Suppose that either
$$
\hbox{argmax} \ \frac{d\rho_3}{d\rho_2} \cap \hbox{argmax} \ \frac{d\rho_2}{d\rho_1} \neq \emptyset
$$
or
$$
\hbox{argmax} \ \frac{d\rho_1}{d\rho_2} \cap \hbox{argmax} \ \frac{d\rho_2}{d\rho_3} \neq \emptyset,
$$
then $\gamma$ is a geodesic.
\end{lemma}

\medskip

\noindent{\bf Proof:} Suppose there exists a $\xi_0 \in \hbox{argmax} \ \frac{d\rho_3}{d\rho_2} \cap \hbox{argmax} \ \frac{d\rho_2}{d\rho_1}$.
To show $\gamma$ is a geodesic, it is enough to show that $d_{\mathcal{M}}(\gamma_1(s), \gamma_2(t)) = t - s$ for all
$s \in [a,b], t \in [b,c]$.

\medskip

Given $s \in [a,b], t \in [b,c]$, let $\alpha = \gamma_1(s), \beta = \gamma_2(t)$. Since $\gamma_1$ and $\gamma_2$ are geodesics,
we have $d_{\mathcal{M}}( \rho_1, \rho_2) = d_{\mathcal{M}}( \rho_1, \alpha) + d_{\mathcal{M}}( \alpha, \rho_2)$ and
$d_{\mathcal{M}}( \rho_2, \rho_3) = d_{\mathcal{M}}( \rho_2, \beta) + d_{\mathcal{M}}( \beta, \rho_3)$.
It follows from Lemma \ref{trieq} that $\hbox{argmax} \ \frac{d\rho_2}{d\rho_1} \subset \hbox{argmax} \ \frac{d\rho_2}{d\alpha}$
and $\hbox{argmax} \ \frac{d\rho_3}{d\rho_2} \subset \hbox{argmax} \ \frac{d\beta}{d\rho_2}$,
thus
$$
\xi_0 \in \hbox{argmax} \ \frac{d\beta}{d\rho_2} \cap \hbox{argmax} \ \frac{d\rho_2}{d\alpha},
$$
which gives
\begin{align*}
d_{\mathcal{M}}( \alpha, \beta) & \geq \log \frac{d\beta}{d\alpha}(\xi_0) \\
                                & = \log \frac{d\beta}{d\rho_2}(\xi_0) + \log \frac{d\rho_2}{d\alpha}(\xi_0) \\
                                & = d_{\mathcal{M}}( \beta, \rho_2) + d_{\mathcal{M}}( \rho_2, \alpha) \\
                                & = (t - b) + (b - s) = t - s,
\end{align*}
while $d_{\mathcal{M}}( \alpha, \beta) \leq (t - b) + (b - s) = t - s$ by the triangle inequality, thus
$$
d_{\mathcal{M}}( \alpha, \beta) = t - s
$$
as required, showing that $\gamma$ is a geodesic. The other case where $\hbox{argmax} \ \frac{d\rho_1}{d\rho_2} \cap \hbox{argmax} \ \frac{d\rho_2}{d\rho_3} \neq \emptyset$
is handled similarly. $\diamond$

\medskip

As a corollary of the above Lemma, we obtain:

\medskip

\begin{prop} \label{geodcomp} The space $\mathcal{M}(Z)$ is geodesically complete. More precisely,
let $\gamma_0 : [a, b] \to \mathcal{M}(Z)$ be a geodesic segment with endpoints $\rho_0 = \gamma_0(a), \rho_1 = \gamma_0(b)$,
let $\xi_0 \in \hbox{argmax} \ \frac{d\rho_1}{d\rho_0}$ and let $\eta_0 \in Z$ be such that $\rho_0(\xi_0, \eta_0) = 1$. Then
there exists a bi-infinite geodesic $\kappa : \mathbb{R} \to \mathcal{M}(Z)$ extending $\gamma_0$ such that
$\kappa(t) \to \xi_0$ in $\widehat{\mathcal{M}(Z)}$ as $t \to +\infty$, and $\kappa(t) \to \eta_0$ in $\widehat{\mathcal{M}(Z)}$ as $t \to -\infty$.
In fact, $\xi_0 \in \hbox{argmax} \ \frac{d\kappa(t)}{d\rho_0}$ for all $t > a$, and $\eta_0 \in \hbox{argmax} \ \frac{d\kappa(t)}{d\rho_0}$ for all $t < a$.
\end{prop}

\medskip

\noindent{\bf Proof:} Let $\gamma_0 : [a, b] \to \mathcal{M}(Z)$ be a given geodesic segment. Let $\rho_0 = \gamma_0(a), \rho_1 = \gamma_0(b)$,
and fix a $\xi_0 \in \hbox{argmax} \ \frac{d\rho_1}{d\rho_0}$ and let $\eta_0 \in Z$ be such that $\rho_0(\xi_0, \eta_0) = 1$.
Then note by Lemma \ref{maxmin} that $\eta_0 \in \hbox{argmin} \ \frac{d\rho_1}{d\rho_0} = \hbox{argmax} \ \frac{d\rho_0}{d\rho_1}$.

\medskip

By Lemma \ref{raypt}, we can find a point $\rho_2 \in \mathcal{M}(Z)$ such that $d_{\mathcal{M}}(\rho_1, \rho_2) = 1$ and such that
$\xi_0 \in \hbox{argmax} \ \frac{d\rho_2}{d\rho_1}$. Let $\gamma_1 : [b,b+1] \to \mathcal{M}(Z)$ be a geodesic such that $\gamma_1(b) = \rho_1,
\gamma_1(b+1) = \rho_2$. Since $\xi_0 \in \hbox{argmax} \ \frac{d\rho_2}{d\rho_1} \cap \hbox{argmax} \ \frac{d\rho_1}{d\rho_0}$, it follows from the
previous Lemma \ref{concatenate} that the concatenation $\gamma_0 \cdot \gamma_1 : [a, b+1] \to \mathcal{M}(Z)$ is a geodesic. From this
and $\xi_0 \in \hbox{argmax} \ \frac{d\rho_2}{d\rho_1} \cap \hbox{argmax} \ \frac{d\rho_1}{d\rho_0}$ we have also
\begin{align*}
\log \frac{d\rho_2}{d\rho_0}(\xi_0) & = \log \frac{d\rho_2}{d\rho_1}(\xi_0) + \log \frac{d\rho_1}{d\rho_0}(\xi_0) \\
                                    & = d_{\mathcal{M}}( \rho_2, \rho_1) + d_{\mathcal{M}}( \rho_1, \rho_0) \\
                                    & = d_{\mathcal{M}}(\rho_2, \rho_0),
\end{align*}
thus $\xi_0 \in \hbox{argmax} \ \frac{d\rho_2}{d\rho_0}$ as well. Since $\rho_0(\xi_0, \eta_0) = 1$, by Lemma \ref{maxmin} this implies
$\eta_0 \in \hbox{argmin} \ \frac{d\rho_2}{d\rho_0} = \hbox{argmax} \ \frac{d\rho_0}{d\rho_2}$, and $\rho_2(\xi_0, \eta_0) = 1$.

\medskip

By Lemma \ref{raypt}, we can find a point $\rho_{-1} \in \mathcal{M}(Z)$ such that $d_{\mathcal{M}}(\rho_0, \rho_{-1}) = 1$ and such that
$\eta_0 \in \hbox{argmax} \ \frac{d\rho_{-1}}{d\rho_0}$. Let $\gamma_{-1} : [a-1,a] \to \mathcal{M}(Z)$ be a geodesic such that $\gamma_{-1}(a-1) = \rho_{-1},
\gamma_{-1}(a) = \rho_0$. Since $\eta_0 \in \hbox{argmax} \ \frac{d\rho_{-1}}{d\rho_0} \cap \hbox{argmax} \ \frac{d\rho_0}{d\rho_2}$, it follows from
Lemma \ref{concatenate} that the concatenation $\kappa_1 := \gamma_{-1} \cdot \gamma_0 \cdot \gamma_1 : [a-1, b+1] \to \mathcal{M}(Z)$ is a geodesic. Also,
from a similar argument as above, we can show $\eta_0 \in \hbox{argmax} \ \frac{d\rho_{-1}}{d\rho_2}$. Since $\rho_2(\xi_0, \eta_0) = 1$, by
Lemma \ref{maxmin} this implies $\xi_0 \in \hbox{argmin} \ \frac{d\rho_{-1}}{d\rho_2} = \hbox{argmax} \ \frac{d\rho_2}{d\rho_{-1}}$, and
$\rho_{-1}(\xi_0, \eta_0) = 1$.

\medskip

Let $\kappa_0 = \gamma_0$. Then the above shows that given the geodesic $\kappa_0 : [a,b] \to \mathcal{M}(Z)$ with endpoints $\rho_0, \rho_1$, and
given $\xi_0 \in \hbox{argmax} \ \frac{d\rho_1}{d\rho_0}$ and $\eta_0 \in Z$ such that $\rho_0(\xi_0, \eta_0) = 1$, we can extend the geodesic
$\kappa_0$ to a geodesic $\kappa_1 : [a-1, b+1] \to \mathcal{M}(Z)$ with endpoints $\rho_{-1}, \rho_2$ which satisfy
$\xi_0 \in \hbox{argmax} \ \frac{d\rho_2}{d\rho_{-1}}$ and $\rho_{-1}(\xi_0, \eta_0) = 1$. Now it is clear that
we can apply this same procedure repeatedly to get a sequence of geodesics $\kappa_n : [a - n, b+n] \to \mathcal{M}(Z)$, such that each geodesic
$\kappa_n$ extends the previous geodesic $\kappa_{n-1}$, and such that the endpoints of $\kappa_n$ are points $\rho_{-n}, \rho_{n+1}$ which
satisfy $\xi_0 \in \hbox{argmax} \ \frac{d\rho_{n+1}}{d\rho_{-n}}$ and $\rho_{-n}(\xi_0, \eta_0) = 1$. Since each $\kappa_n$ extends $\kappa_{n-1}$, the geodesics
$\kappa_n$ define a bi-infinite geodesic $\kappa : \mathbb{R} \to \mathcal{M}(Z)$ such that $\kappa_{|[a-n,b+n]} = \kappa_n$ for
all $n$. In particular, $\kappa$ is an extension of $\gamma_0$.

\medskip

Note that $d_{\mathcal{M}}( \kappa(a-n), \kappa(b+n)) = d_{\mathcal{M}}( \kappa(a-n), \kappa(a)) + d_{\mathcal{M}}( \kappa(a), \kappa(b+n))$,
and so by Lemma \ref{trieq} the fact that $\xi_0 \in \hbox{argmax} \ \frac{d\kappa(b+n)}{d\kappa(a-n)}$ implies
that $\xi_0 \in \hbox{argmax} \ \frac{d\kappa(b+n)}{d\kappa(a)}$. For any $t > a$, we can choose an integer $n$ such that $b+n > t$, then
$d_{\mathcal{M}}( \kappa(a), \kappa(b+n)) = d_{\mathcal{M}}( \kappa(a), \kappa(t)) + d_{\mathcal{M}}( \kappa(t), \kappa(b+n))$
and $\xi_0 \in \hbox{argmax} \ \frac{d\kappa(b+n)}{d\kappa(a)}$ implies again by Lemma \ref{trieq} that $\xi_0 \in \hbox{argmax} \ \frac{d\kappa(t)}{d\kappa(a)}$
as required. Similarly we can show that $\eta_0 \in \hbox{argmax} \ \frac{d\kappa(t)}{d\kappa(a)}$ for all $t < a$.
It follows from Lemma \ref{maximaconv} that $\kappa(t) \to \xi_0$ in $\widehat{\mathcal{M}(Z)}$ as $t \to +\infty$
and $\kappa(t) \to \eta_0$ in $\widehat{\mathcal{M}(Z)}$ as $t \to -\infty$. $\diamond$

\medskip

We can now give the proof of Theorem \ref{mainthm1} from the Introduction:

\medskip

\noindent{\bf Proof of Theorem \ref{mainthm1}:} Let $Z$ be an antipodal space. 
That $\MM(Z)$ is unbounded follows from Lemma \ref{raypt}, the contractibility 
of $\MM(Z)$ is proved in Proposition \ref{contractible}, that $\MM(Z)$ is geodesic is proved in Proposition \ref{spacegeodesic}, while the geodesic completeness of $\MM(Z)$ is proved in Proposition \ref{geodcomp} above. 

\medskip

\section{Gromov product spaces and visual antipodal functions.}

\medskip

In this section we introduce the notion of {\it Gromov product spaces}, and establish the basic properties of such spaces,
which are analogous to the same properties for proper, geodesically complete CAT(-1) spaces.

\medskip

\begin{definition} \label{gromovproductspace} {\bf (Gromov product space)} Let $X$ be a proper, geodesic and geodesically complete metric space. We say that $X$
is a {\it Gromov product space} if there exists a second countable, Hausdorff compactification $\widehat{X}$ of $X$ such that for any $x \in X$,
the Gromov product $(.|.)_x : X \times X \to [0, +\infty)$ extends to a continuous function $(.|.)_x : \widehat{X} \times \widehat{X} \to [0, +\infty]$
with the property that for $\xi, \eta \in \widehat{X} \setminus X$, we have $(\xi|\eta)_x = +\infty$ if and only if $\xi = \eta$.

\medskip

We call such a compactification a {\it Gromov product compactification} of $X$.
\end{definition}

\medskip

For example, any proper, geodesic, geodesically complete, boundary continuous Gromov hyperbolic space is a Gromov product space, with the visual compactification being a
Gromov product compactification. In particular, proper, geodesically complete CAT(-1) spaces are Gromov product spaces.

\medskip

Note that by hypothesis, any Gromov product compactification is a compact metrizable space. We show that if a Gromov product compactification exists, then
it is in fact unique up to equivalence:

\medskip

\begin{prop} \label{bdryuniq} Let $X$ be a Gromov product space. Let $Y_1, Y_2$ be two Gromov product compactifications
of $X$. 

\medskip

Then the identity map $id : X \to X$ extends to a homeomorphism $F : Y_1 \to Y_2$.

\medskip

The boundary map $f : Y_1 \setminus X \to Y_2 \setminus X$ of the homeomorphism $F$ preserves the Gromov inner products on the boundaries,
$$
(f(\xi) | f(\eta))_x = (\xi | \eta)_x
$$
for all $\xi, \eta \in Y_1 \setminus X$, $x \in X$. 

\end{prop}

\medskip

\noindent{\bf Proof:} Since $Y_1, Y_2$ are metrizable, it suffices to show that if a sequence $x_n \in X$ converges in $Y_1$ to some $\xi \in Y_1 \setminus X$,
then $x_n$ converges in $Y_2$ as well, for which it is enough to show that the sequence has a unique limit point in $Y_2$.

\medskip

Thus let $\eta, \eta' \in Y_2 \setminus X$ be limit points in $Y_2$ of the sequence $x_n$. Fix a basepoint $x \in X$. Since $Y_1$ is a
Gromov product compactification and $x_n \to \xi$ in $Y_1$, we have $(x_n | x_m)_{x_0} \to (\xi | \xi)_{x} = +\infty$ as $m,n \to \infty$.
On the other hand, taking subsequences $\{y_k\}, \{z_k\} \subset \{x_n\}$ which converge in $Y_2$ to $\eta$ and $\eta'$ respectively,
we have $(y_k | z_k)_{x} \to (\eta | \eta')_x$ as $k \to \infty$, hence $(\eta | \eta')_x = +\infty$, and so $\eta = \eta'$.

\medskip

This shows that the identity map extends to a homeomorphism $F : Y_1 \to Y_2$. 

\medskip

The second statement of the Proposition follows from the fact that the identity map preserves Gromov inner products and by passing to limits as points in $X$ converge to points in the boundary $Y_1 \setminus X$. 
$\diamond$

\medskip

By virtue of the above Proposition, any Gromov product space $X$ has a unique (up to equivalence) compactification $\widehat{X}^P$
satisfying the hypotheses of Definition \ref{gromovproductspace}, which we call the Gromov product compactification of $X$. We call the set
$\partial_{P} X := \widehat{X}^{P} \setminus X$ the Gromov product boundary of $X$.

\medskip

\begin{lemma} \label{geodesicconverges} Let $X$ be a Gromov product space. Then any geodesic ray $\gamma : [0, \infty) \to X$
converges to a point in the Gromov product boundary, i.e. there exists $\xi \in \partial_{P} X$ such that $\gamma(t) \to \xi$ in
$\widehat{X}^{P}$ as $t \to \infty$.
\end{lemma}

\medskip

\noindent{\bf Proof:} Let $\gamma$ be a geodesic ray, and let $\xi, \xi' \in \partial_{P} X$ be any two limit points in $\partial_{P} X$
of the set $\{ \gamma(t) \ | \ t \geq 0 \}$. letting $x = \gamma(0)$, since $\gamma$ is a geodesic, we have
$$
(\gamma(s) | \gamma(t))_x = \min(s,t) \to +\infty \ \hbox{as} \ s,t \to +\infty,
$$
so choosing sequences $s_n, t_n \in [0, \infty)$ such that $\gamma(s_n) \to \xi, \gamma(t_n) \to \xi'$ as $n \to \infty$, we obtain
$$
(\xi | \xi')_x = \lim_{n \to \infty} (\gamma(s_n) | \gamma(t_n))_x = +\infty,
$$
thus $\xi = \xi'$, and it follows that $\gamma(t) \to \xi$ as $n \to \infty$. $\diamond$

\medskip

For any geodesic ray $\gamma$, we will denote its limit in $\partial_{P} X$ by $\gamma(\infty)$.

\medskip

\begin{lemma} \label{geodesicexists} Let $X$ be a Gromov product space. Then for any $x \in X$ and $\xi \in \partial_{P} X$, there exists a
geodesic ray $\gamma$ such that $\gamma(0) = x$ and $\gamma(\infty) = \xi$.
\end{lemma}

\medskip

\noindent{\bf Proof:} Let $x_n \in X$ be a sequence such that $x_n \to \xi$ as $n \to \infty$. For each $n$, let $\gamma_n$ be a geodesic segment joining
$x$ to $x_n$. By the Arzela-Ascoli Theorem, after passing to a subsequence we may assume that the geodesic segments $\gamma_n$ converge uniformly on compacts
in $[0, \infty)$ to a geodesic ray $\gamma$ such that $\gamma(0) = x$. Let $\xi' = \gamma(\infty)$.

\medskip

By monotonicity of the Gromov product along geodesics starting from $x$, for any $t \geq 0$ and $n$ large such that $d(x, x_n) \geq t$,
we have $(x_n | \gamma(t))_x \geq (\gamma_n(t) | \gamma(t))_x$, and hence
\begin{align*}
(\xi | \gamma(t))_x & = \lim_{n \to \infty} (x_n | \gamma(t))_x \\
                    & \geq \lim_{n \to \infty} (\gamma_n(t) | \gamma(t))_x \\
                    & = (\gamma(t) | \gamma(t))_x \\
                    & = t,
\end{align*}
thus
$$
(\xi | \xi')_x = \lim_{t \to \infty} (\xi | \gamma(t))_x = +\infty.
$$
and so $\xi = \xi'$ as required. $\diamond$

\medskip

\begin{prop} \label{visualantipodal} Let $X$ be a Gromov product space. Then for any $x \in X$, the function
$\rho_x := e^{-(.|.)_x} : \partial_{P} X \times \partial_{P} X \to [0,1]$ is an antipodal function
on the Gromov product boundary $\partial_{P} X$.
\end{prop}

\medskip

\noindent{\bf Proof:} Clearly the function $\rho_x$ is continuous and satisfies symmetry and positivity.
Given $\xi \in \partial_{P} X$, by the previous Lemma we can choose a geodesic ray $\gamma$ such that $\gamma(0) = x$
and $\gamma(\infty) = \xi$. Since $X$ is geodesically complete, we can extend the ray $\gamma$ to a bi-infinite geodesic
$\gamma : (-\infty, \infty) \to X$. Let $\eta = \gamma(-\infty) \in \partial_{P} X$. Since $\gamma$ is a geodesic passing through $x$, we
have $(\xi | \eta)_x = \lim_{t \to \infty} ( \gamma(t) | \gamma(-t))_x = 0$, hence $\rho_x(\xi, \eta) = 1$. $\diamond$

\medskip

We call the functions $\rho_x, x \in X$, {\it visual antipodal functions}.

\medskip

\begin{lemma} \label{visualmoebeqv} Let $X$ be a Gromov product space such that $\partial_{P} X$ has at least four points.
Then for any $x, y \in X$, the visual antipodal functions $\rho_x, \rho_y$ on $\partial_{P} X$ are Moebius equivalent.
\end{lemma}

\medskip

\noindent{\bf Proof:} The proof is the same as for CAT(-1) spaces. For any $x \in X$, the cross-ratio $[.,.,.,.]_{\rho_x}$ can be written as
\begin{align*}
[\xi, \xi', \eta, \eta']_{\rho_x} & = \lim_{a \to \xi, a' \to \xi', b \to \eta, b' \to \eta'} \frac{ \exp\left(-(a|b)_x\right) \exp\left(-(a'|b')_x\right) }
                                                                                                   { \exp\left(-(a|b')_x\right) \exp\left(-(a'|b)_x\right) } \\
                                  & = \lim_{a \to \xi, a' \to \xi', b \to \eta, b' \to \eta'} \exp\left( -(a|b)_x -(a'|b')_x + (a|b')_x + (a'|b)_x \right) \\
                                  & = \lim_{a \to \xi, a' \to \xi', b \to \eta, b' \to \eta'} \exp\left( \frac{1}{2} \left(d(a,b) + d(a',b') - d(a, b') - d(a', b)\right) \right),
\end{align*}
and the last expression above is clearly independent of the choice of $x$. $\diamond$

\medskip

Thus there is a canonical cross-ratio on the Gromov product boundary coming from the visual antipodal functions, and so
the Gromov product boundary $\partial_{P} X$ is an antipodal space.

\medskip

\begin{lemma} \label{isometrytomoebius} Let $X, Y$ be Gromov product spaces, and let $F : X \to Y$ be an isometric embedding. Then $F$ extends
to a Moebius embedding $f : \partial_{P} X \to \partial_{P} Y$ of antipodal spaces. If $F$ is an isometry, then the Moebius map $f$ is a
homeomorphism.
\end{lemma}

\medskip

\noindent{\bf Proof:} Let $\xi \in \partial_{P} X$ and let $x_n \in X$ be a sequence converging to $\xi$. Let $\eta, \eta' \in \partial_{P} Y$
be limit points in $\partial_{P} Y$ of the sequence $F(x_n) \in Y$. Let $\{y_k\}, \{z_k\}$ be subsequences of $\{x_n\}$ such that $F(y_k) \to \eta,
F(z_k) \to \eta'$ as $k \to \infty$. Fixing a basepoint $x \in X$, since $F$ is an isometric embedding we have
\begin{align*}
(\eta | \eta')_{F(x)} & = \lim_{k \to \infty} ( F(y_k) | F(z_k))_{F(x)} \\
                      & = \lim_{k \to \infty} ( y_k | z_k)_x \\
                      & = (\xi | \xi)_x = +\infty,
\end{align*}
hence $\eta = \eta'$. It follows that $F$ has a continuous extension to the boundary $f : \partial_{P} X \to \partial_{P} Y$, which, from an argument similar to that
above, satisfies $(f(\xi) | f(\eta))_{F(x)} = (\xi | \eta)_x$ for all $\xi, \eta \in \partial_{P} X$, in particular $f(\xi) = f(\eta)$ implies
$(\xi | \eta)_x = +\infty$ and so $\xi = \eta$, thus $f$ is one-to-one. Moreover $\rho_{F(x)}(f(\xi), f(\eta)) = \rho_x(\xi, \eta)$ for all
$\xi, \eta \in \partial_{P}X$, from which it follows that $f$ is Moebius.

\medskip

Suppose $F$ is an isometry. For any $\eta \in \partial_{P} Y$, we can choose a sequence $y_n \in Y$ such that $y_n \to \eta$ as $n \to \infty$.
Since $F$ is bijective, we can let $x_n := F^{-1}(y_n) \in X$, and then passing to a subsequence we may assume there exists $\xi \in \partial_{P} X$
such that $x_n \to \xi$ as $n \to \infty$. Then we have $f(\xi) = \lim_{n \to \infty} F(x_n) = \eta$, hence the map $f$ is surjective. Since $f$
is then a continuous bijection between compact Hausdorff spaces, it is a homeomorphism. $\diamond$

\medskip

\begin{lemma} \label{busemanngp} Let $X$ be a Gromov product space. Then for any $x, y \in X$ and $\xi \in \partial_{P} X$, we have
$$
\lim_{z \to \xi} (d(x,z) - d(y,z)) = -d(x, y) + 2(\xi|y)_x = d(x,y) - 2(\xi|x)_y.
$$
\end{lemma}

\medskip

\noindent{\bf Proof:} For any $x,y,z \in X$, we have
$$
d(x,z) - d(y,z) = -d(x, y) + 2(z|y)_x = d(x,y) - 2(z|x)_y.
$$
Since $X$ is a Gromov product space, the limits of the last two expressions above exist as $z$ tends to $\xi \in \partial_{P} X$
and the Lemma follows. $\diamond$

\medskip

We can thus define, analogously to  the case of a CAT(-1) space, the {\it Busemann cocycle} of a Gromov product space to be the
function $B : X \times X \times \partial_{P} X \to \mathbb{R}$ given by
$$
B(x, y, \xi) := \lim_{z \to \xi} (d(x,z)- d(y,z)) \ , \ x,y \in X, \xi \in \partial_{P} X.
$$

We then have the same formula for the derivatives of the visual antipodal functions as in the case of a CAT(-1) space:

\medskip

\begin{prop} \label{visualderiv} Let $X$ be a Gromov product space. For any $x, y \in X$ and $\xi \in \partial_{P} X$, we have
$$
\frac{d\rho_y}{d\rho_x}(\xi) = e^{B(x,y,\xi)}.
$$
\end{prop}

\medskip

\noindent{\bf Proof:} For any $x,y,z,w \in X$, we have
$$
\left(\frac{e^{-(z|w)_y}}{e^{-(z|w)_x}}\right)^2 = e^{d(x,z) - d(y,z)} \cdot e^{d(x,w) - d(y,z)}.
$$
Letting $z,w$ tend to $\xi, \eta \in \partial_{P} X$ respectively gives
$$
\left(\frac{\rho_y(\xi, \eta)}{\rho_x(\xi, \eta)}\right)^2 = e^{B(x,y,\xi)} \cdot e^{B(x,y,\eta)}
$$
for all $\xi, \eta \in \partial_{P} X$, so it follows from Lemma \ref{derivatives} (the Geometric Mean Value Theorem)
that
$$
\frac{d\rho_y}{d\rho_x}(\xi) = e^{B(x,y,\xi)}.
$$
$\diamond$

\medskip

\medskip

From Lemma \ref{busemanngp} it follows that $|B(x,y,\xi)| \leq d(x,y)$ for all $x,y \in X, \xi \in \partial_{P} X$.
It follows also from the definition of the Busemann cocycle that $B(x, y, \xi) = d(x, y)$ if $\xi \in \partial_{P} X$ is the endpoint
at infinity of a geodesic ray starting from $x$ and passing through $y$. Thus we have
$$
\xi \in \hbox{argmax} \ \frac{d\rho_y}{d\rho_x}
$$
if $\xi \in \partial_{P} X$ is the endpoint of a geodesic ray starting from $x$ and passing through $y$.
As in \cite{biswas3}, we have the following isometric embedding of a Gromov product space $X$ into the Moebius space $\mathcal{M}(\partial_{P}X)$:

\medskip

\begin{prop} \label{isomembed} Let $X$ be a Gromov product space.
Then the map
\begin{align*}
i_X : X & \to \mathcal{M}(\partial_{P} X) \\
      x & \mapsto \rho_x
\end{align*}
is an isometric embedding. Moreover the map $i_X : X \to \mathcal{M}(\partial_{P} X)$ extends to
a continuous map $\widehat{i_X} : \widehat{X}^{P} = X \sqcup \partial_{P} X \to \widehat{\mathcal{M}(\partial_{P} X)} = \mathcal{M}(\partial_{P} X) \sqcup \partial_{P} X$
such that $\widehat{i_X}_{|\partial_{P} X} = id : \partial_{P} X \to \partial_{P} X$.
\end{prop}

\medskip

\noindent{\bf Proof:} Let $x, y \in X$. Then for any $\xi \in \partial_{P} X$, as noted above $|B(x,y,\xi)| \leq d(x,y)$. It follows from Lemma \ref{visualderiv} above that
$$
d_{\mathcal{M}}( \rho_x, \rho_y ) \leq d(x, y).
$$
On the other hand, let $\gamma_0$ be a geodesic segment joining $x$ to $y$, then since $X$ is geodesically complete we can
extend $\gamma_0$ to a geodesic ray $\gamma$ starting from $x$ and passing through $y$, then as noted above $B(x, y, \xi) = d(x, y)$,
where $\xi \in \partial X$ is the endpoint of $\gamma$, hence
$$
d_{\mathcal{M}}( \rho_x, \rho_y ) = d(x, y)
$$
as required.

\medskip

Now fix an $x \in X$. Let $\xi \in \partial_{P} X$, and let $x_n \in X$ be a sequence converging to $\xi$.
Note that then $( x_n | \xi)_{x_0} \to +\infty$ as $n \to \infty$. Let $\gamma_n$ be a geodesic
segment joining $x$ to $x_n$, then by geodesic completeness we can extend $\gamma_n$ to a geodesic ray $\widetilde{\gamma_n} : [0, \infty) \to X$.
Let $\xi_n \in \partial_{P} X$ be the endpoint of $\widetilde{\gamma_n}$, then since $x_n$ lies on the ray $\widetilde{\gamma_n}$ which starts from $x$,
we have $\xi_n \in \hbox{argmax} \ \frac{d\rho_{x_n}}{d\rho_{x}}$.

\medskip

Let $\gamma : [0, \infty) \to X$ be a geodesic ray with $\gamma(0) = x, \gamma(\infty) = \xi$. Then for a fixed $n$ if we take $t > d(x, x_n)$
we have the elementary identity
$$
( \widetilde{\gamma_n}(t) | \gamma(t))_{x} = ( x_n | \gamma(t))_{x} + ( \widetilde{\gamma_n}(t) | \gamma(t))_{x_n}.
$$
Keeping $n$ fixed and letting $t$ tend to $\infty$ above gives
\begin{align*}
( \xi_n | \xi )_{x_0} & = ( x_n | \xi )_{x_0} + ( \xi_n | \xi )_{x_n} \\
                      & \geq ( x_n | \xi )_{x_0},
\end{align*}
from which it follows that $( \xi_n | \xi )_{x_0} \to +\infty$ as $n \to \infty$, and hence $\xi_n \to \xi$ in $\partial_{P} X$
as $n \to \infty$. Since $\xi_n \in \hbox{argmax} \ \frac{d\rho_{x_n}}{d\rho_{x}}$ for all $n$ and $d(x, x_n) \to \infty$,
it follows from Lemma \ref{maximaconv} that $i_X(x_n) = \rho_{x_n} \to \xi$ in $\widehat{\mathcal{M}(\partial X)}$
as $n \to \infty$. Therefore $i_X : X \to \mathcal{M}(\partial X)$ extends to a continuous map
$\widehat{i_X} : \widehat{X} \to \widehat{\mathcal{M}(\partial_{P} X)}$ such that $\widehat{i_X}_{|\partial_{P} X} = id : \partial_{P} X \to \partial_{P} X$,
as required.
$\diamond$

\medskip

We will call the isometric embedding $i_X$ the {\it visual embedding} of the Gromov product space $X$. We will need the following Lemma:

\medskip

\begin{lemma} \label{boundid} Let $X$ be a Gromov product space.
Let $F : X \to X$ be an isometric embedding, and let $f : \partial_{P} X \to \partial_{P} X$ be the boundary map of $F$.
If $f = {id}_{\partial_{P} X}$, then $F = {id}_X$.
\end{lemma}

\medskip

\noindent{\bf Proof:} Let $x \in X$. Given $\xi, \eta \in \partial_{P} X$, let $y_n, z_n \in X$ be sequences
converging to $\xi, \eta$ respectively. Then $F(x_n) \to f(\xi), F(z_n) \to f(\eta)$
as $n \to \infty$. Since $F$ is an isometric embedding and $f = {id}_{\partial X}$
we get
\begin{align*}
\rho_{F(x)}(\xi, \eta) & = \rho_{F(x)}(f(\xi), f(\eta)) \\
                       & = \lim_{n \to \infty} \exp\left( - ( F(y_n) | F(z_n) )_{F(x)} \right) \\
                       & = \lim_{n \to \infty} \exp\left( - ( y_n | z_n )_{x} \right) \\
                       & = \rho_x( \xi, \eta),
\end{align*}
thus $\rho_{F(x)} = \rho_x$. From Proposition \ref{isomembed} we then have $d(x, F(x)) = d_{\mathcal{M}}( \rho_x, \rho_{F(x)} ) = 0$,
thus $F(x) = x$ and $F = {id}_X$ as required. $\diamond$

\medskip

The above Lemma will allow us to define a partial order on the set of fillings of an antipodal space as follows:

\medskip

\begin{definition} {\bf (Partial order on fillings of antipodal spaces)} Let $Z$ be an antipodal space. A {\it filling} of $Z$ is a pair
$(X, f)$, where $X$ is a Gromov product space, and $f : \partial_{P} X \to Z$ is a Moebius homeomorphism.

\medskip

Given two fillings $(X, f)$ and $(Y, g)$ of $Z$, we say that the fillings are equivalent if
the Moebius homeomorphism $g^{-1} \circ f : \partial_{P} X \to \partial_{P} Y$ extends to an isometry $H : X \to Y$.
In this case we write $(X, f) \simeq (Y, g)$.

\medskip

We say that the filling $(Y, g)$ {\it dominates} the filling $(X, f)$, if
the Moebius homeomorphism $g^{-1} \circ f : \partial_{P} X \to \partial_{P} Y$ extends to an isometric embedding $H : X \to Y$.
In this case we write $(X, f) \leq (Y, g)$.
\end{definition}

\medskip

We remark that it is not at all clear whether a given antipodal space $Z$ has any fillings or not. The key point is that
the boundary map $f : \partial_{P} X \to Z$ is required to be Moebius, and this is nontrivial to achieve. We will show in the
next section that for any antipodal space $Z$, the Moebius space $\mathcal{M}(Z)$ is a Gromov product space which is a filling of $Z$.

\medskip

The following Proposition shows that we obtain a partial order on fillings:

\medskip

\begin{prop} Let $Z$ be an antipodal space.

\medskip

\noindent (1) If $(X_1, f_1) \leq (X_2, f_2)$ and $(X_2, f_2) \leq (X_3, f_3)$, then $(X_1, f_1) \leq (X_3, f_3)$ (here $(X_i, f_i)$ are fillings of $Z$ for $i = 1,2,3$).

\medskip

\noindent (2) If $(X, f) \leq (Y, g)$ and $(Y, g) \leq (X, f)$, then $(X, f) \simeq (Y, g)$ (here $(X, f)$ and $(Y, g)$ are fillings of $Z$).
\end{prop}

\medskip

\noindent{\bf Proof:} (1) Let $H_1 : X_1 \to X_2$ be an isometric embedding extending the map $f^{-1}_2 \circ f_1$, and let $H_2 : X_2 \to X_3$ be an
isometric embedding extending the map $f^{-1}_3 \circ f_2$, then $H_3 := H_2 \circ H_1 : X_1 \to X_3$ is an isometric embedding extending the map
$f^{-1}_3 \circ f_1$, hence $(X_1, f_1) \leq (X_3, f_3)$.

\medskip

(2) Let $H_1 : X \to Y$ be an isometric embedding extending the map $g^{-1} \circ f$, and let $H_2 : Y \to X$ be an isometric embedding extending the map $f^{-1} \circ g$,
then the map $H := H_2 \circ H_1 : X \to X$ is an isometric embedding whose boundary map satisfies $H_{|\partial_{P} X} = id_{\partial_{P} X}$, hence by
Lemma \ref{boundid} we have $H = id_X$, thus $H_1 : X \to Y$ is surjective, hence is an isometry, and so $(X, f) \simeq (Y, g)$ as required. $\diamond$

\medskip

The previous Proposition leads naturally to the notion of a {\it maximal Gromov product space}:

\medskip

\begin{definition} {\bf (Maximal Gromov product space)} Let $X$ be a Gromov product space and let $Z = \partial_{P} X$ be its
Gromov product boundary. We say that $X$ is {\it maximal} if the filling $(X, id_Z : Z \to Z)$ dominates every other filling of $Z$,
$(Y, f) \leq (X, id_Z)$ for any filling $(Y, f)$ of $Z$.

\medskip

In other words, whenever $Y$ is a Gromov product space equipped with a Moebius homeomorphism $f : \partial_{P} Y \to \partial_{P} X$,
then the Moebius map $f$ extends to an isometric embedding $F : Y \to X$.
\end{definition}

\medskip

The existence of maximal Gromov product spaces is not clear. We will show however in the next section that in fact for any antipodal space
$Z$, the Moebius space $\mathcal{M}(Z)$
is a maximal Gromov product space. Since a Gromov product space comes equipped with its visual embedding $i_X : X \to \mathcal{M}(\partial_{P} X)$,
it will follow that any Gromov product space has a canonical isometric embedding into a maximal Gromov product space with the same boundary.

\medskip

We now turn to CAT(0) spaces. We recall that a proper CAT(0) space $X$ admits a compactification $\overline{X} = X \cup \partial X$ called
the {\it visual compactification}, where $\partial X$ consists of equivalence classes of geodesic rays in $X$,
and the topology on $\overline{X}$ is called the {\it cone topology} (see Chapter II.8 of
\cite{bridsonhaefliger} for details). When the CAT(0) space $X$ is also Gromov hyperbolic, then it is not hard
to show that the visual compactification and the Gromov compactification coincide (\cite{bridsonhaefliger},
Chapter III.H), and so we can and will unambiguously use the notation $\overline{X} = X \cup \partial X$
for both.

\medskip

We also recall briefly the {\it horofunction compactification} of a proper metric space $X$. This can be
described as follows. Let $C(X)$ denote the space of continuous functions on $X$ endowed
with the topology of uniform convergence on compacts. Fix a basepoint $x_0 \in X$, and define for each
$x \in X$ a function $d_x : X \to \mathbb{R}$ by $d_x(y) = d(y, x) - d(x_0, x), y \in X$. Then the map
$\phi: x \in X \mapsto d_x \in C(X)$ is a homeomorphism onto its image, which is contained in the subspace of
$C(X)$ consisting of functions vanishing at $x_0$. The horofunction compactification
is defined to be the closure $\hat{X} \subset C(X)$ in $C(X)$ of the image of $X$ under this map. It is
easy to see that this compactification is independent of the choice of basepoint $x_0$.
The {\it horofunction boundary} of $X$ is defined to be $\partial_{H} X := \hat{X} - X$, and a function
$h \in C(X)$ is said to be a {\it horofunction} if $h+c \in \partial_{H} X$ for some constant $c \in \mathbb{R}$.
It can be shown that for a proper CAT(0) space, the horofunction and visual compactifications are canonically
homeomorphic (\cite{bridsonhaefliger}, Chapter II.8), namely the map $\phi : X \to C(X)$ extends
to a homeomorphism $\overline{\phi} : \overline{X} \to \hat{X}$.

\medskip

The following Proposition gives a large supply of examples of boundary continuous Gromov hyperbolic spaces.

\medskip

\begin{prop} \label{catohyp} Let $X$ be a proper, CAT(0) and Gromov hyperbolic space. Then $X$ is boundary continuous.
\end{prop}

\medskip

\noindent{\bf Proof:} Fix a basepoint $x_0 \in X$. We have to show that given $a, b \in \overline{X}$,
the limit of $(y | z)_{x_0}$ in $[0, +\infty]$ exists as $y, z \in X$ converge to $a,b$ in $\overline{X}$.
Given that $X$ is CAT(0), we let $\overline{\phi} : \overline{X} \to \hat{X}$ be the homeomorphism
between the visual and horofunction compactifications described above.

\medskip

We first consider the case when $a \in X$ and $b = \xi \in \partial X$. Let $h = \overline{\phi}(\xi) \in \partial_{H} X$.
Let $y_n, z_n \in X$ be sequences such that $y_n \to a, z_n \to \xi$ in $\overline{X}$ as $n \to \infty$. Then $\phi(z_n) = d_{z_n} \to \overline{\phi}(\xi) = h$
in $C(X)$ as $n \to \infty$. Thus the functions $d_{z_n}$ converge to $h$ uniformly on compacts in $X$, so since $y_n \to a \in X$ we get
\begin{align*}
2( y_n | z_n)_{x_0} & = d(x_0, y_n) + ( d(x_0, z_n) - d(y_n, z_n)) \\
                    & = d(x_0, y_n) - d_{z_n}(y_n) \\
                    & \to d(x_0, a) - h(a) \ \hbox{as} \ n \to \infty,
\end{align*}
and so it follows that the limit
$$
\lim_{y \to a, z \to \xi} (y | z)_{x_0} = \frac{1}{2} ( d(x_0, a) - h(a))
$$
exists as required.

\medskip

We now consider the case when $a = \xi \in \partial X, b = \eta \in \partial X$. If $\xi = \eta$, then by definition of the
Gromov compactification $(y | z)_{x_0} \to +\infty$ as $y,z$ tend to $\xi = \eta$ in $\overline{X}$, and we are done. So we may assume
$\xi \neq \eta$. Let $h_1 = \overline{\phi}(\xi) \in \partial_{H} X, h_2 = \overline{\phi}(\eta) \in \partial_{H} X$.
Suppose the limit of $(y | z)_{x_0}$ as $y \to \xi, z \to \eta$ does not exist. Then there exist two pairs of sequences
$y_n, z_n \in X$ and $y'_n, z'_n \in X$ such $y_n, y'_n \to \xi, z_n, z'_n \to \eta$ in $\overline{X}$ as $n \to \infty$,
and $( y_n | z_n)_{x_0} \to \lambda, (y'_n | z'_n)_{x_0} \to \lambda'$ as $n \to \infty$, and $\lambda \neq \lambda'$.

\medskip

Since $\xi \neq \eta$, there exists a constant $C > 0$ such that $(y_n | z_n)_{x_0}, (y'_n | z'_n)_{x_0} \leq C$ for all $n$.
Recall that CAT(0) spaces are uniquely geodesic, for any $y,z \in X$ let $[y, z]$ denote the unique geodesic segment joining $y$ to $z$.
Let $\delta \geq 0$ be such that $X$ is $\delta$-hyperbolic, then we have the well-known inequality for $\delta$-hyperbolic spaces,
$d(x_0, [y,z]) \leq (y | z)_{x_0} + 2\delta$ for all $y, z \in X$. It follows that for all $n$ there exist points $p_n \in [y_n, z_n]$
and $p'_n \in [y'_n, z'_n]$ such that $d(x_0, p_n), d(x_0, p'_n) \leq C + 2\delta$. Since $X$ is proper, after passing to a subsequence,
we may assume that $p_n \to p, p'_n \to p'$ as $n \to \infty$ for some $p,p' \in X$. It is then well-known and not hard to show
that after parametrizing the geodesic segments $[y_n, z_n], [y'_n, z'_n]$ such that
they pass through $p_n, p'_n$ at time $t = 0$, these geodesic segments converge uniformly on compacts to bi-infinite geodesics
$\gamma : \mathbb{R} \to X, \gamma' : \mathbb{R} \to X$ respectively, such that $\gamma(0) = p, \gamma'(0) = p'$,
$\gamma(-\infty) = \gamma'(-\infty) = \xi$, and $\gamma(\infty) = \gamma'(\infty) = \eta$.

\medskip

Now noting that $\phi(y_n) = d_{y_n} \to h_1, \phi(z_n) = d_{z_n} \to h_2$ as $n \to \infty$, we have
\begin{align*}
2( y_n | z_n)_{x_0} & = (d(x_0, y_n) - d(p_n, y_n)) + (d(x_0, z_n) - d(p_n, z_n)) \\
                    & = -d_{y_n}(p_n) - d_{z_n}(p_n) \\
                    & \to -h_1(p) - h_2(p) \ \hbox{as} \ n \to \infty,
\end{align*}

similarly we have
$$
2( y'_n | z'_n)_{x_0} \to -h_1(p') - h_2(p') \ \hbox{as} \ n \to \infty.
$$
Thus $2 \lambda = -h_1(p) - h_2(p)$ and $2\lambda' = - h_1(p') - h_2(p')$.

\medskip

Now the bi-infinite geodesics $\gamma, \gamma'$ have the same endpoints $\xi, \eta$ in $\partial X$,
hence are asymptotic, i.e. $d(\gamma(t), \gamma'(t)) \leq K$ for all $t \in \mathbb{R}$ for some $K \geq 0$. Since
$X$ is CAT(0), this implies by the Flat Strip Theorem (\cite{bridsonhaefliger}, Chapter II.2) that the convex hull
$S \subset X$ of $\gamma(\mathbb{R}) \cup \gamma'(\mathbb{R})$ in $X$ is isometric to a flat Euclidean strip
$\mathbb{R} \times [0, d] \subset \mathbb{R}^2$, for some $d \geq 0$. We identify the convex hull $S$ with the
flat strip $\mathbb{R} \times [0, d]$, such that $\gamma(\mathbb{R})$ is identified with $\mathbb{R} \times \{0\}$, and
$\gamma'(\mathbb{R})$ is identified with $\mathbb{R} \times \{d\}$. Then for any sequence $(t_n, r_n) \in \mathbb{R} \times [0, d] \simeq S$,
we have $(t_n, r_n) \to \xi$ in $\overline{X}$ if and only if $t_n \to -\infty$, and $(t_n, r_n) \to \eta$ in $\overline{X}$ if and only if
$t_n \to +\infty$.

\medskip

Suppose that the points $p, p' \in S$ correspond to $(t_0, 0)$ and $(t_1, d)$ respectively.
For $t \in \mathbb{R}$, let $q_t \in S$ correspond to $(t, d/2)$. Then, since $S$ is a Euclidean strip,
\begin{align*}
d(p, q_t) & = \sqrt{(t - t_0)^2 + (d/2)^2} \\
          & = | t - t_0| + O(1/|t|) \ \hbox{as} \ |t| \to \infty,
\end{align*}
similarly
$$
d(p', q_t) = | t - t_1| + O(1/|t|) \ \hbox{as} \ |t| \to \infty.
$$
Since $\phi(q_t) = d_{q_t} \to h_1$ as $t \to -\infty$, and $\phi(q_t) = d_{q_t} \to h_2$ as $t \to +\infty$, we have
\begin{align*}
h_1(p) - h_1(p') & = \lim_{t \to -\infty} \left(d_{q_t}(p) - d_{q_t}(p')\right) \\
                 & = \lim_{t \to -\infty} (d(p, q_t) - d(p', q_t)) \\
                 & = \lim_{t \to -\infty} ( -(t - t_0) + (t - t_1) + O(1/|t|)) \\
                 & = t_0 - t_1,
\end{align*}
and similarly
\begin{align*}
h_2(p) - h_2(p') & = \lim_{t \to +\infty} \left(d_{q_t}(p) - d_{q_t}(p')\right) \\
                 & = \lim_{t \to +\infty} (d(p, q_t) - d(p', q_t)) \\
                 & = \lim_{t \to +\infty} ( (t - t_0) - (t - t_1) + O(1/|t|)) \\
                 & = t_1 - t_0,
\end{align*}
and so
\begin{align*}
2\lambda - 2\lambda' & = -(h_1(p) - h_1(p')) - (h_2(p) - h_2(p')) \\
                     & = -(t_0 - t_1) - (t_1 - t_0) \\
                     & = 0,
\end{align*}
a contradiction to our hypothesis that $\lambda \neq \lambda'$. Thus the limit of $(y | z)_{x_0}$ exists as $y \to \xi, z \to \eta$
in $\overline{X}$, and so $X$ is boundary continuous. $\diamond$

\medskip

\section{Antipodal spaces and maximal Gromov product spaces.}

\medskip

In this section we define a class of Gromov product spaces called {\it maximal Gromov product spaces}, and then
prove an equivalence of categories between antipodal spaces and maximal Gromov product spaces (where the morphisms in the first category are
Moebius homeomorphisms, and the morphisms in the second category are surjective isometries).

\medskip

Let $Z$ be an antipodal space. The next Proposition proves an analogue of Lemma \ref{visualderiv} from the previous section, for the Moebius space $\mathcal{M}(Z)$
associated to $Z$, showing that in this setting too we can define a Busemann cocycle, and it is again
given by the log of the derivative. In this case, the Busemann cocycle is a function $B : \mathcal{M}(Z) \times \mathcal{M}(Z) \times Z \to \mathbb{R}$.

\medskip

\begin{prop} \label{mzbuse} Let $Z$ be an antipodal space. For any $\rho_1, \rho_2 \in \mathcal{M}(Z)$ and $\xi \in Z$, we have
$$
\lim_{\alpha \to \xi} (d_{\mathcal{M}}(\rho_1, \alpha) - d_{\mathcal{M}}(\rho_2, \alpha)) = \log \frac{d\rho_2}{d\rho_1}(\xi)
$$
(where the above limit is taken as $\alpha \in \mathcal{M}(Z)$ converges in $\widehat{\mathcal{M}(Z)}$ to $\xi \in Z$).
\end{prop}

\medskip

\noindent{\bf Proof:} Let $\rho_1, \rho_2 \in \mathcal{M}(Z)$ and $\xi \in Z$. Given $\epsilon > 0$, let $U$ be an
open neighbourhood of $\xi$ such that
$$
\left| \log \frac{d\rho_2}{d\rho_1}(\eta) - \log \frac{d\rho_2}{d\rho_1}(\xi) \right| < \epsilon
$$
for all $\eta \in U$.

\medskip

Let $r > 0$, then $V := N(U, \rho_1, r) \cap N(U, \rho_2, r)$ is an open neighbourhood of $\xi$ in
$\widehat{\mathcal{M}(Z)}$. Given $\alpha \in V \cap \mathcal{M}(Z)$, let $\eta_i \in \hbox{argmax} \ \frac{d\alpha}{d\rho_i}$
for $i = 1,2$, then $\eta_1, \eta_2 \in U$. From the inequalities
$$
\frac{d\alpha}{d\rho_1}(\eta_1) \geq \frac{d\alpha}{d\rho_1}(\eta_2) \ , \ \frac{d\alpha}{d\rho_2}(\eta_2) \geq \frac{d\alpha}{d\rho_2}(\eta_1),
$$
we get
\begin{align*}
d_{\mathcal{M}}(\rho_1, \alpha) - d_{\mathcal{M}}(\rho_2, \alpha) & = \log \frac{d\alpha}{d\rho_1}(\eta_1) - \log \frac{d\alpha}{d\rho_2}(\eta_2) \\
                                                                  & \geq \log \frac{d\alpha}{d\rho_1}(\eta_2) - \log \frac{d\alpha}{d\rho_2}(\eta_2) \\
                                                                  & = \log \frac{d\rho_2}{d\rho_1}(\eta_2) \\
                                                                  & > \log \frac{d\rho_2}{d\rho_1}(\xi) - \epsilon,
\end{align*}

and similarly

\begin{align*}
d_{\mathcal{M}}(\rho_1, \alpha) - d_{\mathcal{M}}(\rho_2, \alpha) & = \log \frac{d\alpha}{d\rho_1}(\eta_1) - \log \frac{d\alpha}{d\rho_2}(\eta_2) \\
                                                                  & \leq \log \frac{d\alpha}{d\rho_1}(\eta_1) - \log \frac{d\alpha}{d\rho_2}(\eta_1) \\
                                                                  & = \log \frac{d\rho_2}{d\rho_1}(\eta_1) \\
                                                                  & < \log \frac{d\rho_2}{d\rho_1}(\xi) + \epsilon,
\end{align*}

thus
$$
\lim_{\alpha \to \xi} (d_{\mathcal{M}}(\rho_1, \alpha) - d_{\mathcal{M}}(\rho_2, \alpha)) = \log \frac{d\rho_2}{d\rho_1}(\xi)
$$
as required. $\diamond$

\medskip

The following Lemma will be useful:

\medskip

\begin{lemma} \label{giproupper} Let $\rho_1, \alpha, \beta \in \mathcal{M}(Z)$. Let
$\xi \in \hbox{argmax} \ \frac{d\alpha}{d\rho_1}, \eta \in \hbox{argmax} \ \frac{d\beta}{d\rho_1}$. Then
$$
( \alpha | \beta )_{\rho_1} \leq \log \frac{1}{\rho_1(\xi, \eta)}
$$
(with the convention that the right-hand side above is $+\infty$ if $\xi = \eta$).
\end{lemma}

\medskip

\noindent{\bf Proof:} We may as well assume $\xi \neq \eta$, otherwise there is nothing to prove.
Let $\rho \in \mathcal{M}(Z)$ be a midpoint of
$\alpha$ and $\beta$ (for example given by Lemma \ref{tcomb}). Then we have
\begin{align*}
2(\alpha|\beta)_{\rho_1} & = d_{\mathcal{M}}(\rho_1, \alpha) + d_{\mathcal{M}}(\rho_1, \beta) - d_{\mathcal{M}}(\alpha, \beta) \\
                         & = d_{\mathcal{M}}(\rho_1, \alpha) + d_{\mathcal{M}}(\rho_1, \beta) - d_{\mathcal{M}}(\alpha, \rho) - d_{\mathcal{M}}(\beta, \rho) \\
                         & = \log \frac{d\alpha}{d\rho_1}(\xi) + \log \frac{d\beta}{d\rho_1}(\eta) - d_{\mathcal{M}}(\alpha, \rho) - d_{\mathcal{M}}(\beta, \rho) \\
                         & \leq \log \frac{d\alpha}{d\rho_1}(\xi) + \log \frac{d\beta}{d\rho_1}(\eta) - \log \frac{d\alpha}{d\rho}(\xi) - \log \frac{d\beta}{d\rho}(\eta) \\
                         & = \log \frac{d\rho}{d\rho_1}(\xi) + \log \frac{d\rho}{d\rho_1}(\eta) \\
                         & = \log \frac{\rho(\xi,\eta)^2}{\rho_1(\xi, \eta)^2} \\
                         & \leq \log \frac{1}{\rho_1(\xi, \eta)^2}.
\end{align*}
$\diamond$

\medskip

The next theorem says that the Moebius space $\mathcal{M}(Z)$ is a Gromov product space, with Gromov product compactification
given by the maxima compactification, and Gromov product boundary given by the antipodal space $Z$. Moreover, it shows that any
Moebius antipodal function on $Z$ is a visual antipodal function for the space $\mathcal{M}(Z)$:

\medskip

\begin{theorem} \label{gromovip} Let $\rho_1 \in \mathcal{M}(Z)$. Then for any $\xi_0, \eta_0 \in Z$, we have
$$
\rho_1(\xi_0, \eta_0) = \lim_{\alpha \to \xi_0, \beta \to \eta_0} \exp\left(-(\alpha | \beta )_{\rho_1}\right)
$$
(where in the above limit the points $\alpha, \beta \in \mathcal{M}(Z)$ converge to $\xi_0, \eta_0 \in Z$ in the
maxima compactification $\widehat{\mathcal{M}(Z)}$).

\medskip

In particular, the space $\mathcal{M}(Z)$ is a Gromov product space, with Gromov product compactification given by the
maxima compactification, $\widehat{\mathcal{M}(Z)}^{P} = \widehat{\mathcal{M}(Z)}$, and Gromov product boundary
given by the antipodal space $Z$, $\partial_{P} \mathcal{M}(Z) = Z$.
\end{theorem}

\medskip

\noindent{\bf Proof:} Let $\rho_1 \in \mathcal{M}(Z)$ and let $\xi_0, \eta_0 \in Z$. We consider two cases:

\medskip

\noindent{\bf Case (i)} $\xi_0 \neq \eta_0$:

\medskip

Given $\epsilon > 0$, by continuity of $\rho_1$ we can choose open neighbourhoods $U_1, U_2$ of $\xi_0, \eta_0$ respectively in $Z$
with disjoint closures such that
\begin{equation} \label{u1u2}
\left|\log \frac{1}{\rho_1(\xi,\eta)} - \log \frac{1}{\rho_1(\xi_0, \eta_0)}\right| < \epsilon
\end{equation}
for all $\xi \in U_1, \eta \in U_2$. Let $r > 0$ and let $\alpha, \beta \in \mathcal{M}(Z)$ be
such that $\alpha \in N(U_1, \rho_1, r), \beta \in N(U_2, \rho_1, r)$. Let $\xi \in \hbox{argmax} \ \frac{d\alpha}{d\rho_1},
\eta \in \hbox{argmax} \ \frac{d\beta}{d\rho_1}$, then $\xi \in U_1, \eta \in U_2$. Then by Lemma \ref{giproupper}
\begin{align*}
(\alpha|\beta)_{\rho_1} & \leq \log \frac{1}{\rho_1(\xi, \eta)} \\
                        & < \log \frac{1}{\rho_1(\xi_0, \eta_0)} + \epsilon,
\end{align*}
and so it follows that
$$
\limsup_{\alpha \to \xi_0, \beta \to \eta_0} \ (\alpha | \beta)_{\rho_1} \ \leq \ \log \frac{1}{\rho_1(\xi_0, \eta_0)}.
$$

\medskip

Again, let $\epsilon > 0$ be fixed and let $U_1, U_2$ be neighbourhoods of $\xi_0, \eta_0$ with disjoint closures
as above satisfying (\ref{u1u2}). Let $V_1, V_2$ be neighbourhoods of $\xi_0, \eta_0$
respectively such that $\overline{V_1} \subset U_1, \overline{V_2} \subset U_2$. Then since $\overline{V_i}$ and $Z - U_i$ are
disjoint compact sets, it is easy to show by positivity and continuity of $\rho_1$ that there is a $\delta > 0$ such that, for $i = 1,2$,
if $\zeta \in \overline{V_i}$ and $\rho_1(\zeta, \zeta') < \delta$, then $\zeta' \in U_i$. Now let $r > 0$ be such that
$$
\frac{e^{-r}}{\rho_1(\xi, \eta)} < \delta
$$
for all $\xi \in V_1, \eta \in V_2$ (this is possible by positivity of $\rho_1$, since $V_1, V_2$ have disjoint closures).

\medskip

Let $\alpha, \beta \in \mathcal{M}(Z)$ be such that $\alpha \in N(V_1, \rho_1, r), \beta \in N(V_2, \rho_1, r)$.
Let $\xi \in \hbox{argmax} \ \frac{d\alpha}{d\rho_1},
\eta \in \hbox{argmax} \ \frac{d\beta}{d\rho_1}$, then $\xi \in V_1, \eta \in V_2$. Then by Lemma \ref{giproupper} we have
\begin{equation} \label{gipupper}
2(\alpha|\beta)_{\rho_1} \leq \log \frac{1}{\rho_1(\xi, \eta)^2}.
\end{equation}

Let $\xi_1 \in \hbox{argmax} \ \frac{d\alpha}{d\beta}$ and let $\eta_1 \in Z$ be such that $\beta(\xi_1, \eta_1) = 1$,
then by Lemma \ref{maxmin} we have $\eta_1 \in \hbox{argmin} \ \frac{d\alpha}{d\beta} = \hbox{argmax} \ \frac{d\beta}{d\alpha}$, and
$\alpha(\xi_1, \eta_1) = 1$.

\medskip

Let $\rho \in \mathcal{M}(Z)$ be a midpoint of $\alpha, \beta$, and let $t_1 = d_{\mathcal{M}}( \rho_1, \alpha), t_2 = d_{\mathcal{M}}( \rho_1, \beta)$
and $s = d_{\mathcal{M}}( \alpha, \beta)$. Then we have
\begin{align*}
1 & \geq \alpha(\xi, \xi_1)^2 \\
  & = \frac{d\alpha}{d\rho_1}(\xi) \cdot \frac{d\alpha}{d\rho_1}(\xi_1) \cdot \rho_1(\xi, \xi_1)^2 \\
  & = e^{t_1} \cdot \frac{d\alpha}{d\rho_1}(\xi_1) \cdot \rho_1(\xi, \xi_1)^2 \\
  & = e^{t_1} \cdot \frac{d\alpha}{d\beta}(\xi_1) \cdot \frac{d\beta}{d\rho_1}(\xi_1) \cdot \rho_1(\xi, \xi_1)^2 \\
  & \geq e^{t_1} \cdot e^s \cdot e^{-t_2} \cdot \rho_1(\xi, \xi_1)^2.
\end{align*}
Using the inequality (\ref{gipupper}), and noting that $t_1 > r$, the above inequality gives us
\begin{align*}
\rho_1(\xi, \xi_1)^2 & \leq e^{t_2 - t_1 - s} \\
                     & = e^{(t_2 + t_1 - s) - 2t_1} \\
                     & = e^{2(\alpha|\beta)_{\rho_1} - 2t_1} \\
                     & \leq \frac{e^{-2t_1}}{\rho_1(\xi, \eta)^2} \\
                     & < \frac{e^{-2r}}{\rho_1(\xi, \eta)^2} \\
                     & < \delta^2,
\end{align*}
thus, since $\xi \in V_1$, by our choice of $\delta$ it follows that $\xi_1 \in U_1$. Similarly we can show that $\eta_1 \in U_2$.

\medskip

Since $d_{\mathcal{M}}(\alpha, \beta) = d_{\mathcal{M}}(\alpha, \rho) + d_{\mathcal{M}}(\rho, \beta)$, it follows from
Lemma \ref{trieq} that $\xi_1 \in \hbox{argmax} \ \frac{d\alpha}{d\rho}$ and $\eta_1 \in \hbox{argmax} \ \frac{d\beta}{d\rho}$.
From Lemma \ref{trieq} we also have $\xi_1 \in \hbox{argmax} \ \frac{d\rho}{d\beta}$. Since $\beta(\xi_1, \eta_1) = 1$, by
Lemma \ref{maxmin} this implies that $\rho(\xi_1, \eta_1) = \beta(\xi_1, \eta_1) = 1$. Using these facts, we obtain
\begin{align*}
2(\alpha|\beta)_{\rho_1} & = d_{\mathcal{M}}(\rho_1, \alpha) + d_{\mathcal{M}}(\rho_1, \beta) - d_{\mathcal{M}}(\alpha, \rho) - d_{\mathcal{M}}(\beta, \rho) \\
                         & \geq \log \frac{d\alpha}{d\rho_1}(\xi_1) + \log \frac{d\beta}{d\rho_1}(\eta_1) - d_{\mathcal{M}}(\alpha, \rho) - d_{\mathcal{M}}(\beta, \rho) \\
                         & = \log \frac{d\alpha}{d\rho_1}(\xi_1) + \log \frac{d\beta}{d\rho_1}(\eta_1) - \log \frac{d\alpha}{d\rho}(\xi_1) - \log \frac{d\beta}{d\rho}(\eta_1) \\
                         & = \log \frac{d\rho}{d\rho_1}(\xi_1) + \log \frac{d\rho}{d\rho_1}(\eta_1) \\
                         & = \log \frac{\rho(\xi_1,\eta_1)^2}{\rho_1(\xi_1, \eta_1)^2} \\
                         & = \log \frac{1}{\rho_1(\xi_1, \eta_1)^2} \\
                         & \geq 2 \log \frac{1}{\rho_1(\xi_0, \eta_0)} - 2\epsilon.
\end{align*}

It follows that
$$
\liminf_{\alpha \to \xi_0, \beta \to \eta_0} \ (\alpha | \beta)_{\rho_1} \ \geq \ \log \frac{1}{\rho_1(\xi_0, \eta_0)},
$$
and so
$$
\lim_{\alpha \to \xi_0, \beta \to \eta_0} (\alpha | \beta )_{\rho_1} = \log \frac{1}{\rho_1(\xi_0, \eta_0)}
$$
as required.

\medskip

\noindent{\bf Case (ii)} $\xi_0 = \eta_0$: \ In this case we need to show that $(\alpha | \beta)_{\rho_1} \to +\infty$ as
$\alpha, \beta \to \xi_0$ in $\widehat{\mathcal{M(Z)}}$.

\medskip

Let $M > 0$. We can choose an open neighbourhood $U_0$ of $\xi_0$ in $Z$ such that $\log \frac{1}{\rho_1(\xi, \eta)} \geq M$
for all $\xi, \eta \in U_0$. Let $U \subset U_0$ be a neighbourhood of $\xi_0$ in $Z$ such that $\overline{U} \subset U_0$, then
there is an $\epsilon > 0$ such that $\rho_1(\xi, \eta) \geq \epsilon$ for all $\xi \in U, \eta \in Z \setminus U_0$. Let $V$ be a neighbourhood of $\xi_0$ in $\widehat{\mathcal{M(Z)}}$ such that
$\hbox{argmax} \frac{d\alpha}{d\rho_1} \subset U$ and $d_{\mathcal{M}}(\alpha, \rho_1) \geq M + \log \frac{1}{\epsilon^2}$
for all $\alpha \in V \cap \mathcal{M}(Z)$.

\medskip

Let $\alpha, \beta \in V \cap \mathcal{M}(Z)$. Let $\xi_1 \in \hbox{argmax} \ \frac{d\alpha}{d\beta}$ and let $\eta_1 \in Z$ be such that $\beta(\xi_1, \eta_1) = 1$,
then by Lemma \ref{maxmin} we have $\eta_1 \in \hbox{argmin} \ \frac{d\alpha}{d\beta} = \hbox{argmax} \ \frac{d\beta}{d\alpha}$, and
$\alpha(\xi_1, \eta_1) = 1$.

\medskip

By the same argument as before, considering a midpoint $\rho$ of $\alpha, \beta$, we have
$$
(\alpha | \beta)_{\rho_1} \geq \log \frac{1}{\rho_1(\xi_1, \eta_1)},
$$
thus if both $\xi_1, \eta_1$ lie in $U_0$ then $(\alpha | \beta)_{\rho_1} \geq M$.

\medskip

Otherwise, we may assume that one of the two points, say $\eta_1$, lies outside $U_0$. Let $\eta \in \hbox{argmax} \frac{d\beta}{d\rho_1}$.
Then $\eta \in U$, so $\rho_1(\eta, \eta_1) \geq \epsilon$, thus
\begin{align*}
0 & = \log \frac{d\beta}{d\rho_1}(\eta) + \log \frac{d\beta}{d\rho_1}(\eta_1) + \log \rho_1(\eta, \eta_1)^2 \\
  & \geq d_{\mathcal{M}}(\beta, \rho_1) + \log \frac{d\beta}{d\rho_1}(\eta_1) + \log \epsilon^2 \\
  & \geq \log \frac{d\beta}{d\rho_1}(\eta_1) + M
\end{align*}
(using $\beta \in V$), thus
$$
\log \frac{d\beta}{d\rho_1}(\eta_1) \leq -M.
$$
We then have
\begin{align*}
2(\alpha | \beta)_{\rho_1} & = d_{\mathcal{M}}(\alpha, \rho_1) + d_{\mathcal{M}}(\beta, \rho_1) - d_{\mathcal{M}}(\alpha, \beta) \\
                           & \geq -\log \frac{d\alpha}{d\rho_1}(\eta_1) + d_{\mathcal{M}}(\beta, \rho_1) - \log \frac{d\beta}{d\alpha}(\eta_1) \\
                           & = d_{\mathcal{M}}(\beta, \rho_1) - \log \frac{d\beta}{d\rho_1}(\eta_1) \\
                           & \geq 2M.
\end{align*}

Thus in either case we have $(\alpha | \beta)_{\rho_1} \geq M$ for all $\alpha, \beta \in V \cap \mathcal{M}(Z)$, and so
$$
\lim_{\alpha, \beta \to \xi_0} (\alpha | \beta)_{\rho_1} = +\infty
$$
as required. $\diamond$

\medskip

Thus for any antipodal space $Z$, the Moebius space $\mathcal{M}(Z)$ is a Gromov product space which is a filling of $Z$, which we proceed to show is
in fact maximal.

\medskip

We first observe that if $Z_1, Z_2$ are antipodal spaces and $f : Z_1 \to Z_2$ is a Moebius homeomorphism, then $f$ induces
a map $f_* : \mathcal{M}(Z_1) \to \mathcal{M}(Z_2)$, given by push-forward of antipodal Moebius functions. Explicitly,
for $\rho \in \mathcal{M}(Z_1)$, the point $f_* \rho \in \mathcal{M}(Z_2)$ is defined by
$$
(f_*\rho)(\xi, \eta) := \rho(f^{-1}(\xi), f^{-1}(\eta)) \ , \xi, \eta \in Z_2.
$$

We then have:

\medskip

\begin{prop} \label{mbisom1} Let $Z_1, Z_2$ be antipodal spaces, and let $f : Z_1 \to Z_2$ be a Moebius homeomorphism.
Then the associated map $F = f_* : \mathcal{M}(Z_1) \to \mathcal{M}(Z_2)$ is an isometry. Moreover, $f$ is the
boundary map of $F$, in the sense that $F$ extends to a homeomorphism
$\widetilde{F} : \widehat{\mathcal{M}(Z_1)} \to \widehat{\mathcal{M}(Z_2)}$ such that $\widetilde{F}_{|Z_1} = f$.
\end{prop}

\noindent{\bf Proof:} It follows from the definition (\ref{derivdefn}) of the derivative that
\begin{equation} \label{derivcomp}
\frac{d(f_* \rho_2)}{d(f_* \rho_1)}(\xi) = \frac{d\rho_2}{d\rho_1}(f^{-1}(\xi))
\end{equation}
for any $\rho_1, \rho_2 \in \mathcal{M}(Z_1)$, $\xi \in Z_2$. Since $f$ is a bijection, we get
\begin{align*}
d_{\mathcal{M}}( f_* \rho_1, f_* \rho_2) & = \max_{\xi \in Z_2} \log \frac{d(f_* \rho_2)}{d(f_* \rho_1)}(\xi) \\
                                         & = \max_{\xi \in Z_2} \log \frac{d\rho_2}{d\rho_1}(f^{-1}(\xi)) \\
                                         & = \max_{\eta \in Z_1} \log \frac{d\rho_2}{d\rho_1}(\eta) \\
                                         & = d_{\mathcal{M}}( \rho_1, \rho_2),
\end{align*}
so $F = f_*$ is an isometric embedding, which is surjective since it has an inverse $G = g_*$ where $g = f^{-1}$.

\medskip

Given $\xi \in Z_1$, let $\rho_1 \in \mathcal{M}(Z_1)$ and let $N(V, f_* \rho_1, r) \subset \widehat{\mathcal{M}(Z_2)}$ be a
basic open neighbourhood of $f(\xi) \in Z_2$ in $\widehat{\mathcal{M}(Z_2)}$, where $V \subset Z_2$ is an
open neighbourhood of $f(\xi)$ in $Z_2$ and $r > 0$. Let $U \subset Z$ be an open neighbourhood of $\xi$ in $Z_1$ such that $f(U) \subset V$.
Consider the basic open neighbourhood $N(U, \rho_1, r) \subset \widehat{\mathcal{M}(Z_1)}$ of
$\xi$ in $\widehat{\mathcal{M}(Z_1)}$.
Then for any $\rho \in N(U, \rho_1, r) \cap \mathcal{M}(Z)$, it follows from (\ref{derivcomp}) above that
$\hbox{argmax} \ \frac{d(f_* \rho)}{d(f_* \rho_1)} = f \left(\hbox{argmax} \ \frac{d\rho}{d\rho_1}\right)$, hence
$f_* \rho \in N(f(U), f_* \rho_1, r) \subset N(V, f_* \rho_1, r) \subset \widehat{\mathcal{M}(Z_2)}$. It follows that
$$
\lim_{\rho \to \xi} F(\rho) = f(\xi) \in Z_2,
$$
and $F$ extends to a continuous map $\widetilde{F} : \widehat{\mathcal{M}(Z_1)} \to \widehat{\mathcal{M}(Z_2)}$
such that $\widetilde{F}_{|Z_1} = f$. Moreover $\widetilde{F}$ is a homeomorphism since it has a continuous inverse given by
$\widetilde{G}$, which is the similar extension of $G$. $\diamond$

\medskip

We can now prove Theorem \ref{mainthm6} from the Introduction on the maximality of $\mathcal{M}(Z)$:

\medskip

\noindent{\bf Proof of Theorem \ref{mainthm6}:} Let $(X, f)$ be a filling of $Z$. By Proposition \ref{mbisom1}, the Moebius homeomorphism $f : \partial_{P} X \to Z$
extends to an isometry $f_* : \mathcal{M}(\partial_{P} X) \to \mathcal{M}(Z)$.

\medskip

The composition $H := f_* \circ i_X : X \to \mathcal{M}(Z)$ is then an isometric embedding of $X$ into $\mathcal{M}(Z)$, where the map $i_X$
is the visual embedding of $X$. By Propositions \ref{isomembed} and \ref{mbisom1}, the boundary map of the embedding $H$ is given
by $H_{\partial_{P} X} = f \circ id_{\partial_{P} X} = f$, and so $H$ is an extension of $f$ and $(X, f) \leq (\mathcal{M}(Z), id_Z)$. 

\medskip

Thus $\mathcal{M}(Z)$ is a maximal Gromov product space.
$\diamond$

\medskip

We now prove Theorem \ref{mainthm7} from the Introduction, which will show that all maximal Gromov product spaces are in fact Moebius spaces:

%

\medskip

\noindent{\bf Proof of Theorem \ref{mainthm7}:} Suppose the visual embedding $i_X$ is an isometry. Then the maximality of $\mathcal{M}(\partial_{P} X)$
easily implies that of $X$. Indeed, for any filling $(Y, g)$ of $\partial_{P} X$, since $\partial_{P} X$ is also the Gromov
product boundary of $\mathcal{M}(\partial_{P} X)$, by maximality of $\mathcal{M}(\partial_{P} X)$ the Moebius map $g : \partial_{P} Y \to \partial_{P} X$
extends to an isometric embedding $H : Y \to \mathcal{M}(\partial_{P} X)$. Then the map $G := \left( i_X \right)^{-1} \circ H : Y \to X$ is an
isometric embedding of $Y$ into $X$, and by Proposition \ref{isomembed} the boundary map of $G$ is given by
$G_{|\partial_{P} Y} = id_{\partial_{P} X} \circ g = g$, thus $G$ extends $g$, and so $X$ is maximal.

\medskip

Conversely, suppose $X$ is maximal. If we let $f = id_{\partial_{P} X}$, then the pair $(\mathcal{M}(\partial_{P} X), f)$ is a
filling of $\partial_{P} X$ (here we are identifying the Gromov product boundary of $\mathcal{M}(\partial_{P} X)$ with
$\partial_{P} X$). Thus by maximality of $X$, the map $f$ extends to an isometric embedding $F : \mathcal{M}(\partial_{P} X) \to X$.
The composition $G := F \circ i_X : X \to X$ is then an isometric embedding of $X$ into itself, and by Proposition \ref{isomembed}
the boundary map of $G$ is given by $G_{|\partial_{P} X} = f \circ id_{\partial_{P} X} = f = id_{\partial_{P} X}$, and hence
from Lemma \ref{boundid} we have $G = id_X$. It follows that $i_X$ is surjective, as required. $\diamond$

\medskip

To complete the equivalence of categories between antipodal spaces $Z$ and maximal Gromov product spaces,
we give below the proof of Theorem \ref{mainthm8}:

\medskip

\noindent{\bf Proof of Theorem \ref{mainthm8}:} Let $X, Y$ be maximal Gromov product spaces and let $f : \partial_{P} X \to \partial_{P} Y$
be a Moebius homeomorphism. Then the induced map $f_* : \mathcal{M}(\partial_{P} X) \to \mathcal{M}(\partial_{P} Y)$ is an isometry
with boundary value equal to $f$. By Theorem \ref{mainthm7}, since $X, Y$ are maximal the visual embeddings $i_X : X \to \mathcal{M}(\partial_{P} X),
i_Y : Y \to \mathcal{M}(\partial_{P} Y)$ are isometries, and by Proposition 
\ref{isomembed} their boundary values are equal to the identity maps of $\partial_{P} X$ and $\partial_{P} Y$ respectively.
It follows that the map $F := \left( i_Y \right)^{-1} \circ f_* \circ i_X : X \to Y$ is an isometry with boundary value equal to $f$ as required. $\diamond$

\medskip

\section{Gromov product spaces and injective metric spaces.}

\medskip

In this section we describe the close connection between Gromov product spaces and injective metric spaces.

\medskip

We recall that any metric space admits an isometric embedding into a "smallest"  injective metric space, namely the {\it injective hull}
of the metric space. There is a well-known construction due to Isbell \cite{isbell} of the injective hull of a metric space in terms of {\it extremal functions}
on the metric space. Given a metric space $X$, one first considers the set of functions $\Delta(X)$ defined by
$$
\Delta(X) := \{ f : X \to \mathbb{R} \ | \ f(x)+f(y) \geq d(x,y) \ \hbox{for all} \ x,y \in X \}.
$$

On $\Delta(X)$ we define the partial order $f \leq g$ if $f(x) \leq g(x)$ for all $x \in X$. An {\it extremal function} on $X$
is then defined to be a minimal element of the partially ordered set $(\Delta(X), \leq)$. In other words, a function $f \in \Delta(X)$
is an extremal function if, whenever $g \in \Delta(X)$ satisfies $g(x) \leq f(x)$ for all $x \in X$, then $g = f$. The set of extremal functions
is denoted by $E(X)$.

\medskip

For any $f, g \in E(X)$, the function $f - g$ is bounded, and one defines the metric on $E(X)$ by
$d(f, g) := || f - g ||_{\infty}, f,g \in E(X)$. The space of extremal functions $E(X)$ turns out to be the injective hull of $X$. Letting $d_x : X \to \mathbb{R}$ denote the distance function from a point $x$ of $X$, the map
$X \to E(X), x \mapsto d_x$ gives a canonical isometric embedding of $X$ into $E(X)$.

\medskip

The following characterization of extremal functions, formally similar to that of antipodal functions as those with discrepancy zero, will lead to a connection with antipodal spaces:
a function $f : X \to \mathbb{R}$ is an extremal function if and only if
\begin{equation} \label{extremalchar}
\sup_{y \in X} (d(x,y) - f(x) - f(y)) = 0
\end{equation}
for all $x \in X$.

\medskip

Given an extremal function $f$ on a metric space $X$ and $\epsilon > 0$, we define the following relation on $X$:
$$
x \sim_{\epsilon} y \quad \hbox{if} \quad f(x) + f(y) < d(x, y) + \epsilon.
$$
The characterization (\ref{extremalchar}) of extremal functions implies that for ahy $x \in X$ and $\epsilon > 0$, there exists $y \in X$ such that $x \sim_{\epsilon} y$.

\medskip

We first show that for an extremal function on a Gromov product space, one can define an appropriate boundary value of the function:

\medskip

\begin{prop} \label{extremalbdval} Let $X$ be a Gromov product space and let $f \in E(X)$ be an extremal function on $X$. Then for any basepoint $o \in X$, the limit
$$
\tau_o(\xi) := \lim_{x \to \xi} (d(o,x) - f(x))
$$
exists for all $\xi \in \partial_{P} X$.

\medskip

Moreover the function $\tau_o$ has discrepancy zero with respect to the visual antipodal function $\rho_o$:
$$
D_{\rho_o}(\tau_o)(\xi) = 0
$$
for all $\xi \in \partial_{P} X$.

\medskip

In particular, there exists $\alpha_o \in \mathcal{M}(\partial_{P} X)$ such that
$$
\tau_o = \log \frac{d\alpha_o}{d\rho_o}.
$$
\end{prop}

\medskip

\noindent{\bf Proof:} We will be using the following well-known fact: any extremal function is $1$-Lipschitz. Given an extremal function $f$
and a basepoint $o \in X$, we define the function $\phi_o : X \to \mathbb{R}$ by
$$
\phi_o(x) := d(o,x) - f(x) \ , \ x \in X.
$$
Since $f$ is an extremal function, we have
$$
\phi_o(x) = (d(o,x) - f(x) - f(o)) + f(o) \leq f(o)
$$
so $\phi_o$ is bounded above on $X$ by $f(o)$.

\medskip

Let $\xi \in \partial_{P} X$. By Lemma \ref{geodesicexists}, we can choose a geodesic ray $\gamma$ such that $\gamma(0) = o$ and $\gamma(t) \to \xi$
as $t \to \infty$. Then for any $s,t \geq 0$, we have, using the fact that $f$ is $1$-Lipschitz,
\begin{align*}
\phi_o(\gamma(t+s)) - \phi_o(\gamma(t)) & = (d(o, \gamma(t+s)) - d(o, \gamma(t))) - (f(\gamma(t+s)) - f(\gamma(t))) \\
                                        & \geq (d(o, \gamma(t+s)) - d(o, \gamma(t))) - d(\gamma(t+s), \gamma(t)) \\
                                        & = 0.
\end{align*}

It follows that the function $t \in [0,\infty) \mapsto \phi_o(\gamma(t))$ is monotone increasing. Since it is bounded above by $f(o)$, it follows that the limit
$$
l := \lim_{t \to \infty} \phi_o(\gamma(t))
$$
exists.

\medskip

Let $\epsilon > 0$ be given. We can choose $T > 0$ such that $\phi_o(\gamma(T)) > l - \epsilon/2$.
Note that since $\gamma$ is a geodesic, we have $(\gamma(t) | o)_{\gamma(T)} = 0$ for any $t \geq T$, hence letting $t$ tend to infinity
gives $(\xi | o)_{\gamma(T)} = 0$. Thus we can choose a neighbourhood $U$ of $\xi$ in $\widehat{X}^{P}$ such that
$(z | o)_{\gamma(T)} < \epsilon/4$ for all $z \in U \cap X$. Then for any $z \in U \cap X$ we have, again using the fact that $f$ is $1$-Lipschitz, that
\begin{align*}
\phi_o(z) & = (d(o,z) - f(z)) - (d(o, \gamma(T)) - f(\gamma(T))) + \phi_o(\gamma(T)) \\
          & > d(o,z) - d(o, \gamma(T)) - (f(z) - f(\gamma(T)) + l - \epsilon/2 \\
          & \geq d(o,z) - d(o,\gamma(T)) - d(z, \gamma(T)) + l - \epsilon/2 \\
          & = -2(z | o)_{\gamma(T)} + l - \epsilon/2 \\
          & > l - \epsilon.
\end{align*}

It follows that
$$
\liminf_{z \to \xi} \phi_o(z) \geq l.
$$

\medskip

Now let $t_n > 0$ be a sequence such that $t_n \to \infty$ as $n \to \infty$, and let $x_n = \gamma(t_n)$,
so $x_n \to \xi$, and $\phi_o(x_n) \leq l$ for all $n$ (since $\phi_o$ is monotone increasing along $\gamma$ with limit $l$).
Since $f$ is extremal, for each $n$ we can choose $y_n \in X$
such that $x_n \sim_{\epsilon/2} y_n$. Note that if we choose a geodesic segment from $x_n$ to $y_n$
and then extend it to a geodesic ray $\gamma_n$, then for any $s \geq d(x_n, y_n)$ we have
\begin{align*}
f(x_n) + f(\gamma_n(s)) & \leq f(x_n) + f(y_n) + d(y_n, \gamma_n(s)) \\
                        & < d(x_n, y_n) + d(y_n, \gamma_n(s)) + \epsilon/2 \\
                        & = d(x_n, \gamma_n(s)) + \epsilon/2,
\end{align*}
thus $x_n \sim_{\epsilon/2} \gamma_n(s)$, so replacing $y_n$ with $\gamma_n(s)$ if necessary with $s$ chosen large enough, we may assume without loss of generality
that $d(o, y_n) \geq n$ for all $n$. Then passing to a subsequence we may assume that $y_n \to \eta$ as $n \to \infty$ for some
$\eta \in \partial_{P} X$. Note that
\begin{align*}
2(x_n | y_n)_o &= d(o,x_n) + d(o, y_n) - d(x_n, y_n) \\
               & = \phi_o(x_n) + \phi_o(y_n) + f(x_n) + f(y_n) - d(x_n, y_n) \\
               & < 2f(o) + \epsilon/2
\end{align*}
for all $n$, and hence $(\xi | \eta)_o \leq f(o) + \epsilon/4 < +\infty$. It follows that $\eta \neq \xi$.

\medskip

We can choose neighbourhoods $V, W$ of $\xi, \eta$ respectively in $\widehat{X}^{P}$ such that
$$
|(z|y)_o - (\xi|\eta)_o| < \epsilon/8
$$
for all $z \in V \cap X, y \in W \cap X$. For any $z \in V \cap X$, and for $n$ large such that $x_n \in V \cap X, y_n \in W \cap X$, we have, using $\phi_o(x_n) \leq l$
and $x_n \sim_{\epsilon/2} y_n$, that
\begin{align*}
\phi_o(z) & = d(o,z) + d(o, y_n) - (f(z)+f(y_n)) - \phi_o(y_n) \\
          & \leq d(o,z) + d(o,y_n) - d(z, y_n) -\phi_o(y_n) \\
          & = 2(z|y_n)_o - \phi_o(y_n) \\
          & < 2(\xi | \eta)_o + \epsilon/4 - \phi_o(y_n) \\
          & < 2(x_n | y_n)_o + \epsilon/2 - \phi_o(y_n) \\
          & = d(o,x_n) - d(x_n, y_n) + f(y_n) + \epsilon/2 \\
          & < d(o, x_n) - f(x_n) + \epsilon/2 + \epsilon/2 \\
          & = \phi_o(x_n) + \epsilon \\
          & \leq l + \epsilon,
\end{align*}
thus
$$
\limsup_{z \to \xi} \phi_o(z) \leq l.
$$
It follows that
$$
\lim_{z \to \xi} \phi_o(z) = l,
$$
so the limit exists as required.

\medskip

Since $\phi_o$ is continuous, so is the function $\tau_o : \partial_{P} X \to \mathbb{R}$ defined by
$$
\tau_o(\xi) := \lim_{z \to \xi} \phi_o(z).
$$
Let $\rho_o$ be the visual antipodal function $\rho_o(.,.) = e^{-(.|.)_o}$. Let $\xi \in \partial_{P} X$, then for any $\eta \in \partial_{P} X$ we
have
\begin{align*}
\tau_o(\xi) + \tau_o(\eta) + \log \rho_o(\xi, \eta)^2 & = \lim_{x \to \xi, y \to \eta} \phi_o(x) + \phi_o(y) - 2(x|y)_o \\
                                                      & = \lim_{x \to \xi, y \to \eta} d(x,y) - (f(x) + f(y)) \\
                                                      & \leq 0,
\end{align*}
thus
$$
D_{\rho_o}(\tau_o)(\xi) \leq 0.
$$
On the other hand, given $\epsilon > 0$, we can choose, as we did before, sequences $x_n, y_n \in X$ such that
$x_n \sim_{\epsilon} y_n$ for all $n$, and such that $x_n \to \xi$ and $y_n \to \eta$ for some $\eta \in \partial_{P} X$
with $\eta \neq \xi$. Then
\begin{align*}
\tau_o(\xi) + \tau_o(\eta) + \log \rho_o(\xi, \eta)^2 & = \lim_{n \to \infty} \phi_o(x_n) + \phi_o(y_n) - 2(x_n|y_n)_o \\
                                                      & = \lim_{n \to \infty} d(x_n,y_n) - (f(x_n) + f(y_n)) \\
                                                      & \geq -\epsilon.
\end{align*}
It follows that $-\epsilon \leq D_{\rho_o}(\tau_o)(\xi) \leq 0$ for any $\epsilon > 0$, and hence $D_{\rho_o}(\tau_o)(\xi) = 0$
as required.
$\diamond$

\begin{prop} \label{formofextremal} Let $X$ be a Gromov product space and let $f \in E(X)$ be an extremal function.
Then there exists a unique $\alpha \in \mathcal{M}(\partial_{P} X)$ such that
$$
f(x) = d_{\mathcal{M}}( \alpha, \rho_x)
$$
for all $x \in X$.
\end{prop}

\medskip

\noindent{\bf Proof:} Let $f \in E(X)$ be an extremal function.
By the previous Proposition, for any $o \in X$, there exists $\alpha_o \in \mathcal{M}(\partial_{P} X)$ such that
$$
\log \frac{d\alpha_o}{d\rho_o}(\xi) = \lim_{z \to \xi} \phi_o(z) \ , \ \xi \in \partial_{P} X,
$$
where $\phi_o(z) = d(o,z) - f(z), z \in X$.

\medskip

We claim that in fact $\alpha_o$ does not depend on the choice of basepoint $o \in X$. If we choose a different basepoint $p \in X$,
then using Proposition \ref{visualderiv}, we have, for any $\xi \in \partial_{P} X$,
\begin{align*}
\log \frac{d\alpha_p}{d\rho_p}(\xi) & = \lim_{z \to \xi} \phi_p(z) \\
                                    & = \lim_{z \to \xi} \phi_o(z) + (d(p,z) - d(o,z)) \\
                                    & = \log \frac{d\alpha_o}{d\rho_o}(\xi) + B(p,o,\xi) \\
                                    & = \log \frac{d\alpha_o}{d\rho_o}(\xi) + \log \frac{d\rho_o}{d\rho_p}(\xi) \\
                                    & = \log \frac{d\alpha_o}{d\rho_p}(\xi),
\end{align*}
which implies $\log \frac{d\alpha_p}{d\alpha_o} \equiv 0$, and so $\alpha_p = \alpha_o$.

\medskip

Thus there exists $\alpha \in \mathcal{M}(\partial_{P} X)$ such that
$$
\log \frac{d\alpha}{d\rho_o}(\xi) = \lim_{z \to \xi} \phi_o(z) \ , \ \xi \in \partial_{P} X,
$$
for any $o \in X$. As observed in the proof of the previous Proposition, $\phi_o$ is bounded above by $f(o)$, and hence
$$
\log \frac{d\alpha}{d\rho_o}(\xi) \leq f(o)
$$
for all $\xi \in \partial_{P} X$, and thus
$$
d_{\mathcal{M}}(\alpha, \rho_o) \leq f(o)
$$
for all $o \in X$. Since the visual embedding $o \mapsto \rho_o$ is an isometric embedding, the function $g : X \to \mathbb{R}$ defined by
$g(x) := d_{\mathcal{M}}(\alpha, \rho_x)$ satisfies $g(x) + g(y) \geq d(x,y)$ for all $x,y \in X$, so $g \in \Delta(X)$. Moreover as shown above, $g(x) \leq f(x)$ for all $x \in X$.
Since $f$ is extremal, this implies $f \equiv g$ as required. Finally, the point $\alpha \in \mathcal{M}(\partial_{P} X)$ is unique,
since if $\beta \in \mathcal{M}(\partial_{P} X)$ satisfies
$$
d_{\mathcal{M}}(\alpha, \rho_x) = d_{\mathcal{M}}(\beta, \rho_x)
$$
for all $x \in X$, then fixing a basepoint $o \in X$, given $\xi \in X$ we can subtract $d_{\mathcal{M}}(\rho_o, \rho_x)$ from
both sides of the above equality and let $x$ tend to $\xi$. Using Propositions \ref{isomembed} and \ref{mzbuse}, we obtain
$$
\log \frac{d\rho_o}{d\alpha}(\xi) = \log \frac{d\rho_o}{d\beta}(\xi)
$$
for all $\xi \in \partial_{P} X$, thus $\log \frac{d\beta}{d\alpha} \equiv 0$, and so $\alpha = \beta$. $\diamond$

\medskip

\begin{lemma} \label{drhoextremal} Let $X$ be a Gromov product space, and let $\alpha \in \mathcal{M}(\partial_{P} X)$. Then the
function $d_{\alpha} : X \to \mathbb{R}$ defined by
$$
d_{\alpha}(x) := d_{\mathcal{M}}( \alpha, \rho_x) \ , \ x \in X,
$$
is an extremal function.
\end{lemma}

\medskip

\noindent{\bf Proof:} As observed earlier, $d_{\alpha} \in \Delta(X)$.

\medskip

Given $x \in X$, let $\xi \in \hbox{argmax} \frac{d\alpha}{d\rho_x}$. Let $y_n \in X$ be a sequence converging to $\xi$.
Then by Proposition \ref{mzbuse}, we have
\begin{align*}
\lim_{n \to \infty} d(x,y_n) - (d_{\alpha}(x) + d_{\alpha}(y_n)) & = \lim_{n \to \infty} (d_{\mathcal{M}}(\rho_x, \rho_{y_n}) - d_{\mathcal{M}}(\alpha, \rho_{y_n})) +
d_{\mathcal{M}}(\alpha, \rho_x) \\
                                                                 & = \log \frac{d\alpha}{d\rho_x}(\xi) - d_{\mathcal{M}}(\alpha, \rho_x) \\
                                                                 & = 0,
\end{align*}
so it follows from the characterization (\ref{extremalchar}) of extremal functions that $d_{\alpha}$ is extremal. $\diamond$

\medskip

We can now prove that the injective hull of a Gromov product space $X$ is isometric to the Moebius space $\mathcal{M}(\partial_{P} X)$, in particular this gives a proof of Theorem \ref{mainthm9} from the Introduction:

\medskip

\begin{theorem} \label{gromovinjectivehull} Let $X$ be a Gromov product space. Then the map
\begin{align*}
\Phi : \mathcal{M}(\partial_{P} X) & \to E(X) \\
               \alpha              & \mapsto d_{\alpha}
\end{align*}
is an isometry.

\medskip

In particular, the injective hull of $X$ is isometric to the Moebius space $\mathcal{M}(\partial_{P} X)$.
\end{theorem}

\medskip

\noindent{\bf Proof:} Let $\alpha, \beta \in \mathcal{M}(\partial_{P} X)$. Then for any $x \in X$ we have
$|d_{\alpha}(x) - d_{\beta}(x)| \leq d_{\mathcal{M}}(\alpha, \beta)$ by the triangle inequality. On the other
hand, if we take a $\xi \in \hbox{argmax} \frac{d\beta}{d\alpha}$ and let $x_n \in X$ be a sequence converging to $\xi$,
then by Proposition \ref{mzbuse} we have
\begin{align*}
\lim_{n \to \infty} d_{\alpha}(x_n) - d_{\beta}(x_n) & = \lim_{n \to \infty} d_{\mathcal{M}}(\alpha, \rho_{x_n}) - d_{\mathcal{M}}(\beta, \rho_{x_n}) \\
                                                     & = \log \frac{d\beta}{d\alpha}(\xi) \\
                                                     & = d_{\mathcal{M}}(\alpha, \beta).
\end{align*}
It follows that
$$
d(\Phi(\alpha), \Phi(\beta)) = || d_{\alpha} - d_{\beta} ||_{\infty} = d_{\mathcal{M}}(\alpha, \beta),
$$
and so $\Phi$ is an isometric embedding. By Proposition \ref{formofextremal}, the map $\Phi$ is surjective, hence is an isometry. $\diamond$

\medskip

As an immediate corollary we obtain a proof of Theorem \ref{mainthm10} from the Introduction, giving the equivalence of maximality and injectivity for Gromov product spaces:

\medskip

\noindent{\bf Proof of Theorem \ref{mainthm10}:} Let $X$ be a Gromov product space. By the previous Theorem, the injective hull $E(X)$
is isometric to the Moebius space $\mathcal{M}(\partial_{P} X)$. Now $X$ is injective if and only if $X$ is isometric to $E(X)$,
while by Theorem \ref{mainthm7}, $X$ is maximal if and only if $X$ is isometric to $\mathcal{M}(\partial_{P}X)$. Since $E(X)$ and $\MM(\partial_{P} X)$ are isometric, it follows
that $X$ is injective if and only if $X$ is maximal. $\diamond$

\medskip

\section{Gromov hyperbolicity of $\mathcal{M}(Z)$ for quasi-metric $Z$.}

\medskip

Let $Z$ be an antipodal space. In this section we show that if $Z$ is a quasi-metric space, i.e.
some $\rho_0 \in \mathcal{M}(Z)$ is a quasi-metric, then the Moebius space $\mathcal{M}(Z)$ is
Gromov hyperbolic.

\medskip

We will call an antipodal function $\rho \in \mathcal{M}(Z)$ a {\it quasi-metric}, more
specifically a $K$-{\it quasi-metric} if it satisfies the inequality
$$
\rho(\xi, \eta) \leq K \max( \rho(\xi, \zeta), \rho(\zeta, \eta) )
$$
for all $\xi, \eta, \zeta \in Z$, for some fixed $K \geq 1$. We will refer to the above inequality
as the "quasi-metric inequality". A well-known construction of Frink shows that if $\rho$ is a $K$-quasi-metric
with $K \leq 2$, then there exists a metric $\alpha$ on $Z$ such that $\frac{1}{2K} \rho \leq \alpha \leq \rho$
(for a proof see \cite{schroederquasi}).

\medskip

It follows that for any quasi-metric $\rho$, if $\epsilon > 0$ is sufficiently small, then there is a metric
$\alpha$ such that
$$
\frac{1}{C} \rho^{\epsilon} \leq \alpha \leq C \rho^{\epsilon}
$$
for some $C \geq 1$.

\medskip

We note that for any two $\rho_1, \rho_2 \in \mathcal{M}(Z)$, by the Geometric Mean-Value Theorem
we have $e^{-d} \rho_2 \leq \rho_1 \leq e^{d} \rho_2$, where $d = d_{\mathcal{M}}(\rho_1, \rho_2)$.
From this it follows that if some $\rho_1 \in \mathcal{M}(Z)$ is a quasi-metric, then all
$\rho \in \mathcal{M}(Z)$ are quasi-metrics. We will say that the antipodal space
$Z$ is quasi-metric if all elements of $\mathcal{M}(Z)$
are quasi-metrics.

\medskip

The following Lemma will be useful:

\medskip

\begin{lemma} \label{derivcont} Suppose that $\rho_0 \in \mathcal{M}(Z)$ is a
quasi-metric. Then there exist $\delta, M > 0$ only depending on $\rho_0$, such that
for all $\rho_1 \in \mathcal{M}(Z)$, if $\xi_1 \in \hbox{argmax} \ \frac{d\rho_1}{d\rho_0}$,
and $\xi_2 \in Z$ is such that $\rho_1(\xi_1, \xi_2) < \delta$, then
$$
\log \frac{d\rho_1}{d\rho_0}(\xi_2) \geq \log \frac{d\rho_1}{d\rho_0}(\xi_1) - M = d_{\mathcal{M}}(\rho_0, \rho_1) - M.
$$
\end{lemma}

\medskip

\noindent{\bf Proof:} Let $\alpha$ be a metric on $Z$ such that $\frac{1}{C} {\rho_0}^{\epsilon} \leq \alpha \leq C {\rho_0}^{\epsilon}$
for some $C \geq 1$ and $\epsilon > 0$. We let
$$
\delta = \left( \frac{1}{2C^2} \right)^{1/\epsilon}.
$$
Now let $\rho_1 \in \mathcal{M}(Z)$, and
let $\tau = \log \frac{d\rho_1}{d\rho_0}$. Let $\xi_1 \in \hbox{argmax} \ \frac{d\rho_1}{d\rho_0}$,
and let $\xi_2 \in Z$ be such that $\rho_1(\xi_1, \xi_2) < \delta$. If $\xi_2 = \xi_1$ then there is nothing to prove,
so we may assume $\xi_2 \neq \xi_1$.

\medskip

Let $\eta_2 \in Z$ be such that $\rho_1(\xi_2, \eta_2) = 1$. Now noting that $\tau(\eta_2) \leq \tau(\xi_1)$ and
$\rho_1(\xi_1, \xi_2)^{\epsilon} \cdot 2C^2 < 1$ by our choice of $\delta$, using the
Geometric Mean-Value Theorem we have
\begin{align*}
\alpha(\eta_2, \xi_2) & \geq \frac{1}{C} \rho_0(\eta_2, \xi_2)^{\epsilon} \\
                      & = \frac{1}{C} \exp(- \epsilon \tau(\eta_2)/2) \cdot \exp(- \epsilon \tau(\xi_2)/2) \cdot \rho_1(\eta_2, \xi_2)^{\epsilon} \\
                      & \geq \frac{1}{C} \exp(- \epsilon \tau(\xi_1)/2) \cdot \exp(- \epsilon \tau(\xi_2)/2) \\
                      & > \frac{1}{C} \exp(- \epsilon \tau(\xi_1)/2) \cdot \exp(- \epsilon \tau(\xi_2)/2) \cdot \rho_1(\xi_1, \xi_2)^{\epsilon} \cdot 2C^2 \\
                      & = \frac{1}{C} \rho_0(\xi_1, \xi_2)^{\epsilon} \cdot 2C^2 \\
                      & \geq 2 \alpha(\xi_1, \xi_2).
\end{align*}
From this inequality and the triangle inequality for $\alpha$ we get
$$
\alpha(\eta_2, \xi_1) \geq \alpha(\eta_2, \xi_2) - \alpha(\xi_1, \xi_2) \geq \alpha(\xi_1, \xi_2) > 0,
$$
and hence
\begin{align*}
\frac{\alpha(\eta_2, \xi_2)}{\alpha(\eta_2, \xi_1)} & \leq \frac{ \alpha(\eta_2, \xi_2) } { \alpha(\eta_2, \xi_2) - \alpha(\xi_1, \xi_2) } \\
                                                    & = \frac{1}{1 - \alpha(\xi_1, \xi_2)/\alpha(\eta_2, \xi_2)} \\
                                                    & \leq \frac{1}{1 - 1/2} = 2.
\end{align*}
Now using the Geometric Mean-Value Theorem we have
\begin{align*}
\tau(\eta_2) + \tau(\xi_1) + \log \rho_0(\eta_2, \xi_1)^2 & = \log \rho_1(\eta_2, \xi_1)^2 \\
                                                          & \leq 0 = \log \rho_1(\eta_2, \xi_2)^2 \\
                                                          & = \tau(\eta_2) + \tau(\xi_2) + \log \rho_0(\eta_2, \xi_2)^2,
\end{align*}
thus
\begin{align*}
\tau(\xi_1) & \leq \tau(\xi_2) + \log \frac{\rho_0(\eta_2, \xi_2)^2}{\rho_0(\eta_2, \xi_1)^2} \\
            & = \tau(\xi_2) + \frac{2}{\epsilon} \log \frac{\rho_0(\eta_2, \xi_2)^{\epsilon}}{\rho_0(\eta_2, \xi_1)^{\epsilon}} \\
            & \leq \tau(\xi_2) + \frac{2}{\epsilon} \log \left( C^2 \cdot \frac{\alpha(\eta_2, \xi_2)}{\alpha(\eta_2, \xi_1)} \right) \\
            & \leq \tau(\xi_2) + \frac{2}{\epsilon} \log(2C^2),
\end{align*}
and so the Lemma follows by taking $M = (2/\epsilon) \log(2C^2)$. $\diamond$

\medskip

We will need also the following Lemmas.

\medskip

\begin{lemma} \label{gipbyrhoest} Let $\rho_0, \alpha, \beta \in \mathcal{M}(Z)$. Let $\xi \in \hbox{argmax} \ \frac{d\alpha}{d\rho_0},
\eta \in \hbox{argmax} \ \frac{d\beta}{d\rho_0}$, and suppose that $\xi \neq \eta$. Let $\tilde{\eta} \in \hbox{argmax} \ \frac{d\beta}{d\alpha}$.
Then
$$
\frac{\exp\left( - (\alpha | \beta)_{\rho_0} \right)}{\rho_0(\xi, \eta)} \leq \frac{1}{\alpha(\xi, \tilde{\eta}) \beta(\xi, \eta)}
$$
(with the usual convention that the right-hand side is $+\infty$ if $\xi = \tilde{\eta}$).
\end{lemma}

\medskip

\noindent{\bf Proof:} By Proposition \ref{geodcomp}, we can extend geodesic segments joining $\rho_0$ to $\alpha, \beta$ to
geodesic rays $\kappa, \kappa' : [0, \infty) \to \mathcal{M}(Z)$ such that $\kappa(0) = \kappa'(0) = \rho_0$ and
$\kappa(t) \to \xi$, $\kappa'(t) \to \eta$, as $t \to +\infty$. Moreover from Proposition \ref{geodcomp} we have
$\xi \in \hbox{argmax} \ \frac{d\kappa(t)}{d\rho_0}, \eta \in \hbox{argmax} \ \frac{d\kappa'(t)}{d\rho_0}$ for all $t > 0$,
so letting $t_0 := \max(d_{\mathcal{M}}( \rho_0, \alpha), d_{\mathcal{M}}( \rho_0, \beta))$, for $t > t_0$ we can apply
Lemma \ref{trieq} (to the obvious triples of points) to get $\xi \in \hbox{argmax} \ \frac{d\kappa(t)}{d\alpha}, \eta \in \hbox{argmax} \ \frac{d\kappa'(t)}{d\beta}$. Then we have the elementary identity
$$
( \kappa(t) | \kappa'(t))_{\rho_0} - ( \alpha | \beta)_{\rho_0} = ( \kappa(t) | \beta)_{\alpha} + ( \kappa(t) | \kappa'(t))_{\beta}
$$
for $t > t_0$. Since $\xi \in \hbox{argmax} \ \frac{d\kappa(t)}{d\alpha}$ and
$\tilde{\eta} \in \hbox{argmax} \ \frac{d\beta}{d\alpha}$, by Lemma \ref{giproupper} we get
$$
( \kappa(t) | \kappa'(t))_{\rho_0} - ( \alpha | \beta)_{\rho_0} \leq \log \frac{1}{\alpha( \xi, \tilde{\eta})} + ( \kappa(t) | \kappa'(t))_{\beta}
$$
for $t > t_0$.
Now letting $t \to \infty$ in the above inequality, since $\xi \neq \eta$ we can apply Theorem \ref{gromovip} to obtain
$$
\log \frac{1}{ \rho_0(\xi, \eta) } - ( \alpha | \beta)_{\rho_0} \leq \log \frac{1}{\alpha( \xi, \tilde{\eta})} + \log \frac{1}{\beta(\xi, \eta)},
$$
from which the Lemma follows by applying the exponential to both sides above. $\diamond$

\medskip

\begin{lemma} \label{xineqeta} Let $\rho_0 \in \mathcal{M}(Z)$ and suppose that $\rho_0$ is a quasi-metric.
Suppose $\alpha_n, \beta_n, \rho_n \in \mathcal{M}(Z)$ and $C_n > 0$
are sequences such that $C_n \to +\infty$ as $n \to \infty$ and
\begin{equation} \label{hypfail}
(\alpha_n | \beta_n)_{\rho_0} \leq \min\left( (\alpha_n | \rho_n)_{\rho_0}, (\rho_n | \beta_n)_{\rho_0} \right) - C_n
\end{equation}
for all $n$. For all $n$, let $\xi_n, \eta_n \in Z$ be such that
$\xi_n \in \hbox{argmax} \ \frac{d\alpha}{d\rho_0}, \eta_n \in \hbox{argmax} \ \frac{d\beta}{d\rho_0}$.
Then $\xi_n \neq \eta_n$ for all $n$ sufficiently large.
 \end{lemma}

\medskip

\noindent{\bf Proof:} Suppose not, so $\xi_n = \eta_n$ for infinitely many $n$, then passing to a
subsequence we may assume $\xi_n = \eta_n$ for all $n$.

\medskip

Let $s_n = d_{\mathcal{M}}( \rho_0, \alpha_n), t_n = d_{\mathcal{M}}( \rho_0, \beta_n)$ and $r_n = d_{\mathcal{M}}( \alpha_n, \beta_n)$.
Interchanging $\alpha_n, \beta_n$
if necessary, we may assume $s_n \leq t_n$ for all $n$. Let $\widetilde{\eta_n} \in \hbox{argmax} \ \frac{d\beta_n}{d\alpha_n}$.
Applying the elementary inequality $(y|z)_x \leq \min(d(x,y), d(x,z))$, the inequality (\ref{hypfail}) gives us
$$
\frac{1}{2} ( s_n + t_n - r_n) \leq s_n - C_n,
$$
and hence
\begin{align*}
\log \frac{d\beta_n}{d\alpha_n}(\widetilde{\eta_n}) = r_n & \geq (t_n - s_n) + 2C_n \\
                                                          & = \log \frac{d\beta_n}{d\alpha_n}(\xi_n) + 2C_n,
\end{align*}
where we have used $\xi_n = \eta_n \in \hbox{argmax} \ \frac{d\alpha_n}{d\rho_0} \cap \hbox{argmax} \ \frac{d\beta_n}{d\rho_0}$,
so that $\frac{d\alpha_n}{d\rho_0}(\xi_n) = e^{s_n}, \frac{d\beta_n}{d\rho_0}(\xi_n) = e^{t_n}$. Now using the fact
that $\xi_n \in \hbox{argmax} \ \frac{d\beta_n}{d\rho_0}$, we get
\begin{align*}
\log \frac{d\beta_n}{d\rho_0}(\xi_n) & \geq \log \frac{d\beta_n}{d\rho_0}(\widetilde{\eta_n}) \\
                                     & = \log \frac{d\beta_n}{d\alpha_n}(\widetilde{\eta_n}) + \log \frac{d\alpha_n}{d\rho_0}(\widetilde{\eta_n}) \\
                                     & \geq \log \frac{d\beta_n}{d\alpha_n}(\xi_n) + 2C_n + \log \frac{d\alpha_n}{d\rho_0}(\widetilde{\eta_n}),
\end{align*}
thus
\begin{align*}
\log \frac{d\alpha_n}{d\rho_0}(\widetilde{\eta_n}) & \leq \log \frac{d\beta_n}{d\rho_0}(\xi_n) - \log \frac{d\beta_n}{d\alpha_n}(\xi_n) - 2C_n \\
                                                   & =  \log \frac{d\alpha_n}{d\rho_0}(\xi_n) - 2C_n.
\end{align*}

Since $C_n \to +\infty$ and $\rho_0$ is a quasi-metric, by Lemma \ref{derivcont} the above inequality implies that we must have
$\alpha_n( \xi_n, \widetilde{\eta_n} ) \geq \delta$ for all $n$ large enough, where $\delta > 0$ is the constant given by
Lemma \ref{derivcont}.

\medskip

Fix an $n$ large enough so that $\alpha_n( \xi_n, \widetilde{\eta_n} ) \geq \delta$. For each $t > 0$, by Lemma \ref{raypt}
we can choose a point $\rho_t \in \mathcal{M}(Z)$ such that $\xi_n \in \hbox{argmax} \ \frac{d\rho_t}{d\alpha_n}$ and
$d_{\mathcal{M}}( \alpha_n, \rho_t) = t$. Then by Lemma \ref{giproupper}, we have
\begin{align*}
2(\rho_t | \beta_n)_{\alpha_n} & = d_{\mathcal{M}}( \alpha_n, \rho_t) + d_{\mathcal{M}}( \alpha_n, \beta_n ) - d_{\mathcal{M}}( \beta_n, \rho_t) \\
                               & \leq \log \frac{1}{\alpha_n( \xi_n, \widetilde{\eta_n})^2} \\
                               & \leq \log \frac{1}{\delta^2},
\end{align*}
hence
$$
d_{\mathcal{M}}( \alpha_n, \rho_t) - d_{\mathcal{M}}( \beta_n, \rho_t) \leq -r_n + \log \frac{1}{\delta^2}.
$$

Now letting $t \to \infty$ in the above inequality (with $n$ fixed), by Lemma \ref{maximaconv} we know
$\rho_t \to \xi_n$ in $\widehat{\mathcal{M}(Z)}$ as $t \to \infty$, and so applying Proposition \ref{mzbuse}
we get
\begin{align*}
t_n - s_n = \log \frac{d\beta_n}{d\alpha_n}(\xi_n) & \leq -r_n + \log \frac{1}{\delta^2} \\
                                                   & \leq -( (t_n - s_n) + 2C_n) + \log \frac{1}{\delta^2},
\end{align*}
hence
$$
0 \leq 2(t_n - s_n) \leq -2C_n + \log \frac{1}{\delta^2},
$$
which gives us a contradiction for $n$ large since $C_n \to +\infty$ as $n \to \infty$. $\diamond$

\medskip

We can now prove the main theorem of this section:

\medskip

\begin{theorem} \label{gromovhyp} Let $Z$ be an antipodal space and suppose that some (and hence every) $\rho_0 \in \mathcal{M}(Z)$ is
a quasi-metric. Then the space $\mathcal{M}(Z)$ is Gromov hyperbolic.
\end{theorem}

\medskip

\noindent{\bf Proof:} Given $\rho_0 \in \mathcal{M}(Z)$ such that $\rho_0$ is a quasi-metric, we fix $\delta, M > 0$ as given by
Lemma \ref{derivcont}.

\medskip

To show that $\mathcal{M}(Z)$ is Gromov hyperbolic, as remarked in section 5 it
suffices to show that there is a constant $C > 0$ such that
$$
(\alpha | \beta)_{\rho_0} \geq \min\left( (\alpha | \rho)_{\rho_0}, (\rho | \beta)_{\rho_0} \right) - C
$$
for all $\alpha, \beta, \rho \in \mathcal{M}(Z)$ (in that case as remarked earlier $\mathcal{M}(Z)$ will be $2C$-hyperbolic).
Suppose not, then there are sequences $\alpha_n, \beta_n, \rho_n \in \mathcal{M}(Z)$ and $C_n > 0$ such that
\begin{equation} \label{tendinf}
( \alpha_n | \beta_n )_{\rho_0} \leq \min\left( ( \alpha_n | \rho_n )_{\rho_0}, ( \rho_n | \beta_n )_{\rho_0} \right) -  C_n
\end{equation}
for all $n$ and such that $C_n \to +\infty$ as $n \to \infty$.


\medskip

Let $\xi_n \in \hbox{argmax} \ \frac{d\alpha_n}{d\rho_0},
\eta_n \in \hbox{argmax} \ \frac{d\beta_n}{d\rho_0}, \lambda_n \in \hbox{argmax} \ \frac{d\rho_n}{d\rho_0}$ and $\widetilde{\eta_n} \in \hbox{argmax} \ \frac{d\beta_n}{d\alpha_n}$.
We then note that by Lemma \ref{xineqeta}, our hypotheses imply that $\xi_n \neq \eta_n$ for all $n$ large enough, and so we may as well
assume that $\xi_n \neq \eta_n$ for all $n$.

\medskip

Applying Lemma \ref{giproupper} and using the fact that $\rho_0$ is a $K$-quasi-metric for some $K \geq 1$, we have
\begin{align*}
\min\left( ( \alpha_n | \rho_n )_{\rho_0}, ( \rho_n | \beta_n )_{\rho_0} \right) & \leq \log \left( \min\left( \frac{1}{\rho_0(\xi_n, \lambda_n)},
                                                                                                               \frac{1}{\rho_0(\lambda_n, \eta_n)} \right) \right) \\
                                                                                 & = \log \left( \frac{1}{\max( \rho_0(\xi_n, \lambda_n), \rho_0(\lambda_n, \eta_n) )} \right) \\
                                                                                 & \leq \log \left( \frac{K}{\rho_0(\xi_n, \eta_n)} \right),
\end{align*}
so the above together with (\ref{tendinf}) gives, after exponentiating, that
\begin{equation} \label{ratioest}
\frac{\exp\left( - (\alpha_n | \beta_n)_{\rho_0} \right)}{\rho_0(\xi_n, \eta_n)} \geq \frac{e^{C_n}}{K}.
\end{equation}

\medskip

Let $\zeta_n \in Z$ be such that $\rho_0( \xi_n, \zeta_n) = 1$, so that
$\zeta_n \in \hbox{argmin} \ \frac{d\alpha_n}{d\rho_0} = \hbox{argmax} \ \frac{d\rho_0}{d\alpha_n}$ by Lemma \ref{maxmin}.
Let $s_n = d_{\mathcal{M}}( \rho_0, \alpha_n), t_n = d_{\mathcal{M}}( \rho_0, \beta_n), r_n = d_{\mathcal{M}}( \alpha_n, \beta_n)$. Then by
(\ref{tendinf}) we have
$$
2(\alpha_n | \beta_n)_{\rho_0} = s_n + t_n - r_n \leq 2 \min(s_n, t_n) - 2C_n,
$$
and so
\begin{align*}
2(\rho_0 | \beta_n)_{\alpha_n} = s_n + r_n - t_n & \geq (s_n - t_n) + (s_n + t_n - 2 \min(s_n, t_n) + 2C_n) \\
                                                 & \geq (s_n - t_n) + (s_n + t_n - 2 s_n + 2C_n) = 2C_n,
\end{align*}
and by Lemma \ref{giproupper} we have
$$
( \rho_0 | \beta_n)_{\alpha_n} \leq \log \frac{1}{\alpha_n( \zeta_n, \widetilde{\eta_n} )}.
$$
Combining these inequalities gives
\begin{align*}
\alpha_n( \zeta_n, \widetilde{\eta_n} ) & \leq \exp\left( - ( \rho_0 | \beta_n)_{\alpha_n} \right) \\
                                        & \leq e^{-C_n} \\
                                        & \to 0 \ \hbox{as} \ n \to \infty.
\end{align*}

We now make the following claim:

\medskip

\noindent{\bf Claim.} There exists $\kappa > 0$ such that $\alpha_n(\xi_n, \widetilde{\eta_n}) \geq \kappa$ for all $n$ large enough.

\medskip

\noindent{\bf Proof of Claim:} Suppose not, then passing to a subsequence we may assume $\alpha_n(\xi_n, \widetilde{\eta_n}) \to 0$ as $n \to \infty$.
Then for all $n$ large we have $\alpha_n(\xi_n, \widetilde{\eta_n}) < \delta$, where $\delta > 0$ is the constant given by Lemma \ref{derivcont}.
By Lemma \ref{derivcont} this implies that for all $n$ large
$$
\log \frac{d\alpha_n}{d\rho_0}(\widetilde{\eta_n}) \geq \log \frac{d\alpha_n}{d\rho_0}(\xi_n) - M = s_n - M.
$$
Using the Geometric Mean-Value Theorem this gives us
\begin{align*}
\rho_0(\zeta_n, \widetilde{\eta_n})^2 & = \frac{d\rho_0}{d\alpha_n}(\zeta_n) \cdot \frac{d\rho_0}{d\alpha_n}(\widetilde{\eta_n}) \cdot \alpha_n(\zeta_n, \widetilde{\eta_n})^2 \\
                                      & \leq e^{s_n} \cdot e^{-s_n + M} \cdot e^{-2C_n} \\
                                      & = e^{M - 2C_n} \\
                                      & \to 0 \ \hbox{as} \ n \to \infty.
\end{align*}

Now note that the fact that $\rho_0$ is a $K$-quasi-metric implies that
\begin{align*}
\rho_0(\xi_n, \widetilde{\eta_n}) & \geq \frac{ \rho_0(\xi_n, \zeta_n) }{K} - \rho_0(\zeta_n, \widetilde{\eta_n}) \\
                                  & = \frac{1}{K} - \rho_0(\zeta_n, \widetilde{\eta_n}) \\
                                  & \geq \frac{1}{2K} \ \hbox{for all } \ n \ \hbox{large enough}.
\end{align*}

Again using the Geometric Mean-Value Theorem we then get, for all $n$ large enough, that
\begin{align*}
\alpha_n(\xi_n, \widetilde{\eta_n})^2 & =  \frac{d\alpha_n}{d\rho_0}(\xi_n) \cdot \frac{d\alpha_n}{d\rho_0}(\widetilde{\eta_n}) \cdot \rho_0(\xi_n, \widetilde{\eta_n})^2 \\
                                      & \geq e^{s_n} \cdot e^{-s_n} \cdot \left( \frac{1}{2K} \right)^2 \\
                                      & = \left( \frac{1}{2K} \right)^2,
\end{align*}
which contradicts our assumption that $\alpha_n(\xi_n, \widetilde{\eta_n}) \to 0$ as $n \to \infty$.
This finishes the proof of the Claim.

\medskip

Now, since $\xi_n \neq \eta_n$ for all $n$, we can apply Lemma \ref{gipbyrhoest} and the Claim to get
\begin{align*}
\frac{\exp\left( - (\alpha_n | \beta_n)_{\rho_0} \right)}{\rho_0(\xi_n, \eta_n)} & \leq \frac{1}{\alpha_n(\xi_n, \widetilde{\eta_n}) \beta_n(\xi_n, \eta_n)} \\
                                                                                 & \leq \frac{1}{\kappa} \frac{1}{ \beta_n(\xi_n, \eta_n)},
\end{align*}
and the above together with the inequality (\ref{ratioest}) gives
\begin{align*}
\beta_n(\xi_n, \eta_n) & \leq \frac{1}{\kappa} \cdot K \cdot e^{-C_n} \\
                       & \to 0 \ \hbox{as} \ n \to \infty.
\end{align*}
So for all $n$ large enough we have $\beta_n(\xi_n, \eta_n) < \delta$, with $\delta > 0$ being as before the constant given by
Lemma \ref{derivcont}, therefore by Lemma \ref{derivcont} we have
\begin{equation} \label{betaest}
\log \frac{d\beta_n}{d\rho_0}(\xi_n) \geq \log \frac{d\beta_n}{d\rho_0}(\eta_n) - M = t_n - M,
\end{equation}
for all $n$ large enough. Now using the inequality
$$
\frac{\exp\left( - (\alpha_n | \beta_n)_{\rho_0} \right)}{\rho_0(\xi_n, \eta_n)} \leq \frac{1}{\kappa} \frac{1}{ \beta_n(\xi_n, \eta_n)}
$$
which we proved above, we get
$$
\exp\left( \alpha_n | \beta_n \right)_{\rho_0} \geq \kappa \cdot \frac{\beta_n(\xi_n, \eta_n)}{\rho_0(\xi_n, \eta_n)},
$$
and hence, using the Geometric Mean-Value Theorem, the inequality (\ref{betaest}) above and the fact that $C_n \to +\infty$, we have
\begin{align*}
2( \alpha_n | \beta_n)_{\rho_0} & \geq \log \frac{\beta_n(\xi_n, \eta_n)^2}{\rho_0(\xi_n, \eta_n)^2} + \log \kappa^2 \\
                                & = \log \frac{d\beta_n}{d\rho_0}(\xi_n) + \log \frac{d\beta_n}{d\rho_0}(\eta_n) + \log \kappa^2 \\
                                & \geq (t_n-M) + t_n + \log \kappa^2 \\
                                & \geq 2(\rho_n | \beta_n)_{\rho_0} - M + \log \kappa^2 \\
                                & \geq 2 \min( (\alpha_n | \rho_n)_{\rho_0}, (\rho_n | \beta_n)_{\rho_0} ) - M + \log \kappa^2 \\
                                & > 2 \min( (\alpha_n | \rho_n)_{\rho_0}, (\rho_n | \beta_n)_{\rho_0} ) - 2C_n, \ \ \hbox{for all} \ n \ \hbox{large enough},
\end{align*}
which contradicts our hypothesis (\ref{tendinf}). This finishes the proof of the Theorem. $\diamond$

\medskip

Conversely, we have:

\medskip

\begin{theorem} \label{hypimpqm} Let $Z$ be an antipodal space. If $\mathcal{M}(Z)$ is Gromov hyperbolic, say $\delta$-hyperbolic,
then every $\rho_1 \in \mathcal{M}(Z)$ is a $K$-quasi-metric, where $K = e^{\delta}$.
\end{theorem}

\medskip

\noindent{\bf Proof:} Let $\rho_1 \in \mathcal{M}(Z)$, then we are given that
$$
( \alpha | \beta)_{\rho_1} \geq \min( (\alpha | \rho)_{\rho_1}, ( \rho | \beta)_{\rho_1} ) - \delta
$$
for all $\alpha, \beta, \rho \in \mathcal{M}(Z)$, or equivalently
$$
\exp\left( - ( \alpha | \beta)_{\rho_1} \right) \leq K \max\left( \exp\left( - ( \alpha | \rho)_{\rho_1} \right), \exp\left( - ( \rho | \beta)_{\rho_1} \right) \right),
$$
where $K = e^{\delta}$. Let $\xi, \eta, \zeta \in Z$ be three distinct points, then
letting $\alpha, \beta, \rho$ tend to $\xi, \eta, \zeta$ in $\widehat{\mathcal{M}(Z)}$ respectively in the above inequality
and using Theorem \ref{gromovip} gives
$$
\rho_1(\xi, \eta) \leq K \max( \rho_1(\xi, \zeta), \rho_1(\zeta, \eta) )
$$
for all $\xi, \eta, \zeta$ distinct points in $Z$, and the above inequality holds trivially if any two of these three points are equal,
thus $\rho_1$ is a $K$-quasi-metric. $\diamond$

\medskip

Combining the above two theorems gives immediately

\medskip

\begin{theorem} \label{hypiffqm} Let $Z$ be an antipodal space. Then $\mathcal{M}(Z)$ is Gromov hyperbolic
if and only if some $\rho_0 \in \mathcal{M}(Z)$ is quasi-metric. In this case, in fact all $\rho \in \mathcal{M}(Z)$
are $K$-quasi-metric, where $K = e^{\delta}$ and $\delta \geq 0$ is such that $\mathcal{M}(Z)$ is $\delta$-hyperbolic.
In particular if $\mathcal{M}(Z)$ is a tree then $K = 1$ and all $\rho \in \mathcal{M}(Z)$ are ultrametrics.
\end{theorem}

\medskip

Note that the last statement of the above Theorem follows from the fact that trees are $0$-hyperbolic.
We observe also that from the above Theorem it follows that if $X$ is a proper, geodesically complete, boundary
continuous Gromov hyperbolic space, then the space $\mathcal{M}(\partial X)$ is Gromov hyperbolic,
since all the elements $\rho_x \in \mathcal{M}(\partial X), x \in X,$ are quasi-metrics.

\medskip

\begin{lemma} \label{gpbdryghyp} Let $X$ be a Gromov product space. If $X$ is Gromov hyperbolic, then the Gromov product compactification of $X$ is equivalent to the visual compactification, and $X$ is boundary continuous.

\medskip

Moreover, the identification $f : \partial_{P} X \to \partial P$ between the Gromov product boundary and the visual boundary is a Moebius homeomorphism.
\end{lemma}

\medskip

\noindent{\bf Proof:} Let $X$ be a Gromov product space and suppose that 
$X$ is Gromov hyperbolic. We fix a basepoint $o \in X$. 

\medskip

Let $\xi \in \partial_{P} X$ and let $x_n \in X$ be a sequemce converging to 
$\xi$ in $\widehat{X}^{P}$. Let $\eta, \eta' \in \partial X$ be any two limit points of the sequence in the visual compactification, and let $y_k, z_k \in X$ be subsequences of $(x_n)$ such that $y_k \to \eta, z_k \to \eta'$ as $k \to \infty$. 
If $\eta \neq \eta'$, then since $X$ is Gromov hyperbolic there exists $M > 0$ such that $(y_k | z_k)_o \leq M$ for all $k$. On the other hand, since $X$ is a Gromov product space, we have $(y_k | z_k)_o \to (\xi | \xi)_o = +\infty$ as 
$k \to \infty$, a contradiction. This implies that $\eta = \eta'$. It follows that  the sequence $(x_n)$ has a unique limit point in the visual compactification, hence is convergent in the visual compactification. It follows that the identity map $id : X \to X$ has a continuous extension to a map $F : \widehat{X}^{P} \to \overline{X}$ (where $\overline{X}$ is the visual compactification). 

\medskip

Similarly, let $\zeta \in \partial X$, and let $p_n \in X$ be a sequence converging to $\zeta$ in $\overline{X}$. Since $X$ is Gromov hyperbolic, this implies that $(p_n | p_m)_o \to +\infty$ as $m,n \to \infty$. Since $X$ is a Gromov product space, this implies that if $\zeta', \zeta'' \in \partial_{P} X$ are any two limit points of $(p_n)$ in $\widehat{X}^{P}$, then $(\zeta' | \zeta'')_o = +\infty$, and hence $\zeta' = \zeta''$. As above, this implies that the sequence $(p_n)$ converges in $\widehat{X}^{P}$. It follows that $id : X \to X$ extends to a continuous map $G : \widehat{X} \to \widehat{X}^{P}$. 

\medskip

Since $F \circ G_{|X} = id_X$ and $X$ is dense in $\overline{X}$, it follows that $F \circ G = id_{\overline{X}}$. Similarly $G \circ F = id_{\widehat{X}^{P}}$. It follows that $F$ and $G$ are homeomorphisms, and so the compactifications $\widehat{X}^{P}$ and $\overline{X}$ are equivalent. 

\medskip

Since the Gromov inner products $(.|.)_x$ extend continuously to $\widehat{X}^{P} \times \widehat{X}^{P}$ and $\widehat{X}^{P}$ is equivalent to $\overline{X}$, the Gromov inner products extend continuously to $\overline{X} \times \overline{X}$, and so $X$ is boundary continuous. 

\medskip

Since the identification $f : \partial_{P} X \to \partial X$ is the boundary map of the identity map of $X$, for any $x \in X$ the map $f : (\partial_{P} X, \rho_x) \to (\partial X, \rho_x)$ is an isometry of quasi-metric spaces, in particular it is Moebius. 
$\diamond$
      
\medskip

As an immediate corollary of the previous Lemma, we obtain a proof of Theorem \ref{mainthm2} from the Introduction.

\medskip

\noindent{\bf Proof of Theorem \ref{mainthm2}:} Let $Z$ be an antipodal space. 
In Theorem \ref{hypiffqm} it is proved that $\mathcal{M}(Z)$ is Gromov hyperbolic 
if and only if $Z$ is quasi-metric.

\medskip

Suppose that $Z$ is quasi-metric, so that $\mathcal{M}(Z)$ is Gromov hyperbolic. 
By Theorem \ref{gromovip} we know that $\MM(Z)$ is a Gromov product space with Gromov product compactification given by the maxima compactification, and with Gromov product boundary equal to $Z$. It then follows from Lemma \ref{gpbdryghyp} that $\MM(Z)$ is a boundary continuous Gromov hyperbolic space, and the identity map $id : \MM(Z) \to \MM(Z)$ extends to a homeomorphism $F : \overline{\MM(Z)} \to \widehat{\MM(Z)}$. For any $\alpha \in \MM(Z)$, since the identity map of $\MM(Z)$ preserves Gromov inner products with respect to $\alpha$, 
it follows from Theorem \ref{gromovip} that the boundary map of $F$ is an isometry $f : (\partial \MM(Z), \rho_{\alpha}) \to (Z, \alpha)$ of quasi-metric spaces, in particular $\MM(Z)$ is a Moebius hyperbolic filling of $Z$. $\diamond$

\medskip

We recall the definition of maximal Gromov hyperbolic spaces mentioned in the Introduction.

\medskip

\begin{definition} \label{maxspace} {\bf (Maximal Gromov hyperbolic spaces)} Let $X$ be a proper, geodesically complete, boundary continuous
Gromov hyperbolic space. We say that $X$ is {\it maximal} if, given a proper, geodesically complete, boundary continuous
Gromov hyperbolic space $Y$, if there is a Moebius homeomorphism $g : \partial Y \to \partial X$, then $g$ extends
continuously to an isometric embedding $G : Y \to X$.
\end{definition}

\medskip

\begin{lemma} \label{maxgpmaxghyp} Let $X$ be a proper, geodesically complete, boundary continuous Gromov hyperbolic space. Then $X$ is a maximal Gromov product  space if and only if $X$ is a maximal Gromov hyperbolic space.
\end{lemma} 

\medskip

\noindent{\bf Proof:} We note first that $X$ is a Gromov product space whose 
Gromov product space boundary we may identify with its visual boundary, 
$\partial_{P} X = \partial X$. 

\medskip

Suppose that $X$ is a maximal Gromov product space. If $(Y, f : \partial Y \to \partial X)$ is a Moebius hyperbolic filling of $\partial X$, where $Y$ is a proper, geodesically complete, boundary continuous Gromov hyperbolic space and $f$ is a Moebius homeomorphism, then since $Y$ is also a Gromov product space with Gromov product boundary given by $\partial Y$, and $X$ is a maximal Gromov product space, the map $f$ extends to an isometric embedding $F : X \to Y$. Thus $X$ is a maximal Gromov hyperbolic space.

\medskip

Conversely, suppose that $X$ is a maximal Gromov hyperbolic space. Let $(Y, f :\partial_{P} Y \to \partial X)$ be a filling of the antipodal space $\partial X$, where $Y$ is a Gromov product space and $f$ is a Moebius homeomorphism. For any visual quasi-metric $\rho_x$ on $\partial X$, the pull-back $f^* \rho_x \in \MM(\partial_{P} Y)$ is a quasi-metric, so $\partial_{P} Y$ is a quasi-metric antipodal space, and so by Theorem \ref{hypiffqm}, the space $\MM(\partial_{P} Y)$ is Gromov hyperbolic. Since $Y$ embeds isometrically into $\MM(\partial_{P}Y)$ via the visual embedding, it follows that $Y$ is also Gromov hyperbolic. From Lemma \ref{gpbdryghyp}, we then have that $Y$ is a proper, geodesically complete, boundary continuous Gromov hyperbolic space, with Gromov product boundary identified with the visual boundary, $\partial_{P} Y = \partial Y$, and we may identify the map $f$ with a Moebius homeomorphism $f : \partial Y \to \partial X$. Since $X$ is a maximal Gromov hyperbolic space, the map $f$ extends to an isometric embedding $F : Y \to X$. Thus $X$ is a maximal Gromov product space. 
$\diamond$

\medskip

We can now prove Theorem \ref{mainthm3} from the Introduction:

\medskip

%
%
\noindent{\bf Proof of Theorem \ref{mainthm3}:} Let $Z$ be a quasi-metric antipodal space, then we know from Theorems \ref{mainthm1} and \ref{mainthm2} that $\MM(Z)$ is a proper, geodesically complete, boundary continuous Gromov hyperbolic  space. By Theorem \ref{mainthm6} (which we have already proved in section 6), $\MM(Z)$ is a maximal Gromov product space. It follows from Lemma \ref{maxgpmaxghyp} that $X$ is a maximal Gromov hyperbolic space. $\diamond$

\medskip
We now study the visual embedding $i_X : X \to \mathcal{M}(\partial X)$ when $X$ is a proper, geodesically complete, boundary continuous Gromov hyperbolic space.

\medskip

\begin{lemma} \label{buselwrbd} Let $X$ be a proper, geodesically complete, boundary continuous Gromov hyperbolic space. Let $\delta > 0$ be such that all
geodesic triangles in $X$ are $\delta$-thin. Let $x_0 \in X$, and let $\xi, \xi_0 \in \partial X$ be such that
$(\xi | \xi_0)_{x_0} \geq 5\delta$. Let $\gamma_0 : [0, \infty) \to X$ be a geodesic ray such that $\gamma_0(0) = x_0,
\gamma_0(\infty) = \xi_0$. Then for all $t$ such that $0 \leq t \leq 3\delta$, we have
$$
B( x_0, \gamma_0(t), \xi) \geq t - 2\delta,
$$
where $B : X \times X \times \partial X$ is the Busemann cocyle of $X$.
\end{lemma}

\medskip

\noindent{\bf Proof:} Let $\gamma_1 : [0, \infty) \to X$ be a geodesic ray such that
$\gamma_1(0) = x_0, \gamma_1(\infty) = \xi$. For each $r > 0$, let $\alpha_r \subset X$ be a geodesic
segment joining $\gamma_0(r)$ to $\gamma_1(r)$. Since $X$ is boundary continuous there is an $R > 3\delta$
such that $( \gamma_0(r) | \gamma_1(r) )_{x_0} > ( \xi | \xi_0) - \delta$ for all $r \geq R$.

\medskip

We recall that for any $y, z \in X$, we have the well-known elementary estimate $(y | z)_{x_0} \leq d(x_0, \alpha)$
where $\alpha \subset X$ is any geodesic segment joining $y$ to $z$. We fix a $t \in [0, 3\delta]$. Then for $r \geq R$,
this estimate together with the triangle inequality gives
\begin{align*}
d(\gamma_0(t), \alpha_r) & \geq d(x_0, \alpha_r) - t \\
                         & \geq ( \gamma_0(r) | \gamma_1(r) )_{x_0} - t \\
                         & > (5 \delta - \delta) - 3\delta \\
                         & = \delta.
\end{align*}

On the other hand, the geodesic triangle with vertices $x_0, \gamma_0(r), \gamma_1(r)$ given by the geodesics $\gamma_0([0, r]),
\gamma_1([0, r]), \alpha_r \subset X$ is $\delta$-thin, so the point $\gamma_0(t)$ on this triangle is contained in the
union of the $\delta$-neighbourhood of the two sides $\alpha_r, \gamma_1([0, r])$. The above inequality implies then that
$\gamma_0(t)$ must be contained in the $\delta$-neighbourhood of the side $\gamma_1([0, r])$. Thus there exists a
point $z_r \in \gamma_1([0, r])$ such that $d(\gamma_0(t), z_r) \leq \delta$. Applying the triangle inequality we get
\begin{align*}
d(x_0, \gamma_1(r) ) - d( \gamma_0(t), \gamma_1(r) ) & = d(x_0, z_r ) + d( z_r, \gamma_1(r) ) - d( \gamma_0(t), \gamma_1(r) ) \\
                                                     & \geq d(x_0, z_r) - d( z_r, \gamma_0(t)) \\
                                                     & \geq d(x_0, \gamma_0(t)) - 2 d( z_r, \gamma_0(t) ) \\
                                                     & \geq t - 2\delta.
\end{align*}
Now letting $r$ tend to $\infty$ in the above inequality and using the definition of the Busemann cocycle we get
$$
B(x_0, \gamma_0(t), \xi) \geq t - 2\delta
$$
as required. $\diamond$

\medskip

As a consequence of the above Lemma, we can prove the following:

\medskip

\begin{theorem} \label{embeddeldense} Let $X$ be a proper, geodesically complete, boundary continuous Gromov hyperbolic
space. Let $\delta > 0$ be such that all geodesic triangles in $X$ are $\delta$-thin. Then the image of the
isometric embedding $i_X : X \to \mathcal{M}(\partial X)$ is $7\delta$-dense in $\mathcal{M}(\partial X)$.
\end{theorem}

\medskip

\noindent{\bf Proof:} Let $\rho \in \mathcal{M}(\partial X)$. Then we need to show that there is an $x \in X$ such that
$d_{\mathcal{M}}( \rho, \rho_x) \leq 7\delta$. Define $\phi : X \to [0,\infty)$ by $\phi(x) = d_{\mathcal{M}}( \rho, \rho_x), x \in X$.
Since $i_X$ is an isometric embedding, it is clear that the function $\phi$ is proper, i.e. $\phi(x) \to +\infty$ as $d(x, y_0) \to +\infty$
for any fixed $y_0 \in X$. Since the space $X$ is proper, it follows that the function $\phi$ attains a global
minimum at some $x_0 \in X$. Let $r = \phi(x_0) = d_{\mathcal{M}}( \rho, \rho_{x_0})$.

\medskip

Let $\xi_0 \in \hbox{argmax} \ \frac{d\rho}{d\rho_{x_0}}$, and let $U = \{ \xi \in \partial X | ( \xi | \xi_0)_{x_0} \geq 5\delta \}$.
Let $\gamma_0 : [0, \infty)$ be a geodesic ray such that $\gamma_0(0) = x_0, \gamma_0(\infty) = \xi_0$. We let $x = \gamma_0(3\delta)$.
Then we claim that there exists $\xi \in \partial X - U$ such that $\log \frac{d\rho}{d\rho_{x_0}}(\xi) \geq r - 4\delta$. If not,
then for all $\xi \in \partial X - U$ we have, since $d_{\mathcal{M}}( \rho_{x_0}, \rho_x) = d(x_0, x) = 3\delta$, that
\begin{align*}
\log \frac{d\rho}{d\rho_x}(\xi) & = \log \frac{d\rho}{d\rho_{x_0}}(\xi) - \log \frac{d\rho_x}{d\rho_{x_0}}(\xi) \\
                                & \leq (r - 4\delta) + 3\delta \\
                                & = r - \delta,
\end{align*}
while for all $\xi \in U$ we have, by Lemma \ref{buselwrbd} and Lemma \ref{visualderiv}, that
\begin{align*}
\log \frac{d\rho}{d\rho_x}(\xi) & = \log \frac{d\rho}{d\rho_{x_0}}(\xi) - \log \frac{d\rho_x}{d\rho_{x_0}}(\xi) \\
                                & \leq r - B( x_0, x, \xi) \\
                                & \leq r - (3\delta - 2\delta) \\
                                & = r - \delta.
\end{align*}
The above two inequalities then imply that
\begin{align*}
\phi(x) = d_{\mathcal{M}}( \rho, \rho_x) & = \max_{\xi \in \partial X} \log \frac{d\rho}{d\rho_x}(\xi) \\
                                         & \leq r - \delta \\
                                         & < r = \phi(x_0),
\end{align*}
which contradicts the fact that $x_0$ is a global minimum of the function $\phi$. Thus our claim holds,
and so there is a $\xi \in \partial X - U$ such that $\log \frac{d\rho}{d\rho_{x_0}}(\xi) \geq r - 4\delta$.
Since $\xi \in \partial X - U$, we have $\rho_{x_0}(\xi_0, \xi) = \exp( - (\xi|\xi_0)_{x_0} ) > e^{-5\delta}$.
Applying the Geometric Mean-Value Theorem then gives
\begin{align*}
1 & \geq \rho(\xi_0, \xi)^2 \\
  & = \frac{d\rho}{d\rho_{x_0}}(\xi_0) \cdot \frac{d\rho}{d\rho_{x_0}}(\xi) \cdot \rho_{x_0}(\xi_0, \xi)^2 \\
  & \geq e^r \cdot e^{r - 4\delta} \cdot e^{-10\delta} \\
  & = e^{2r - 14\delta},
\end{align*}
thus
$$
d_{\mathcal{M}}( \rho, \rho_{x_0}) = r \leq 7\delta
$$
as required. $\diamond$

\medskip

We can now prove Theorem \ref{mainthm4} from the Introduction:

\medskip

\noindent{\bf Proof of Theorem \ref{mainthm4}:} Let $X$ be a proper, 
geodesically complete, boundary continuous Gromov hyperbolic space. 

\medskip

The first statement of Theorem \ref{mainthm4} is proved in Theorem 
\ref{embeddeldense} above. 

\medskip

For the second statement, we note that by Theorem \ref{mainthm7} (which was already proved in section 6) $X$ is a maximal Gromov product space if and only if the visual embedding $i_X : X \to \MM(\partial X)$ is an isometry, so by Lemma \ref{maxgpmaxghyp} it follows that $X$ is a maximal Gromov hyperbolic space if and only if the visual embedding is an isometry, as required. $\diamond$

\medskip

We remark that any tree is boundary continuous, since it is CAT(-1).

\medskip

\begin{cor} \label{treesurj} Let $X$ be a proper, geodesically complete tree. Then the visual embedding $i_X : X \to \mathcal{M}(\partial X)$ is an isometry.
\end{cor}

\medskip

\noindent{\bf Proof:} Let $\rho \in \mathcal{M}(\partial X)$. Let $\delta_n > 0$ be a sequence such that $\delta_n \to 0$ as $n \to \infty$.
Since $X$ is a tree, geodesic triangles in $X$ are $\delta_n$-thin for all $n$, hence by Theorem \ref{embeddeldense} for all $n$
there exists $x_n \in X$ such that $d_{\mathcal{M}}(\rho, \rho_{x_n}) \leq 7\delta_n$. Since $i_X$ is an isometric embedding, this implies
in particular that the sequence $x_n \in X$ is bounded, hence after passing to a subsequence we may assume that $x_n \to x$ as $n \to \infty$
for some $x \in X$. Since $\delta_n \to 0$ and $d_{\mathcal{M}}(\rho, \rho_{x_n}) \leq 7\delta_n$ for all $n$, letting $n$ tend to $\infty$ in
this inequality gives $d_{\mathcal{M}}( \rho, \rho_x) = 0$, hence $\rho = \rho_x$ and $i_X$ is surjective. $\diamond$

\medskip

We can now use one of the main results of \cite{beyrer-schroeder} to prove the following:

\medskip

\begin{prop} \label{mztree} Let $Z$ be an antipodal space, and suppose that some $\rho_0 \in \mathcal{M}(Z)$ is
an ultrametric. Then $\mathcal{M}(Z)$ is a tree, and all elements $\rho$ of $\mathcal{M}(Z)$ are ultrametrics.
\end{prop}

\medskip

\noindent{\bf Proof:} Let $\rho_0 \in \mathcal{M}(Z)$ be an ultrametric. Since $Z$ is compact, the ultrametric
space $(Z, \rho_0)$ is complete, and hence, as remarked in \cite{beyrer-schroeder}, defines, in the
terminology of \cite{beyrer-schroeder}, an {\it ultrametric Moebius space}. Then by Theorem 4.1 of
\cite{beyrer-schroeder}, there exists a geodesically complete tree $X$ such that $\partial X$ is Moebius
homeomorphic to $Z$. Moreover, it is easy to see from the explicit definition of the tree $X$ in terms of $Z$ given in
\cite{beyrer-schroeder}, that $X$ is proper because $Z$ is compact. Note that $X$ is boundary continuous since
it is CAT(-1).

\medskip

Let $g : \partial X \to Z$ be a Moebius homeomorphism. By Theorem \ref{mbisom1}, the associated map
$G = g_* : \mathcal{M}(\partial X) \to \mathcal{M}(Z)$ is a surjective isometry, while by Corollary \ref{treesurj}
the map $i_X : X \to \mathcal{M}(\partial X)$ is a surjective isometry, thus the composition
$H = G \circ i_X : X \to \mathcal{M}(Z)$ is a surjective isometry, and so $\mathcal{M}(Z)$ is a tree.
It follows from Theorem \ref{hypiffqm} that all $\rho \in \mathcal{M}(Z)$ are ultrametrics. $\diamond$

\medskip

Finally, we prove Theorem \ref{mainthm5} from the Introduction:

\medskip

\noindent{\bf Proof of Theorem \ref{mainthm5}:} Let $X, Y$ be maximal Gromov hyperbolic spaces, and let $f : \partial X \to \partial Y$ be a Moebius homeomorphism. It follows from the proof of Theorem \ref{mainthm4} that $X$ and $Y$ are also maximal Gromov product spaces, and by Lemma \ref{gpbdryghyp} the Gromov product boundaries of $X$ and $Y$ can be identified with the visual boundaries $\partial X$ and $\partial Y$. It then follows from Theorem \ref{mainthm8} (which was already proved in section 6) that $f$ extends to an isometry $F : X \to Y$. $\diamond$

\medskip

\section{Tangent spaces of $\mathcal{M}(Z)$.}

\medskip

Let $Z$ be an antipodal space. In this section we study the infinitesimal structure of the space $\mathcal{M}(Z)$.
Recall that any $\rho_1 \in \mathcal{M}(Z)$ gives rise to an isometric embedding $\phi_{\rho_1} : \mathcal{M}(Z) \to C(Z), \rho \mapsto \log \frac{d\rho}{d\rho_1}$, of $\mathcal{M}(Z)$ into the Banach space $C(Z)$. For any other $\rho_2 \in \mathcal{M}(Z)$, the Chain Rule implies that
$$
\phi_{\rho_1}(\rho) = \phi_{\rho_2}(\rho) + \log \frac{d\rho_2}{d\rho_1}
$$
for all $\rho \in \mathcal{M}(Z)$, in other words the isometric embeddings $\phi_{\rho_1}$ and $\phi_{\rho_2}$ differ by a translation
of the Banach space $C(Z)$. It follows that if we have a curve $\gamma : t \in (a, b) \mapsto \rho(t) \in \mathcal{M}(Z)$ in $\mathcal{M}(Z)$,
then for the curves
$\phi_{\rho_1} \circ \gamma, \phi_{\rho_2} \circ \gamma$ in the Banach space $C(Z)$, one is differentiable if and only if the other is,
and moreover in this case
$$
\frac{d}{dt} \phi_{\rho_1} \circ \gamma = \frac{d}{dt} \phi_{\rho_2} \circ \gamma.
$$
We can therefore talk unambiguously of differentiable curves in $\mathcal{M}(Z)$, namely we can define a curve $\gamma$ in $\mathcal{M}(Z)$
to be differentiable if the curve $\phi_{\rho_1} \circ \gamma$ in $C(Z)$ is differentiable for some (and hence for all) $\rho_1 \in \mathcal{M}(Z)$.
Similarly one can talk of differentiable maps $\psi : M \to \mathcal{M}(Z)$ from a manifold $M$ into the space $\mathcal{M}(Z)$. These
considerations lead naturally to the following definition of tangent spaces at points of $\mathcal{M}(Z)$:

\medskip

\begin{definition} \label{tangentdefn} {\bf (Tangent spaces of $\mathcal{M}(Z)$)} Let $\rho \in \mathcal{M}(Z)$. We say that a curve
$t \in (-\epsilon, \epsilon) \mapsto \rho(t) \in \mathcal{M}(Z)$ is admissible at $\rho$ if $\rho(0) = \rho$ and
$\frac{d}{dt}_{|t = 0} \phi_{\rho}(\rho(t)) \in C(Z)$ exists.
Then the tangent space
to $\mathcal{M}(Z)$ at $\rho$ is defined to be the subset $T_{\rho} \mathcal{M}(Z)$ of $C(Z)$ given by
$$
T_{\rho} \mathcal{M}(Z) := \{ \ v = \frac{d}{dt}_{|t = 0} \phi_{\rho}(\rho(t)) \in C(Z) \ | \ t \mapsto \rho(t) \ \hbox{is admissible at} \ \rho \ \}.
$$
\end{definition}

\medskip

We remark that a priori it is not even clear that the tangent space $T_{\rho} \mathcal{M}(Z)$ is a vector space.
We will show that this is indeed the case, and identify explicitly the vector space $T_{\rho} \mathcal{M}(Z)$.
We need the following definition:

\medskip

Let $\rho \in \mathcal{M}(Z)$. Then a continuous function $v \in C(Z)$ is said to be {\it $\rho$-odd} if, for all $\xi, \eta \in Z$,
whenever $\rho(\xi, \eta) = 1$ then $v(\xi) + v(\eta) = 0$. We will denote by $\mathcal{O}_{\rho}(Z) \subset C(Z)$ the set of $\rho$-odd
functions, then it is clear that $\mathcal{O}_{\rho}(Z)$ is a closed linear subspace of $C(Z)$.

\medskip

Recall that for $\rho \in \mathcal{M}(Z)$ and $v \in C(Z)$ the discrepancy function $D_{\rho}(v) \in C(Z)$ was defined by
$$
D_{\rho}(v)(\xi) = \max_{\eta \in Z - \{\xi\}} v(\xi) + v(\eta) + \log \rho(\xi, \eta)^2, \ \xi \in Z.
$$
Note that if $v$ is $\rho$-odd, then given $\xi \in Z$ we can choose $\eta \in Z$ such that $\rho(\xi, \eta) = 1$, and so
\begin{align*}
D_{\rho}(v)(\xi) & \geq v(\xi) + v(\eta) + \log \rho(\xi, \eta)^2 \\
                 & = v(\xi) + v(\eta) \\
                 & = 0,
\end{align*}
thus $D_{\rho}(v)(\xi) \geq 0$ for all $\xi \in Z$ if $v$ is $\rho$-odd.

\medskip

\begin{lemma} \label{discrepodd} Let $\rho \in \mathcal{M}(Z)$ and let $v \in \mathcal{O}_{\rho}(Z)$ be a
$\rho$-odd function. Then
$$
\lim_{t \to 0} \frac{ || D_{\rho}( tv ) ||_{\infty} }{|t|} = 0.
$$
\end{lemma}

\medskip

\noindent{\bf Proof:} Suppose not. Then there exists $\kappa > 0$ and a sequence $t_n \neq 0$ such that $t_n \to 0$ as $n \to \infty$,
and
$$
\frac{ || D_{\rho}( t_n v ) ||_{\infty} }{|t_n|} \geq \kappa
$$
for all $n$. Let $\xi_n \in Z$ be such that $|| D_{\rho}( t_n v ) ||_{\infty} = | D_{\rho}( t_n v )(\xi_n) |$. As noted above,
since $v$ is $\rho$-odd, $D_{\rho}( t_n v )(\xi_n) \geq 0$, and so $|| D_{\rho}( t_n v ) ||_{\infty} = D_{\rho}( t_n v )(\xi_n)$.
Now let $\eta_n \in Z$ be such that
$$
D_{\rho}( t_n v)(\xi_n) = t_n v(\xi_n) + t_n v(\eta_n) + \log \rho(\xi_n, \eta_n)^2.
$$
Passing to a subsequence, we may assume there are $\xi, \eta \in Z$ such that $\xi_n \to \xi, \eta_n \to \eta$ as $n \to \infty$.
Now by our hypothesis, for all $n$ we have
\begin{align*}
|v(\xi_n) + v(\eta_n)| & = \frac{| t_n v(\xi_n) + t_n v(\eta_n) |}{|t_n|} \\
                       & \geq \frac{t_n v(\xi_n) + t_n v(\eta_n) + \log \rho(\xi_n, \eta_n)^2}{|t_n|} \\
                       & = \frac{ || D_{\rho}( t_n v ) ||_{\infty} }{|t_n|} \\
                       & \geq \kappa.
\end{align*}
Letting $n$ tend to $\infty$ in the above inequality gives $|v(\xi)+v(\eta)| \geq \kappa > 0$. Since $v$ is $\rho$-odd this
implies that $\rho(\xi, \eta) < 1$. Since $\rho(\xi_n, \eta_n) \to \rho(\xi, \eta)$ as $n \to \infty$, we can then
choose $\epsilon > 0$ and $N \geq 1$ such that for all $n \geq 1$ we have $\log \rho(\xi_n, \eta_n)^2 \leq -\epsilon$.
We may also assume that $N \geq 1$ is chosen such that $2|t_n| ||v||_{\infty} < \epsilon$ for all $n \geq N$. Then for
any $n \geq N$ we have
\begin{align*}
|| D_{\rho} (t_n v) ||_{\infty} & = D_{\rho}( t_n v)(\xi_n) \\
                                & = t_n v(\xi_n) + t_n v(\eta_n) + \log \rho(\xi_n, \eta_n)^2 \\
                                & \leq 2|t_n| ||v||_{\infty} - \epsilon \\
                                & < 0,
\end{align*}
a contradiction. The Lemma follows. $\diamond$

\medskip

We can now prove Theorem \ref{mainthm11} from the Introduction, which we restate here:

\medskip

\begin{theorem} \label{tangentspace} Let $\rho \in \mathcal{M}(Z)$. Then the tangent space to $\mathcal{M}(Z)$ at $\rho$ is equal
to the space of $\rho$-odd functions, i.e.
$$
T_{\rho} \mathcal{M}(Z) = \mathcal{O}_{\rho}(Z).
$$
In particular, $T_{\rho} \mathcal{M}(Z)$ is a closed linear subspace of $C(Z)$.
\end{theorem}

\medskip

\noindent{\bf Proof:} Let $v \in T_{\rho} \mathcal{M}(Z)$. Then $v = \frac{d}{dt}_{|t = 0} \phi_{\rho}(\rho(t))$, where
$t \in (-\epsilon, \epsilon) \mapsto \rho(t) \in \mathcal{M}(Z)$ is a curve such that $\rho(0) = \rho$ and the derivative
$\frac{d}{dt}_{|t = 0} \phi_{\rho}(\rho(t))$ exists in $C(Z)$. This means that
$$
\lim_{t \to 0} \ \left|\left| \frac{1}{t} \log \frac{d\rho(t)}{d\rho} - v \right|\right|_{\infty} = 0.
$$
In particular, for a fixed $\xi \in Z$, the function $u_{\xi}(t) =  \log \frac{d\rho(t)}{d\rho}(\xi)$ is differentiable at $t = 0$,
and $u_{\xi}'(0) = v(\xi)$. Now let $\xi, \eta \in Z$ be such that $\rho(\xi, \eta) = 1$. Then by the Geometric Mean-Value Theorem
\begin{align*}
 \log \frac{d\rho(t)}{d\rho}(\xi) + \log \frac{d\rho(t)}{d\rho}(\eta) & = \log \rho(t)(\xi, \eta)^2 \\
                                                                      & \leq 0
\end{align*}
for all $t \in (-\epsilon, \epsilon)$, with equality for $t = 0$. Thus the function $u_{\xi}(t) + u_{\eta}(t)$ has a local
maximum at $t = 0$, and since it is differentiable at $t = 0$ this implies
$$
v(\xi) + v(\eta) = u_{\xi}'(0) + u_{\eta}'(0) = 0,
$$
and so it follows that $v$ is $\rho$-odd. Thus $T_{\rho} \mathcal{M}(Z) \subset O_{\rho}(Z)$.

\medskip

To show the reverse inclusion, let $v \in O_{\rho}(Z)$ be a given $\rho$-odd function. We define
$$
\rho(t) = \mathcal{P}_{\infty}( tv, \rho) \in \mathcal{M}(Z), \ t \in \mathbb{R}.
$$
Note that then $\rho(0) = \rho$.

\medskip

Let $\epsilon > 0$ be given. Then by Lemma \ref{discrepodd}, there exists $t_0 > 0$ such that
$$
\frac{ || D_{\rho}( tv ) ||_{\infty} }{|t|} < \epsilon/2
$$
for all $t$ such that $0 < |t| < t_0$. Now fix a $t$ such that $0 < |t| < t_0$.
Let $\xi \in Z$, and let $\eta \in Z$ be such that
$$
D_{\rho}( tv )(\xi) = tv(\xi) + tv(\eta) + \log \rho(\xi, \eta)^2.
$$
Now using Proposition \ref{pinfest}, we obtain the following estimates for $\log \frac{d\rho(t)}{d\rho}$,
\begin{align*}
\log \frac{d\rho(t)}{d\rho}(\xi) & \geq t v(\xi) - \frac{1}{2} D_{\rho}( tv )(\xi) \\
                                 & \geq t v(\xi) - |t| \frac{\epsilon}{4},
\end{align*}
and
\begin{align*}
\log \frac{d\rho(t)}{d\rho}(\xi) & \leq t v(\xi) - D_{\rho}( tv )(\xi) + \frac{1}{2} D_{\rho} ( tv )(\eta) \\
                                 & \leq t v(\xi) + |t| \frac{3\epsilon}{4}.
\end{align*}
Dividing the above inequalities by $t$, we get, if $t > 0$,
$$
-\frac{\epsilon}{4} \leq \frac{1}{t} \log \frac{d\rho(t)}{d\rho}(\xi) - v(\xi) \leq \frac{3\epsilon}{4}
$$
for all $\xi \in Z$, while if $t < 0$ then we get
$$
-\frac{3\epsilon}{4} \leq \frac{1}{t} \log \frac{d\rho(t)}{d\rho}(\xi) - v(\xi) \leq \frac{\epsilon}{4}
$$
for all $\xi \in Z$.

\medskip

Either way, we have
$$
\left|\left| \frac{1}{t} \log \frac{d\rho(t)}{d\rho} - v \right|\right|_{\infty} < \epsilon
$$
for all $t$ such that $0 < |t| < t_0$. Thus
$$
\lim_{t \to 0} \left|\left| \frac{1}{t} \log \frac{d\rho(t)}{d\rho} - v \right|\right|_{\infty} = 0,
$$
hence the derivative $\frac{d}{dt}_{|t = 0} \phi_{\rho}(\rho(t))$ exists in $C(Z)$
and equals $v$, so by definition $v \in T_{\rho} \mathcal{M}(Z)$. Thus $T_{\rho} \mathcal{M}(Z) \supset \mathcal{O}_{\rho}(Z)$,
and so $T_{\rho} \mathcal{M}(Z) = \mathcal{O}_{\rho}(Z)$. $\diamond$

\medskip

We apply the above Theorem to some examples:

\medskip

\noindent{\bf Example 1. Trees.} Let $Z$ be an ultrametric antipodal space, so that $\mathcal{M}(Z)$ is a tree by Proposition \ref{mztree}.
By Theorem \ref{mainthm2}, we can identify the visual boundary $\partial \mathcal{M}(Z)$ of the tree $\mathcal{M}(Z)$ with the space $Z$. For any pair of distinct points $\xi, \eta \in Z$ there is a unique bi-infinite geodesic in $\mathcal{M}(Z)$
with endpoints at infinity $\xi, \eta$, we will denote this bi-infinite geodesic by $(\xi, \eta) \subset \mathcal{M}(Z)$. Similarly
for $\rho \in \mathcal{M}(Z)$ and $\xi \in Z$ we will denote by $[\rho, \xi) \subset \mathcal{M}(Z)$ the
unique geodesic ray in $\mathcal{M}(Z)$ starting from $\rho$ with endpoint at infinity $\xi$.

\medskip

We say that $\rho \in \mathcal{M}(Z)$ is a {\it vertex} of this tree if there exist three distinct points $\xi_1, \xi_2, \xi_3 \in Z$ such that
$\rho \in (\xi_1, \xi_2) \cap (\xi_2, \xi_3) \cap (\xi_3, \xi_1)$. In this case $\rho(\xi_i, \xi_j) = 1$ for all $i \neq j$. Thus if
$v \in \mathcal{O}_{\rho}(Z)$ is a $\rho$-odd function on $Z$, then $v(\xi_i) + v(\xi_j) = 0$ for all $i \neq j$, from which it
follows that $v(\xi_i) = 0$ for all $i$. Now for any $\xi \in Z$, if $\rho(\xi, \xi_i) = 1$ for all $i$ then $v(\xi) + v(\xi_i) = 0$ for all
$i$ and so $v(\xi) = - v(\xi_i) = 0$. Otherwise there exists $i$ such that $\rho(\xi, \xi_i) < 1$, which means that the geodesic rays
$[\rho, \xi)$ and $[\rho, \xi_i)$ intersect in a geodesic segment of positive length. It then follows easily that $\rho(\xi, \xi_j) = 1$ for all
$j \neq i$, hence $v(\xi) = - v(\xi_j) = 0$. It follows that any $\rho$-odd function is identically zero if $\rho$ is a vertex of the tree,
and so the tangent space $T_{\rho} \mathcal{M}(Z)$ is trivial in this case,
$$
T_{\rho} \mathcal{M}(Z) = \mathcal{O}_{\rho}(Z) = \{ 0 \}
$$
when $\rho$ is a vertex.

\medskip

On the other hand, if $\rho \in \mathcal{M}(Z)$ is not a vertex, then it is not hard to show that $\mathcal{M}(Z) - \{\rho\}$ has exactly
two connected components, say $A, B$. If we denote by $\partial A, \partial B \subset Z$ the boundaries of $A, B$ in $\overline{\mathcal{M}(Z)}$ (minus the point $\rho$),
then we have $Z = A \sqcup B$, and $\rho(\xi, \eta) = 1$ for all $\xi \in A, \eta \in B$. Fix a $\xi_0 \in A$ and an $\eta_0 \in B$.
Let $v \in \mathcal{O}_{\rho}(Z)$ be a $\rho$-odd function. Then for all $\eta \in B$, we have $v(\eta) = - v(\xi_0)$, while for all
$\xi \in A$, we have $v(\xi) = -v(\eta_0) = v(\xi_0)$. Thus the $\rho$-odd functions are precisely the functions $v$ on $Z$ such that
$v \equiv c$ on $A$ and $v \equiv -c$ on $B$ for some $c \in \mathbb{R}$. It follows that the linear map
\begin{align*}
T_{\rho} \mathcal{M}(Z) & \to \mathbb{R} \\
           v            & \mapsto v(\xi_0)
\end{align*}
is an isomorphism, and so the tangent space $T_{\rho} \mathcal{M}(Z)$ is one-dimensional when $\rho$ is not a vertex, as expected.

\medskip

\noindent{\bf Example 2. Uniquely antipodal elements of $\mathcal{M}(Z)$.} For an antipodal space $Z$, we say that an element
$\rho \in \mathcal{M}(Z)$ is {\it uniquely antipodal} if for every $\xi \in Z$ there exists a {\bf unique} $\eta \in Z$
such that $\rho(\xi, \eta) = 1$. A well-known example of such uniquely antipodal elements is given by the visual
metrics $\rho = \rho_x$ on the boundary $Z = \partial X$ of a complete, simply connected Riemannian manifold $X$ of sectional curvature bounded
above by $-1$. For a uniquely antipodal $\rho$ on an antipodal space $Z$, by definition we obtain an antipodal map $i_{\rho} : Z \to Z$
such that $\rho(\xi, \eta) = 1$ if and only if $\eta = i_{\rho}(\xi)$. It is not hard to see from compactness of $Z$ and uniqueness
of antipodal points that the antipodal map $i_{\rho} : Z \to Z$ must in fact be continuous, and is hence a homeomorphism since
$i_{\rho} \circ i_{\rho} = id$.

\medskip

We can thus define in this case a linear projection $O_{\rho} : C(Z) \to \mathcal{O}_{\rho}(Z)$ from $C(Z)$ onto the closed linear subspace
of $\rho$-odd functions by the usual formula
$$
O_{\rho}(u) := \frac{1}{2} \left( u - u \circ i_{\rho} \right) \ , \ u \in C(Z).
$$
It is then easy to see that if $Z$ is infinite,
then the space of $\rho$-odd functions is infinite dimensional. For this, one simply chooses a small neighbourhood $U$ of a non-isolated point
$\xi_0 \in Z$ such that the closures of $U$ and $i_{\rho}(U)$ are disjoint. Then for any $u \in C_c(U)$ (the space of continuous functions on
$Z$ with support contained in $U$), we can define a $\rho$-odd function $v$ by $v = O_{\rho}(u)$. Noting that $v_{|U} = u_{|U}$, this
implies that the restriction $\left(O_{\rho}\right)_{|C_c(U)} : C_c(U) \to \mathcal{O}_{\rho}(Z)$ is injective, hence $\mathcal{O}_{\rho}(Z)$
is infinite-dimensional.

\medskip

Thus the tangent space $T_{\rho} \mathcal{M}(Z)$ is infinite-dimensional whenever $\rho$ is uniquely antipodal (in fact one only needs
$\rho$ to be uniquely antipodal in a small enough neighbourhood of some non-isolated point in $Z$ and then the same argument goes through).
In particular, taking $\rho = \rho_x$ on $\partial X$ where $X$ is a complete, simply connected Riemannian manifold of sectional
curvature bounded above by $-1$, using the fact that the tangent spaces of $X$ are finite dimensional while the tangent space to $\MM(\partial X)$ at $\rho_x$ is infinite-dimensional, it is not hard to show that in this case the visual embedding $i_X : X \to \mathcal{M}(\partial X)$ is {\bf not}
surjective, and in a sense the image of $X$ in $\MM(\partial X)$ has infinite 
codimension.

\medskip

\bibliography{moeb}
\bibliographystyle{alpha}

\end{document}